\documentstyle[12pt]{amsart}

\newcommand{\nc}{\newcommand}

\nc{\ad}{\operatorname{ad}}
\nc{\bSt}{\mbox{\bf{St}}}
\nc{\card}{\operatorname{card}}
\nc{\cd}{\operatorname{cd}}
\nc{\Ch}{\operatorname{Ch}}
\nc{\chara}{\operatorname{char}}
\nc{\CHom}{\cal{H}om}
\nc{\codim}{\operatorname{codim}}
\nc{\Cone}{\operatorname{Cone}}
\nc{\depth}{\operatorname{depth}}
\nc{\dirlim}{\underset{\rightarrow}{\operatorname{lim}}}
\nc{\emp}{\emptyset}
\nc{\Fac}{\cal{F}ac}
\nc{\Hom}{\operatorname{Hom}}
\nc{\Id}{\operatorname{Id}}
\nc{\Ima}{\operatorname{Im}}
\nc{\Ind}{\operatorname{Ind}}
\nc{\invlim}{\underset{\leftarrow}{\operatorname{lim}}}
\nc{\Ker}{\operatorname{Ker}}
\nc{\Ob}{\operatorname{Ob}}
\nc{\one}{\mbox{\bf{1}}}
\nc{\Or}{\cal{O}r}
\nc{\Ord}{\cal{O}rd}
\nc{\Part}{\cal{P}art}
\nc{\sgn}{\operatorname{sgn}}
\nc{\Sh}{\cal{S}h}
\nc{\Tor}{\operatorname{Tor}}
\nc{\Vect}{\cal{V}ect}

\nc{\BA}{\Bbb A}
\nc{\ba}{\mbox{\bf{a}}}
\nc{\baJ}{\bar{J}}
\nc{\BAO}{\overset{\circ}{\Bbb A}}
\nc{\BB}{\Bbb B}
\nc{\BC}{\Bbb C}
\nc{\BE}{\overline{E}}
\nc{\BF}{\overline{F}}
\nc{\bF}{\mbox{\bf{F}}}
\nc{\bL}{\mbox{\bf{L}}}
\nc{\bM}{\mbox{\bf{M}}}
\nc{\bmu}{\vec{\mu}}
\nc{\BN}{\Bbb N}
\nc{\bnu}{\vec{\nu}}
\nc{\BP}{\Bbb P}
\nc{\bP}{\mbox{\bf{P}}}
\nc{\bq}{\mbox{\bf{q}}}
\nc{\BR}{\Bbb R}
\nc{\br}{\mbox{\bf{r}}}
\nc{\bs}{\mbox{\bf{s}}}
\nc{\bt}{\mbox{\bf{t}}}
\nc{\bU}{\mbox{\bf{U}}}
\nc{\BUpsilon}{\bar{\Upsilon}}
\nc{\bw}{\mbox{\bf{w}}}
\nc{\bx}{\mbox{\bf{x}}}
\nc{\BZ}{\Bbb Z}
\nc{\bz}{\mbox{\bf{z}}}

\nc{\CA}{\cal A}
\nc{\CAO}{\overset{\circ}{\cal{A}}}
\nc{\CB}{\cal B}
\nc{\CC}{\cal C}
\nc{\CD}{\cal D}
\nc{\CE}{\cal E}
\nc{\CF}{\cal F}
\nc{\CH}{\cal H}
\nc{\CI}{\cal I}
\nc{\CJ}{\cal J}
\nc{\CK}{\cal K}
\nc{\CL}{\cal L}
\nc{\CM}{\cal M}
\nc{\CN}{\cal N}
\nc{\CP}{\cal P}
\nc{\CQ}{\cal Q}
\nc{\CS}{\cal S}
\nc{\CX}{\cal X}

\nc{\DO}{\overset{\circ}{D}}
\nc{\dpar}{\partial}

\nc{\fE}{\frak{E}}
\nc{\fF}{\frak F}
\nc{\ff}{\frak f}
\nc{\fu}{\frak{u}}

\nc{\HO}{\overset{\circ}{H}}

\nc{\tC}{\tilde{C}}
\nc{\tc}{\tilde{c}}
\nc{\tCC}{\tilde{\cal{C}}}
\nc{\tD}{\tilde{D}}
\nc{\tE}{\tilde E}
\nc{\tF}{\tilde F}
\nc{\tfF}{\tilde{\frak{F}}}
\nc{\tJ}{\tilde{J}}
\nc{\tK}{\tilde K}
\nc{\tM}{\tilde{M}}
\nc{\tP}{\tilde{P}}
\nc{\tS}{\tilde S}
\nc{\ttheta}{\tilde{\theta}}
\nc{\tU}{\tilde{U}}
\nc{\tUpsilon}{\tilde{\Upsilon}}
\nc{\ty}{\tilde y}
\nc{\tY}{\tilde Y}
\nc{\txi}{\tilde{\xi}}

\nc{\nen}{\newenvironment}
\nc{\ol}{\overline}
\nc{\ul}{\underline}
\nc{\ra}{\rightarrow}
\nc{\lra}{\longrightarrow}
\nc{\Lra}{\Longrightarrow}
\nc{\Lla}{\Longleftarrow}
\nc{\Llra}{\Longleftrightarrow}
\nc{\hra}{\hookrightarrow}
\nc{\iso}{\overset{\sim}{\lra}}

\nc{\Thm}[1]{Theorem~\ref{#1}}
\nc{\Prop}[1]{Proposition~\ref{#1}}
\nc{\Lem}[1]{Lemma~\ref{#1}}
\nc{\Cor}[1]{Corollary~\ref{#1}}
\nc{\Conj}[1]{Conjecture~\ref{#1}}
\nc{\Claim}[1]{Claim~\ref{#1}}
\nc{\Defn}[1]{Definition~\ref{#1}}
\nc{\Exa}[1]{Example~\ref{#1}}
\nc{\Rem}[1]{Remark~\ref{#1}}
\nc{\Note}[1]{Note~\ref{#1}}

\nen{thm}[1]{\label{#1}{\bf Theorem.\ } \em}{}
\nen{prop}[1]{\label{#1}{\bf Proposition.\ } \em}{}
\nen{lem}[1]{\label{#1}{\bf Lemma.\ } \em}{}
\nen{cor}[1]{\label{#1}{\bf Corollary.\ } \em}{}
\nen{conj}[1]{\label{#1}{\bf Conjecture.\ } \em}{}
\nen{claim}[1]{\label{#1}{\bf Claim.\ } \em}{}

\nen{defn}[1]{\label{#1}{\bf Definition.\ } }{}
\nen{exa}[1]{\label{#1}{\bf Example.\ } }{}

\nen{rem}[1]{\label{#1}{\em Remark.\ } }{}
\nen{note}[1]{\label{#1}{\em Note.\ } }{}
\nen{exer}[1]{\label{#1}{\em Exercise.\ } }{}

\begin{document}

\title[]{Localization of $\fu$-modules. II.\\ Configuration spaces and
quantum groups}
\author{Michael Finkelberg}
\address{Independent Moscow University, 65-3 Mikloukho-Maklai St.,
apt. 86, Moscow 117342 Russia}
\author{Vadim Schechtman}
\address{Dept. of Mathematics, SUNY at Stony Brook, Stony Brook,
NY 11794-3651 USA}
\thanks{The second author was supported in part by NSF grant DMS-9202280}
\date{December 31, 1994; revised: January 1995\\
q-alg/9412017}
\maketitle

\section{Introduction}

\subsection{} This paper is a sequel to ~\cite{fs}.
We are starting here the geometric study of the tensor category $\CC$
associated with a quantum group (corresponding to a
Cartan matrix of finite type) at a root of unity (see ~\cite{ajs}, 1.3 and
the present paper, ~\ref{C} for the precise definitions).

The main results of this paper are Theorems ~\ref{shaposym},
{}~\ref{stalks}, ~\ref{shsymnsl} and ~\ref{stalknsl} which

---  establish isomorphisms between homogeneous
components of irreducible objects in $\CC$ and spaces of vanishing
cycles at the origin of certain Goresky-MacPherson sheaves on
configuration spaces;

--- establish isomorphisms of the stalks at the origin
of the above GM sheaves with certain Hochschild complexes
(which compute the Hochschild homology of a certain "triangular"
subalgebra of our quantum group with coefficients in the coresponding
irreducible representation);

--- establish the analogous results for tensor products of irreducibles.
In geometry,
the tensor product of representations corresponds to a "fusion"
of sheaves on configuration spaces --- operation defined
using the functor of nearby cycles, see Section ~\ref{fus}.

We must mention that the assumption that we are dealing with a
Cartan matrix of finite type and a root of unity appears only at the very end
(see Chapter 4). We need these assumptions in order
to compare our representations with the conventional
definition of the category $\CC$. All previous results are valid in more
general assumptions. In particular a
Cartan matrix could be arbitrary and
a deformation parameter $\zeta$ not necessarily a root of unity.

\subsection{} A part of the results of this paper constitutes the description
of the cohomology of certain "standard" local systems over configuration
spaces in terms of quantum groups. These results, due to Varchenko
and one of us, were announced several years ago in ~\cite{sv2}.
The proofs may be found in ~\cite{v}.
Our proof of these results uses completely different approach.
Some results close to this paper were discussed in ~\cite{s}.

Certain results of a similar geometric spirit are discussed in
{}~\cite{fw}.

\subsection{} We are grateful to A.Shen who made our communication during
the writing of this paper possible.

\subsection{}
\label{notations} {\em Notations.} We will use all the notations from
{}~\cite{fs}. References to
{\em loc. cit.} will look like I.1.1.
If $a,b$ are two integers, we will denote by $[a,b]$ the
set of all integers $c$ such that $a\leq c\leq b$; $[1,a]$ will be denoted
by $[a]$.
$\Bbb N$ will denote the set of non-negative integers.
For $r\in\BN$, $\Sigma_r$ will denote the group of all bijections $[r]\iso
[r]$.

We suppose that our ground field $B$ has characteristic $0$,
and fix an element $\zeta\in B$, $\zeta\neq 0$. For $a\in\BZ$ we will use
the notation
\begin{equation}
\label{azeta}
[a]_{\zeta}=1-\zeta^{-2a}
\end{equation}

The word "$t$-exact" will allways mean $t$-exactness with respect to
the middle perversity.

\newpage
\begin{center}
{\bf CHAPTER 1. Algebraic discussion.}
\end{center}
\vspace{1cm}

\section{Free algebras and bilinear forms}

Most definitions of this section follow ~\cite{l} and ~\cite{sv2}
(with slight modifications). We also add some new definitions and
computations important for the sequel. Cf. also ~\cite{v}, Section 4.

\subsection{} Untill the end of this paper, let us fix
a finite set $I$ and a symmetric $\BZ$-valued bilinear form $\nu,\nu'\mapsto
\nu\cdot\nu'$ on the free abelian group $\BZ[I]$ (cf. ~\cite{l}, 1.1).
We will denote
by $X$ the dual abelian group $\Hom(\BZ[I],\BZ)$. Its elements will be
called {\em weights}. Given $\nu\in\BZ[I]$,
we will denote by $\lambda_{\nu}\in X$ the functional $i\mapsto
i\cdot \nu$. Thus we have
\begin{equation}
\label{lamnu}
\langle\lambda_{\nu},\mu\rangle=\nu\cdot\mu
\end{equation}
for all $\nu,\mu\in\BN[I]$.

\subsection{} Let $\fF$ denote a free associative $B$-algebra with
$1$ with generators
$\theta_i,\ i\in I$. Let $\BN[I]$ be a submonoid of $\BZ[I]$ consisting
of all linear combinations of elements of $I$ with coefficients in $\BN$.
For $\nu=\sum\nu_ii\in\BN[I]$ we denote by $\fF_{\nu}$ the $B$-subspace
of $\fF$ spanned by all monomials $\theta_{i_1}\theta_{i_2}\cdot\ldots\cdot
\theta_{i_p}$ such that for any $i\in I$, the number of occurences of $i$
in the sequence $i_1,\ldots,i_p$ is equal to $\nu_i$.

We have a direct sum decomposition $\fF=\oplus_{\nu\in\BN[I]}\fF_{\nu}$,
all spaces $\fF_{\nu}$ are finite dimensional, and we have $\fF_0=B\cdot 1$,
$\fF_{\nu}\cdot\fF_{\nu'}\subset\fF_{\nu+\nu'}$.

Let $\epsilon:\fF\lra B$ denote the augmentation ---
a unique $B$-algebra
map such that $\epsilon(1)=1$ and $\epsilon(\theta_i)=0$ for all $i$. Set
$\fF^+:=\Ker(\epsilon)$. We have $\fF^+=\oplus_{\nu\neq 0}\fF_{\nu}$.

An element $x\in\fF$ is called {\em homogeneous} if it belongs to
$\fF_{\nu}$ for some $\nu$. We then set $|x|=\nu$.
We will use the notation $\depth(x)$ for the number $\sum_i\nu_i$
if $\nu=\sum_i\nu_ii$; it will be called {\em the depth} of $x$.

\subsection{}
\label{twist} Given a sequence
$\vec{K}=(i_1,\ldots,i_N),\ i_j\in I$, let us denote by $\theta_{\vec{K}}$
the monomial
$\theta_{i_1}\cdot\ldots\cdot\theta_{i_N}$. For an empty sequence we set
$\theta_{\emp}=1$.

For
$\tau\in\Sigma_N$
let us introduce the number
\begin{equation}
\label{zetaitau}
\zeta(\vec{K};\tau)=\prod\zeta^{i_a\cdot i_b},
\end{equation}
the product over all $a,b$ such that $1\leq a<b\leq N$ and $\tau(a)>\tau(b)$.

We will call this number {\em the twisting number of the sequence $\vec{K}$
with respect to the permutation $\tau$}.

We will use the notation
\begin{equation}
\label{permut}
\tau(\vec{K})=(i_{\tau(1)},i_{\tau(2)},\ldots,i_{\tau(N)})
\end{equation}

\subsection{}
\label{algff} Let us regard the tensor product $\fF\otimes\fF$ (in
the sequel $\otimes$ will mean $\otimes_B$ unless specified otherwise)
as a $B$-algebra with multiplication
\begin{equation}
\label{mult}
(x_1\otimes x_2)\cdot (x'_1\otimes x_2')=
\zeta^{|x_2|\cdot |x'_1|}x_1x'_1\otimes x_2x'_2
\end{equation}
for homogeneous $x_2,\ x'_1$. Let us define a map
\begin{equation}
\label{comult}
\Delta:\fF\lra\fF\otimes\fF
\end{equation}
as a unique algebra homomorphism carrying $\theta_i$ to $\theta_i\otimes 1+
1\otimes\theta_i$.

\subsection{}
\label{coalgff} Let us define a coalgebra structure on $\fF\otimes\fF$ as
follows. Let us introduce the braiding isomorphism
\begin{equation}
\label{braid}
r:\fF\otimes\fF\iso\fF\otimes\fF
\end{equation}
by the rule
\begin{equation}
\label{braidform}
r(x\otimes y)=\zeta^{|x|\cdot |y|}y\otimes x
\end{equation}
for homogeneous $x,y$. By definition,
\begin{equation}
\label{deltaff}
\Delta_{\fF\otimes\fF}:\fF\otimes\fF\lra(\fF\otimes\fF)\otimes(\fF\otimes\fF)
\end{equation}
coincides with the composition $(1_{\fF}\otimes r\otimes 1_{\fF})\circ
(\Delta_{\fF}\otimes\Delta_{\fF})$.

The multiplication
\begin{equation}
\label{multiplic}
\fF\otimes\fF\lra\fF
\end{equation}
is a coalgebra morphism.

\subsection{}
\label{deltaform} Let us describe $\Delta$ more explicitely. Suppose a sequence
$\vec{K}=(i_1,\ldots,i_N),\ i_j\in I,$ is given. For a  subset
$A=\{ j_1,\ldots, j_a\}\subset [N],\ j_1<\ldots <j_a$, let
$A'=[N]-A=\{k_1,\ldots,k_{N-a}\},\ k_1<\ldots <k_{N-a}$. Define a permutation
$\tau_A$ by the formula
\begin{equation}
\label{taua}
(\tau(1),\ldots,\tau(N))=(j_1,j_2,\ldots,j_a,k_1,k_2,\ldots,k_{N-a})
\end{equation}
Set $\vec{K}_A:=(i_{j_1},i_{j_2},\ldots,i_{j_a}),\
\vec{K}_{A'}:=(i_{k_1},i_{k_2},\ldots,i_{k_{N-a}})$.

\subsubsection{}
\label{delta-expl} {\bf Lemma.} {\em
$$
\Delta(\theta_{\vec{K}})=\sum_{A\subset K}\zeta(\vec{K};\tau_A)
\theta_{\vec{K}_A}\otimes\theta_{\vec{K}_{A'}},
$$
the summation ranging over all subsets $A\subset [N]$.}

{\bf Proof} follows immediately from the definitions. $\Box$

\subsection{}
\label{powers} Let us denote by
\begin{equation}
\label{deltan}
\Delta^{(N)}:\fF\lra\fF^{\otimes N}
\end{equation}
iterated coproducts; by the coassociativity they are well defined.

Let us define a structure of an algebra on $\fF^{\otimes N}$ as follows:
\begin{equation}
\label{product}
(x_1\otimes\ldots\otimes x_N)\cdot (y_1\otimes\ldots\otimes y_N)=
\zeta^{\sum_{j<i}|x_i|\cdot |y_j|} x_1y_1\otimes\ldots\otimes x_Ny_N
\end{equation}
for homogeneous $x_1,\ldots,x_N;y_1,\ldots,y_N$. The map $\Delta^{(N)}$ is an
algebra morphism.

\subsection{}
Suppose we have a sequence $\vec{K}=(i_1,\ldots,i_N)$.
Let us consider an element $\Delta^{(N)}(\theta_{\vec{K}})$;
let $\Delta^{(N)}(\theta_{\vec{K}})^+$ denote its projection to the
subspace $\fF^{+\otimes N}$.

\subsubsection{}
\label{deltamon}
{\bf Lemma.} {\em
$$
\Delta^{(N)}(\theta_{\vec{K}})^+=\sum_{\tau\in\Sigma_N}\zeta(\vec{K};\tau)
\theta_{i_{\tau(1)}}\otimes\ldots\otimes\theta_{i_{\tau(N)}}
$$}

{\bf Proof} follows from ~\ref{delta-expl} by induction on $N$. $\Box$

\subsection{}
\label{duals}
For each component $\fF_{\nu}$ consider the dual $B$-space $\fF^*_{\nu}$,
and set $\fF^*:=\oplus\fF^*_{\nu}$. Graded components
$\Delta_{\nu,\nu'}:\fF_{\nu+\nu'}\lra\fF_{\nu}\otimes\fF_{\nu'}$ define dual
maps $\fF^*_{\nu}\otimes\fF_{\nu'}^*\lra\fF^*_{\nu+\nu'}$ which
give rise to a multiplication
\begin{equation}
\label{mult*}
\fF^*\otimes\fF^*\lra\fF^*
\end{equation}
making $\fF^*$ a graded associative algebra with $1$ (dual to the
augmentation of $\fF$).
This follows from the coassociativity of $\Delta$, cf. ~\cite{l}, 1.2.2.

Here and in the sequel, we will use identifications $(V\otimes W)^*=
V^*\otimes W^*$ (for finite dimensional spaces $V,W$) by the rule
$\langle\phi\otimes\psi,x\otimes y\rangle=\langle\phi,x\rangle\cdot
\langle\psi,y\rangle$.

The dual to ~(\ref{multiplic}) defines a comultiplication
\begin{equation}
\delta:\fF^*\lra\fF^*\otimes\fF^*
\end{equation}
It makes $\fF^*$ a graded coassociative coalgebra with a counit.

The constructions dual to ~\ref{algff} and ~\ref{coalgff}
equip $\fF^*\otimes\fF^*$ with a structure of a coalgebra and an algebra.
It follows from {\em loc. cit} that ~(\ref{mult*}) is a coalgebra
morphism, and $\delta$ is an algebra morphism.

By iterating $\delta$ we get maps
\begin{equation}
\label{deltait}
\delta^{(N)}:\fF^*\lra\fF^{*\otimes N}
\end{equation}
If we regard $\fF^{*\otimes N}$ as an algebra by the same construction as in
{}~(\ref{product}), $\delta^{(N)}$ is an algebra morphism.

\subsection{Lemma.}
\label{form} {\em There exists a unique bilinear form
$$
S(\ ,\ ):\fF\otimes\fF\lra B
$$
such that

(a) $S(1,1)=1$ and $(\theta_i,\theta_j)=\delta_{i,j}$ for all $i,j\in I$;\\
(b) $S(x,y'y'')=S(\Delta(x),y'\otimes y'')$ for all $x,y',y''\in\fF$;\\
(c) $S(xx',y'')=S(x\otimes x',\Delta(y''))$ for all $x,x',y''\in\fF$.

(The bilinear form
$$
(\fF\otimes\fF)\otimes(\fF\otimes\fF)\lra B
$$
given by
$$
(x_1\otimes x_2)\otimes (y_1\otimes y_2)\mapsto S(x_1,y_1)S(x_2,y_2)
$$
is denoted again by $S(\ ,\ )$.)

The bilinear form $S(\ ,\ )$ on $\fF$ is symmetric. The different
homogeneous components $\fF_{\nu}$ are mutually orthogonal. }

{\bf Proof.} See ~\cite{l}, 1.2.3. Cf. also ~\cite{sv2}, (1.8)-(1.11). $\Box$

\subsection{}
\label{epsilons} Following ~\cite{l}, 1.2.13 and ~\cite{sv2}, (1.10)-(1.11),
let us introduce operators $\delta_i:\fF\lra\fF,\ i\in I,$ as unique linear
mappings satisfying
\begin{equation}
\label{formeps}
\delta_i(1)=0;\ \delta_i(\theta_j)=\delta_{i,j},\ j\in I;\
\delta_i(xy)=\delta_i(x)y+\zeta^{|x|\cdot i}x\delta_i(y)
\end{equation}
for homogeneous $x$.

It follows from ~\ref{form} (c) that
\begin{equation}
\label{contrav}
S(\theta_ix,y)=S(x,\delta_i(y))
\end{equation}
for all $i\in I,\ x,y\in \fF$, and obviously $S$ is determined uniquely
by this property, together with the requirement $S(1,1)=1$.

\subsection{}
\label{formulas}
{\bf Lemma.} {\em For any two sequences $\vec{K},\ \vec{K}'$ of $N$ elements
from $I$ we have
$$
S(\theta_{\vec{K}},\theta_{\vec{K}'})=\sum_{\tau\in\Sigma_N:\ \tau(\vec{K})=
\vec{K}'}\zeta(\vec{K};\tau).
$$}

{\bf Proof} follows from ~(\ref{contrav}) by induction on $N$, or else from
{}~\ref{deltamon}. $\Box$

\subsection{} Let us define elements $\theta^*_i\in\fF_i^*$ by the rule
$<\theta^*_i,\theta_i>=1$. The form $S$ defines a homomomorphism of graded
algebras
\begin{equation}
\label{formap}
S:\fF\lra\fF^*
\end{equation}
carrying $\theta_i$ to $\theta^*_i$. $S$ is determined uniquely by
this property.

\subsection{}
\label{morcoalg} {\bf Lemma.} {\em The map $S$ is
a morphism of coalgebras.}

{\bf Proof.} This follows from the symmetry of $S$. $\Box$

\vspace{1cm}
{\em VERMA MODULES}
\vspace{1cm}

\subsection{}
\label{verma}
Let us pick a weight $\Lambda$. Our aim now will be to define
certain $X$-graded vector space $V(\Lambda)$ equipped with the following
structures.

(i) A structure of left $\fF$-module $\fF\otimes V(\Lambda)\lra V(\Lambda)$;

(ii) a structure of left $\fF$-comodule $V(\Lambda)\lra\fF\otimes
V(\Lambda)$;

(iii) a symmetric bilinear form $S_{\Lambda}$ on $V(\Lambda)$.

As a vector space, we set $V(\Lambda)=\fF$. We will define on $V(\Lambda)$
two gradings. The first one, $\BN[I]$-grading coincides with the grading
on $\fF$. If $x\in V(\Lambda)$ is a homogeneous element, we will
denote by $\depth(x)$ its depth as an element of $\fF$.

The second grading --- $X$-grading --- is defined
as follows. By definition, we set
$$
V(\Lambda)_{\lambda}=\oplus_{\nu\in\BN[I]|\Lambda-\lambda_{\nu}
=\lambda}\fF_{\nu}
$$
for $\lambda\in X$. In particular, $V(\Lambda)_{\Lambda}=\fF_0=B\cdot 1$.
We will denote the element $1$ in $V(\Lambda)$ by $v_{\Lambda}$.

By definition, multiplication
\begin{equation}
\label{multv}
\fF\otimes V(\Lambda)\lra V(\Lambda)
\end{equation}
coincides with the  multiplication in $\fF$.

Let us define an $X$-grading in $\fF$ by setting
$$
\fF_{\lambda}=\oplus_{\nu\in\BN[I]|-\lambda_{\nu}=\lambda}\fF_{\nu}
$$
for $\lambda\in X$. The map ~(\ref{multv})
is compatible with both $\BN[I]$ and $X$-gradings (we define gradings
on the tensor product as usually as a sum of gradings of factors).

\subsection{The form $S_{\Lambda}$}
\label{slambda} Let us define linear operators
$\epsilon_i:V(\Lambda)\lra V(\Lambda),\ i\in I,$ as unique
operators such that $\epsilon_i(v_{\Lambda})=0$
and
\begin{equation}
\label{epsilon}
\epsilon_i(\theta_jx)=[\langle\beta,i\rangle]_{\zeta}\delta_{i,j}x +
\zeta^{i\cdot j}\theta_j\epsilon_i(x)
\end{equation}
for $j\in I,\ x\in V(\Lambda)_{\beta}$.

We define $S_{\Lambda}: V(\Lambda)\otimes V(\Lambda)\lra B$ as a unique
linear map such that
$S_{\Lambda}(v_{\Lambda},v_{\Lambda})=1$, and
\begin{equation}
\label{vermacontra}
S_{\Lambda}(\theta_ix,y)=S_{\Lambda}(x,\epsilon_i(y))
\end{equation}
for all $x,y\in V(\Lambda),\ i\in I$.
Let us list elementary properties of $S_{\Lambda}$.

\subsubsection{}
\label{orthog} {\em Different graded components $V(\Lambda)_{\nu},\
\nu\in\BN[I]$, are orthogonal with respect to $S_{\Lambda}$}.

This follows directly from the definition.

\subsubsection{}
\label{ssymmetr} {\em The form $S_{\Lambda}$ is symmetric.}

This is an immediate corollary of the formula
\begin{equation}
\label{epsitheta}
S_{\Lambda}(\epsilon_i(y),x)=S_{\Lambda}(y,\theta_ix)
\end{equation}
which in turn is proved by an easy induction on $\depth(x)$.

\subsubsection{"Quasiclassical" limit}
\label{slclassic} Let us consider restriction
of our form to the homogeneous component $V(\Lambda)_{\lambda}$ of depth
$N$. If we divide our form by $(\zeta -1)^N$ and formally pass
to the limit $\zeta\lra 1$, we get the "Shapovalov" contravariant
form as defined in ~\cite{sv1}, 6.4.1.

The next lemma is similar to ~\ref{formulas}.

\subsection{}
\label{formulasv} {\bf Lemma.} {\em For any $\vec{K},\ \vec{K'}$ as in
{}~\ref{formulas} we have
$$
S_{\Lambda}(\theta_{\vec{K}}v_{\Lambda},\theta_{\vec{K'}}v_{\Lambda})=
\sum_{\tau\in\Sigma_N:\ \tau(\vec{K})=\vec{K'}}
\zeta(\vec{K};\tau)A(\vec{K},\Lambda;\tau)
$$
where
$$
A(\vec{K},\Lambda;\tau)=\prod_{a=1}^N[\langle\Lambda-\sum_{b:\ b<a,\tau(b)<
\tau(a)}\lambda_{i_b},i_a\rangle]_{\zeta}.
$$}

{\bf Proof.} Induction on $N$, using definition of $S_{\Lambda}$. $\Box$

\subsection{Coaction} Let us define a linear map
\begin{equation}
\label{coaction}
\Delta_{\Lambda}:V(\Lambda)\lra\fF\otimes V(\Lambda)
\end{equation}
as follows. Let us introduce linear operators $t_i:\fF^+\otimes V(\Lambda)
\lra\fF^+\otimes V(\Lambda),\ i\in I,$ by the formula
\begin{equation}
\label{ti}
t_i(x\otimes y)=\theta_ix\otimes y-\zeta^{i\cdot\nu-2\langle\lambda,i\rangle}
\cdot x\theta_i\otimes y+\zeta^{i\cdot\nu}x\otimes\theta_iy
\end{equation}
for $x\in\fF_{\nu}$ and $y\in V(\Lambda)_{\lambda}$.

By definition,
\begin{eqnarray}
\label{coactform}
\Delta_{\Lambda}(\theta_{i_N}\cdot\ldots\cdot\theta_{i_1}v_{\Lambda})=
1\otimes\theta_{i_N}\cdot\ldots\cdot\theta_{i_1}v_{\Lambda}\\ \nonumber
+[\langle\Lambda-\lambda_{i_{1}}-\ldots-\lambda_{i_{N-1}},i_N\rangle]_{\zeta}
\cdot\theta_{i_N}\otimes\theta_{i_{N-1}}\cdot\ldots\cdot\theta_{i_1}v_{\Lambda}
\\ \nonumber
+\sum_{j=1}^{N-1}
[\langle\Lambda-\lambda_{i_{1}}-\ldots-\lambda_{i_{j-1}},i_j\rangle]_{\zeta}
\cdot t_{i_N}\circ t_{i_{N-1}}\circ\ldots\circ t_{i_{j+1}}
(\theta_{i_j}\otimes\theta_{i_{j-1}}\cdot\ldots\cdot\theta_{i_1}v_{\Lambda})
\nonumber
\end{eqnarray}

\subsection{}
\label{quantcom} Let us define linear operators
\begin{equation}
\label{qadj}
\ad_{\theta_i,\lambda}:\fF\lra\fF,\ i\in I,\ \lambda\in X
\end{equation}
by the formula
\begin{equation}
\label{quantadj}
\ad_{\theta_i,\lambda}(x)=\theta_ix-\zeta^{i\cdot\nu-2\langle\lambda,i\rangle}
\cdot x\theta_i
\end{equation}
for $x\in\fF_{\nu}$.

Let us note the following relation
\begin{equation}
\label{adjdelta}
(\delta_i\circ\ad_{\theta_j,\lambda}-\zeta^{i\cdot
j}\cdot\ad_{\theta_j,\lambda}
\circ\delta_i)(x)=[\langle\lambda-\lambda_{\nu},i\rangle ]_{\zeta}\delta_{ij}x
\end{equation}
for $x\in\fF_{\nu}$, where $\delta_i$ are operators defined in
{}~\ref{epsilons}, and $\delta_{ij}$ the Kronecker symbol.

\subsection{Formula for coaction}
\label{coactcom} Let us pick a sequence
$\vec{I}=(i_N,i_{N-1},\ldots,i_1)$. To shorten the notations, we set
\begin{equation}
\label{shortcom}
\ad_{j,\lambda}:=\ad_{\theta_{i_j},\lambda},\ j=1,\ldots, N
\end{equation}

\subsubsection{Quantum commutators}
\label{commut} For any non-empty subset $Q\subset [N]$, set
$\theta_{\vec{I},Q}:=\theta_{\vec{I}_Q}$ where $\vec{I}_Q$ denotes the
sequence obtained from $\vec{I}$ by omitting all entries $i_{j},\ j\in Q$.
We will denote $\fF_Q=\fF_{\nu_Q}$ where $\nu_Q:=
\sum_{j\in Q}i_j$.

Let us define an element
$[\theta_{\vec{I},Q,\Lambda}]\in\fF_{Q}$ as follows.
Set
\begin{equation}
\label{comone}
[\theta_{\vec{I},\{j\},\Lambda}]=\zeta^{i_j\cdot(\sum_{k>j}i_k)}
\theta_{i_j}
\end{equation}
for all $j\in [N]$.

Suppose now that $\card(Q)=l+1\geq 2$.
Let $Q=\{j_0,j_1,\ldots,j_l\},\ j_0<j_1<\ldots <j_l$.
Define the weights
$$
\lambda_a=\Lambda-\lambda_{\sum i_k},\ a=1,\ldots, l,
$$
where the summation is over $k$ from $1$ to $j_a-1$, $k\neq j_1,j_2,\ldots,
j_{a-1}$.

Let us define sequences $\vec{N}:=(N,N-1,\ldots,1),\
\vec{Q}=(j_l,j_{l-1},\ldots,j_0)$ and
$\vec{N}_Q$ obtained from $\vec{N}$ by omitting all entries $j\in Q$.
Define the permutation $\tau_Q\in \Sigma_N$ by the requirement
$$
\tau_Q(\vec{N})=\vec{Q}||\vec{N}_Q
$$
where $||$ denotes concatenation.

Set by definition
\begin{equation}
\label{comgen}
[\theta_{\vec{I},Q,\Lambda}]=\zeta(\vec{I},\tau_Q)\cdot
\ad_{j_l,\lambda_l}\circ\ad_{j_{l-1},\lambda_{l-1}}\circ\ldots\circ
\ad_{j_1,\lambda_1}(\theta_{i_{j_0}})
\end{equation}

\subsubsection{}
\label{formcoact} {\bf Lemma.} {\em We have
\begin{equation}
\label{sumcoact}
\Delta_{\Lambda}(\theta_{\vec{I}}v_{\Lambda})=
1\otimes \theta_{\vec{I}}v_{\Lambda}+
\sum_{Q}[\langle\Lambda-\lambda_{i_1}-\lambda_{i_2}-\ldots -
\lambda_{i_{j(Q)-1}},i_{j(Q)}\rangle]_{\zeta}\cdot
[\theta_{\vec{I},Q,\Lambda}]\otimes\theta_{\vec{I},Q}v_{\Lambda},
\end{equation}
the summation over all non-empty subsets $Q\subset [N]$,
$j(Q)$ denotes the minimal element of $Q$.}

{\bf Proof.} The statement of the lemma follows at once from the
inspection of definition ~(\ref{coactform}), after rearranging
the summands. $\Box$

Several remarks are in order.

\subsubsection{} Formula ~(\ref{sumcoact}) as similar to ~\cite{s}, 2.5.4.

\subsubsection{} If all elements $i_j$ are distinct then
the part of the sum  in the rhs of ~(\ref{sumcoact})
corresponding to one-element subsets $Q$ is
equal to
$\sum_{j=1}^N\theta_{i_j}\otimes\epsilon_{i_j}(\theta_{\vec{I}}v_{\Lambda})$.

\subsubsection{"Quasiclassical" limit}
It follows from the definition of quantum commutators that if
we divide the rhs of ~(\ref{sumcoact}) by $(\zeta -1)^N)$
and formally pass to the limit $\zeta\lra 1$, we get the expression for
the coaction
obtained in ~\cite{sv1}, 6.15.3.2.

\subsection{} Let us define the space $V(\Lambda)^*$ as the direct sum
$\oplus_{\nu}V(\Lambda)^*_{\nu}$. We define an $\BN[I]$-grading on it as
$(V(\Lambda)^*)_{\nu}=V(\Lambda)^*_{\nu}$, and an $X$-grading as
$V(\Lambda)^*_{\lambda}=\oplus_{\nu:\ \Lambda-\lambda_{\nu}=\lambda}
V(\Lambda)^*_{\nu}$.

The form $S_{\Lambda}$ induces the map
\begin{equation}
\label{maps}
S_{\Lambda}:V(\Lambda)\lra V(\Lambda)^*
\end{equation}
compatible with both gradings.

\subsection{Tensor products}
\label{tensprod} Suppose we are given $n$ weights
$\Lambda_0,\ldots,\Lambda_{n-1}$.

\subsubsection{} For every $m\in\BN$ we introduce a bilinear
form $S=S_{m;\Lambda_0,\ldots,\Lambda_{n-1}}$ on the tensor product
$\fF^{\otimes m}\otimes V(\Lambda_0)
\otimes\ldots\otimes V(\Lambda_{n-1})$ by the formula
$$
S(x_1\otimes\ldots\otimes x_m\otimes y_0\otimes\ldots
\otimes y_{n-1},
x'_1\otimes\ldots\otimes x'_m\otimes y'_0\otimes
\ldots\otimes y'_{n-1})=\prod_{i=1}^mS(x_i,x'_i)\prod_{j=0}^{n-1}
S_{\Lambda_j}(y_j,y'_j)
$$
(in the evident notations). This form defines mappings
\begin{equation}
\label{stens}
S:\fF^{\otimes m}\otimes V(\Lambda_0)\otimes\ldots\otimes V(\Lambda_{n-1})
\lra
\fF^{*\otimes m}\otimes V(\Lambda_0)^*\otimes\ldots\otimes V(\Lambda_{n-1})^*
\end{equation}

\subsubsection{} We will regard $V(\Lambda_0)
\otimes\ldots\otimes V(\Lambda_{n-1})$ as an $\fF^{\otimes n}$-module with
an action
\begin{equation}
\label{multiter}
(u_0\otimes\ldots\otimes u_{n-1})\cdot (x_1\otimes\ldots\otimes x_{n-1})=
\zeta^{-\sum_{j<i}\langle\lambda_j,\nu_i\rangle}u_0x_0\otimes\ldots\otimes
u_{n-1}x_{n-1}
\end{equation}
for $u_i\in\fF_{\nu_i},\ x_j\in V(\Lambda_j)_{\lambda_j}$, cf.
{}~(\ref{product}). Here we regard $\fF^{\otimes n}$ as an algebra according
to the rule of {\em loc. cit.}; one checks easily using ~(\ref{lamnu})
that we really get a module structure.

Using the iterated comultiplication $\Delta^{(n)}$, we get
a structure of an $\fF$-module on
$V(\Lambda_0)\otimes\ldots\otimes V(\Lambda_{n-1})$.

\subsection{Theorem}
\label{coactshap} {\em We have an identity
\begin{equation}
\label{coactshapo}
S_{\Lambda}(xy,z)=S_{1;\Lambda}(x\otimes y,\Delta_{\Lambda}(z))
\end{equation}
for any $x\in\fF,\ y,z\in V(\Lambda)$ and any weight $\Lambda$.}

\subsection{Proof} We may suppose that $x,y$ and $z$ are monomials. Let
$z=\theta_{\vec{I}}v_{\Lambda}$ where $\vec{I}=(i_N,\ldots,i_1)$.

(a) Let us suppose first that all indices $i_j$ are distinct. We will
use the notations and computations from ~\ref{coactcom}. The sides of
{}~(\ref{coactshapo}) are non-zero only if $y$ is equal to
$z_Q:=\theta_{\vec{I},Q}v_{\Lambda}$ for some subset $Q\subset [N]$.

Therefore, it follows from Lemma ~\ref{formcoact} that it is enough to prove

\subsubsection{} {\bf Lemma.} {\em For every non-empty $Q\subset [N]$ and
$x\in\fF_Q$ we have
\begin{equation}
\label{cominduct}
S_{\Lambda}(xz_Q,z)=
[\langle\Lambda,i_{j(Q)}\rangle-\mu_Q\cdot i_{j(Q)}]_{\zeta}
\cdot S(x,[\theta_{\vec{I},Q,\Lambda}])\cdot S_{\Lambda}(z_Q,z_Q)
\end{equation}
where $j(Q)$ denotes the minimal element of $Q$, and
$$
\mu_Q:=\sum_{a=1}^{j(Q)-1}i_a.
$$}

{\bf Proof.} If $\card(Q)=1$
the statement follows from the definiton ~(\ref{comone}). The proof will
proceed
by the simultaneous induction by $l$ and $N$.
Suppose that $x=\theta_{i_p}\cdot x'$, so $x'\in\fF_{Q'}$ where
$Q'=Q-\{i_p\}$, $p=j_a$ for some $a\in [0,l]$. Let us set $\vec{I}'=
\vec{I}-\{i_p\},\ z'=z_{\{i_p\}}$, so that $z_Q=z'_{Q'}$.

We have
\begin{eqnarray}
S_{\Lambda}(\theta_{i_p}x'\cdot z_Q,z)=
S_{\Lambda}(x'\cdot z_Q,\epsilon_{i_p}(z))=\\ \nonumber
=[\langle\Lambda,i_p\rangle-(\sum_{k<p}i_k)\cdot i_p]_{\zeta}\cdot
\zeta^{(\sum_{k>p}i_k)\cdot i_p}\cdot
S(x'\cdot z_{Q},z')=\\ \nonumber
=[\langle\Lambda,i_p\rangle-(\sum_{k<p}i_k)\cdot i_p]_{\zeta}\cdot
[\langle\Lambda,i_{j(Q')}\rangle-\mu_{Q'}\cdot i_{j(Q')}]_{\zeta}\cdot
\zeta^{(\sum_{k>p}i_k)\cdot i_p}\cdot \\ \nonumber
\cdot S(x,[\theta_{\vec{I},Q,\Lambda}])\cdot S_{\Lambda}(z'_{Q'},z'_{Q'})
\nonumber
\end{eqnarray}
by induction hypothesis.
On the other hand,
$$
S(\theta_{i_p}\cdot x',[\theta_{\vec{I},Q,\Lambda}])=
S(x',\delta_{i_p}([\theta_{\vec{I},Q,\Lambda}])).
$$
Therefore, to complete the induction step it is enough to prove that
\begin{eqnarray}
\label{deltacommut}
[\langle\Lambda,i_{j(Q)}\rangle-\mu_Q\cdot i_{j(Q)}]_{\zeta}\cdot
\delta_{i_p}([\theta_{\vec{I},Q,\Lambda}])=\\ \nonumber
=[\langle\Lambda,i_p\rangle-(\sum_{k<p}i_k)\cdot i_p]_{\zeta}\cdot
[\langle\Lambda,i_{j(Q')}\rangle-\mu_{Q'}\cdot i_{j(Q')}]_{\zeta}\cdot
\zeta^{(\sum_{k>p}i_k)\cdot i_p}
[\theta_{\vec{I}',Q',\Lambda}] \nonumber
\end{eqnarray}
This formula follows directly from the definition of
quantum commutators ~(\ref{comgen}) and formula ~(\ref{adjdelta}).
One has to treat separately two cases: $a>0$, in which case $j(Q)=j(Q')=j_0$
and
$a=0$, in which case $j(Q)=j_0,\ j(Q')=j_1$. Lemma is proven.
$\Box$

This completes the proof of case (a).

(b) There are repeating indices in the sequence $\vec{I}$.
Suppose that $\theta_{\vec{I}}\in\fF_{\nu}$. At this point we will use
symmetrization constructions (and simple facts) from Section ~\ref{symmetr}
below. The reader will readily see that there is no vicious circle.
So, this part of the proof must be read after {\em loc.cit.}

There exists a finite set $J$ and a map $\pi:J\lra I$ such that
$\nu=\nu_{\pi}$.
Using compatibility of the coaction and the forms $S$ with symmetrization ---
cf. Lemmata ~\ref{avercoact} and ~\ref{averterms} below ---
our claim is immediately
reduced to the analogous claim for the algebra $^{\pi}\fF$,
the module $V(^{\pi}\Lambda)$ and homogeneous weight $\chi_J$
which does not contain multiple indices and therefore follows from (a) above.

This completes the proof of the theorem. $\Box$

\subsection{} Let us pick a weight $\Lambda$.
We can consider numbers $q_{ij}:=\zeta^{i\cdot j}$ and $r_i:=
\langle\Lambda,i\rangle$,  $i,j\in I$ as parameters of our bilinear
forms.

More precisely, for a given $\nu\in\BN[I]$ the matrix elements
of the form $S$ (resp., $S_{\Lambda}$) on $\fF_{\nu}$ (resp.,
on $V(\Lambda)_{\Lambda-\lambda_{\nu}}$) in the standard bases of these spaces
are certain universal polynomials of $q_{ij}$ (resp., $q_{ij}$ and
$r_i$). Let us denote their determinants by $\det(S_{\nu})(\bq)$ and
$\det(S_{\Lambda,\nu})(\bq;\br)$ respectively. These determinants are
polynomials of corresponding variables with integer coefficients.

\subsubsection{}
\label{generic} {\bf Lemma.} {\em Polynomials $\det(S_{\nu})(\bq)$ and
$\det(S_{\Lambda,\nu})(\bq;\br)$ are not identically zero.

In other words, bilinear forms $S$ and $S_{\Lambda}$ are non-degenerate
for generic values of parameters  --- "Cartan matrix" $(q_{ij})$ and
"weight" $(r_i)$.}

{\bf Proof.} Let us consider the form $S_{\Lambda}$ first.
The specialization of the matrix of $S_{\Lambda,\nu}$ at $\zeta=1$ is the
identity matrix. It follows easily that $\det(S_{V,\nu})(\bq;\br)\neq 0$.

Similarly, the matrix of $S_{\nu}$ becomes identity at $\zeta=0$,
which implies the generic non-degeneracy.  $\Box$

\subsection{Theorem}
\label{coactas} {\em Coaction $\Delta_{\Lambda}$ is coassociative,
i.e.
\begin{equation}
\label{coasform}
(1_{\fF}\otimes\Delta_{\Lambda})\circ\Delta_{\Lambda}=
(\Delta\otimes 1_{V(\Lambda)})\circ\Delta_{\Lambda}.
\end{equation}}

{\bf Proof.} The equality ~(\ref{coasform}) is a polynomial identity
depending on parameters $q_{ij}$ and $r_i$ of the preceding subsection.
For generic values of these parameters it is true due to associativity
of the action of $\fF$ an $V(\Lambda)$, Theorem ~\ref{coactshap} and
Lemma ~\ref{generic}. Therefore it is true for all values
of parameters. $\Box$

\subsection{} The results Chapter 2 below provide a different,
geometric proof of Theorems ~\ref{coactshap} and ~\ref{coactas}.
Namely, the results of Section ~\ref{standsheaves} summarized in Theorem
{}~\ref{phi*} provide an isomorphism of our algebraic picture
with a geometric one, and in the geometrical language the above theorems
are obvious: they
are nothing but the naturality of the canonical morphism
between the extension by zero and the extension by star,
and the claim that a Cousin complex is a complex. Lemma ~\ref{generic}
also follows from geometric considerations: the extensions by zero and
by star coincide for generic values of monodromy.

\subsection{}
\label{comod} By Theorem ~\ref{coactas} the dual maps
\begin{equation}
\label{deltal*}
\Delta_{\Lambda}^*:\fF\otimes V(\Lambda)^*\lra V(\Lambda)^*
\end{equation}
give rise to a structure
of a $\fF^*$-module on $V(\Lambda)^*$.

More generally, suppose we are given $n$ modules $V(\Lambda_0),\ldots,
V(\Lambda_{n-1})$. We regard the tensor product
$V(\Lambda_0)^*\otimes\ldots\otimes V(\Lambda_{n-1})^*$ as a
$\fF^{*\otimes n}$-module according to the "sign" rule ~(\ref{multiter}).
Using iterated comultiplication ~(\ref{deltait}) we get a structure
of a $\fF^*$-module on $V(\Lambda_0)^*\otimes\ldots\otimes
V(\Lambda_{n-1})^*$.

\subsubsection{}
\label{scomod} The square
$$\begin{array}{ccccc}
\;&\fF\otimes V(\Lambda_0)\otimes\ldots\otimes
V(\Lambda_{n-1})&\lra&V(\Lambda_0)\otimes\ldots\otimes
V(\Lambda_{n-1})&\;\\ \nonumber
\;&S\downarrow&\;&\downarrow S&\; \\ \nonumber
\;&\fF^*\otimes V(\Lambda_0)^*\otimes\ldots\otimes
V(\Lambda_{n-1})^*&\lra&V(\Lambda_0)^*\otimes\ldots\otimes
V(\Lambda_{n-1})^*&\;
\end{array}$$
commutes.

This follows from ~\ref{coactshap} and ~\ref{morcoalg}.

\section{Hochschild complexes}

\subsection{} If $A$ is an augmented $B$-algebra, $A^+$ --- the kernel
of the augmentation, $M$ an $A$-module, let $C_A^{\bullet}(M)$ denote
the following complex. By definition, $C_A^{\bullet}(M)$ is concentrated
in non-positive degrees. For $r\geq 0$
$$
C_A^{-r}(M)=A^{+\otimes r}\otimes M.
$$
We will use a notation $a_r|\ldots|a_1|m$ for $a_r\otimes\ldots a_1\otimes m$.

The differential $d:C_A^{-r}(M)\lra C_A^{-r+1}(M)$ acts as
$$
d(a_r|\ldots|a_1|m)=
\sum_{p=1}^{r-1}(-1)^{p}a_r|\ldots|a_{p+1}a_p|\ldots a_1|m+
a_r|\ldots a_{2}|a_1m.
$$
We have canonically $H^{-r}(C_A^{\bullet}(M))\cong \Tor^A_r(B,M)$ where
$B$ is considered as an $A$-module by means of the augmentation,
cf. ~\cite{m}, Ch. X, \S 2.

We will be interested in the algebras $\fF$ and $\fF^*$. We define
the augmentation
$\fF\lra B$ as being zero on all $\fF_{\nu},\ \nu\in\BN[I],\ \nu\neq 0$,
and identity on $\fF_0$; in the same way it is defined on $\fF^*$.

\subsection{}
\label{gradings} Let $M$ be a $\BN[I]$-graded $\fF$-module.
Each term $C_{\fF}^{-r}(M)$ is $\BN[I]$-graded
by the sum of gradings of tensor factors.
We will denote $_{\nu}C_{\fF}^{-r}(M)$ the weight $\nu$ component.

For $\bnu=(\nu_0,\ldots,\nu_r)\in \BN[I]^{r+1}$ we set
$$
_{\bnu}C_{\fF}^{-r}(M)=\ _{\nu_r,\ldots,\nu_0}C_{\fF}^{-r}(M)=
\fF_{\nu_r}\otimes\ldots\otimes\fF_{\nu_1}\otimes M_{\nu_0}.
$$
Thus,
$$
_{\nu}C_{\fF}^{-r}(M)=\oplus_{\nu_0+\ldots\nu_r=\nu}\
_{\nu_r,\ldots,\nu_0}C_{\fF}^{-r}(M).
$$
Note that all $\nu_p$ must be $>0$ for $p>0$ since tensor factors lie in
$\fF^+$.

The differential $d$ clearly respects the  $\BN[I]$-grading; thus the whole
complex is $\BN[I]$-graded:
$$
C_{\fF}^{\bullet}(M)=\oplus_{\nu\in\BN[I]}\
_{\nu}C_{\fF}^{\bullet}(M).
$$
The same discussion applies to $\BN[I]$-graded $\fF^*$-modules.

\subsection{}
\label{gradehoch} Let us fix weights
$\Lambda_0,\ldots,\Lambda_{n-1}$, $n\geq 1$. We will consider the Hochschild
complex
$C^{\bullet}_{\fF}(V(\Lambda_0)\otimes\ldots\otimes V(\Lambda_{n-1}))$
where the structure of an $\fF$-module on
$V(\Lambda_0)\otimes\ldots\otimes V(\Lambda_{n-1})$ has been introduced
in ~\ref{tensprod}.

\subsubsection{} In the sequel we will use the following notation.
If $K\subset I$ is a subset, we will denote by $\chi_K:=\sum_{i\in K}i
\in \BN[I]$.

\subsubsection{}
\label{gh} Suppose we have a map
\begin{equation}
\label{maprho}
\varrho: I\lra [-n+1, r]
\end{equation}
where $r$ is some non-negative integer. Let us introduce the elements
\begin{equation}
\label{nurho}
\nu_a(\varrho)=\chi_{\varrho^{-1}(a)},
\end{equation}
$a\in [-n+1,r]$.
Let us denote by $\CP_r(I;n)$ the set of all maps ~(\ref{maprho}) such that
$\varrho^{-1}(a)\neq\emp$ for all $a\in[r]$. It is easy to see that this set
is not empty iff $0\leq r\leq N$.

Let us assign to such a $\varrho$ the space
\begin{equation}
\label{homog}
_{\varrho}C^{-r}_{\fF}(V(\Lambda_{0})\otimes\ldots\otimes V(\Lambda_{n-1})):=
\fF_{\nu_r(\varrho)}\otimes\ldots\otimes\fF_{\nu_1(\varrho)}\otimes
V(\Lambda_{0})_{\nu_0(\varrho)}\otimes\ldots\otimes
V(\Lambda_{n-1})_{\nu_{-n+1}(\varrho)}
\end{equation}
For each $\varrho\in\CP_r(I;n)$ this space is non-zero, and we have
\begin{equation}
\label{decompos}
_{\chi_I}C^{-r}_{\fF}(V(\Lambda_{0})\otimes\ldots\otimes V(\Lambda_{n-1}))=
\oplus_{\varrho\in\CP_r(I;n)}\
_{\varrho}C^{-r}_{\fF}(V(\Lambda_{0})\otimes\ldots\otimes V(\Lambda_{n-1}))
\end{equation}

\subsection{Bases}
\label{refine} Let us consider the set $\CP_N(I;n)$. Obviously, if
$\varrho\in\CP_N(I;n)$ then $\varrho(I)=[N]$, and the induced
map $I\lra [N]$ is a bijection; this way we get an isomorphism
between $\CP_N(I;n)$ and the set of all bijections $I\iso [N]$ or,
to put it differently, with the set of all total orders on $I$.

For an arbitrary $r$, let $\varrho\in\CP_r(I;n)$ and $\tau\in\CP_N(I;n)$.
Let us say that $\tau$ is {\em a refinement of $\varrho$}, and
write $\varrho\leq\tau$, if
$\varrho(i)<\varrho(j)$ implies $\tau(i)<\tau(j)$ for each $i,j\in I$.
The map $\tau$ induces total orders on all subsets $\varrho^{-1}(a)$.
We will denote by $\Ord(\varrho)$ the set of all refinements of a given
$\varrho$.

Given $\varrho\leq\tau$ as above, and $a\in [-n+1,r]$, suppose that
$\varrho^{-1}(a)=\{ i_1,\ldots,i_p\}$ and
$\tau(i_1)<\tau(i_2)<\ldots<\tau(i_p)$. Let us define a monomial
$$
\theta_{\varrho\leq\tau;a}=\theta_{i_p}\theta_{i_{p-1}}\cdot\ldots\cdot
\theta_{i_1}\in\fF_{\nu_a(\varrho)}
$$
If $\varrho^{-1}(a)=\emp$, we set $\theta_{\varrho\leq\tau;a}=1$.
This defines a monomial
\begin{equation}
\label{monom}
\theta_{\varrho\leq\tau}=\theta_{\varrho\leq\tau;r}\otimes\ldots\otimes
\theta_{\varrho\leq\tau;1}\otimes\theta_{\varrho\leq\tau;0}v_{\Lambda_0}\otimes
\ldots\otimes\theta_{\varrho\leq\tau;-n+1}v_{\Lambda_{n-1}}\in\
_{\varrho}C^{-r}_{\fF}(V(\Lambda_{0})\otimes\ldots\otimes V(\Lambda_{n-1}))
\end{equation}

\subsubsection{}
\label{basis} {\bf Lemma.} {\em The set  $\{\theta_{\varrho\leq\tau}|
\tau\in\Ord(\varrho)\}$ forms a basis of the space
{}~{$_{\varrho}C^{-r}_{\fF}(V(\Lambda_{0})\otimes\ldots\otimes
V(\Lambda_{n-1}))$.}}

{\bf Proof} is obvious. $\Box$

\subsubsection{}
\label{basisj} {\bf Corollary.} {\em The set  $\{\theta_{\varrho\leq\tau}|
\varrho\in\CP_r(I;n),\ \tau\in\Ord(\varrho)\}$ forms a basis of the space
$_{\chi_I}C^{-r}_{\fF}(V(\Lambda_{0})\otimes\ldots\otimes
V(\Lambda_{n-1}))$.} $\Box$

\subsection{}
\label{dualbases} We will also consider dual Hochschild complexes
$C^{\bullet}_{\fF^*}(V(\Lambda_{0})^*\otimes\ldots\otimes
V(\Lambda_{n-1})^*)$ where $V(\Lambda_{0})^*\otimes\ldots\otimes
V(\Lambda_{n-1})^*$ is regarded as an $\fF^*$-module as in ~\ref{scomod}.

We have obvious isomorphisms
$$
C^{-r}_{\fF^*}(V(\Lambda_{0})^*\otimes\ldots\otimes
V(\Lambda_{n-1})^*)\cong
C^{-r}_{\fF}(V(\Lambda_{0})\otimes\ldots\otimes
V(\Lambda_{n-1}))^*
$$
We define graded components
$$
_{\varrho}C^{-r}_{\fF^*}(V(\Lambda_{0})^*\otimes\ldots\otimes
V(\Lambda_{n-1})^*),\
\varrho\in\CP_r(I;n),
$$
as duals to $_{\varrho}C^{-r}_{\fF}(V(\Lambda_{0})\otimes\ldots\otimes
V(\Lambda_{n-1}))$.

We will denote by $\{\theta_{\varrho\leq\tau}^*|
\varrho\in\CP_r(I;n),\ \tau\in\Ord(\varrho)\}$ the basis of
$_{\chi_I}C^{-r}_{\fF^*}(V(\Lambda_{0})^*\otimes\ldots\otimes
V(\Lambda_{n-1})^*)$ dual to the basis
$\{\theta_{\varrho\leq\tau}|
\varrho\in\CP_r(I;n),\ \tau\in\Ord(\varrho)\}$, ~\ref{basisj}.

\subsection{} The maps $S_{r;\Lambda_0,\ldots,\Lambda_{n-1}}$, cf.
{}~(\ref{stens}), for different $r$ are compatible with differentials
in Hochschild complexes, and therefore induce morphism
of complexes
\begin{equation}
\label{shapohoch}
S:C^{\bullet}_{\fF}(V(\Lambda_{0})\otimes\ldots\otimes
V(\Lambda_{n-1})) \lra
C^{\bullet}_{\fF^*}(V(\Lambda_{0})^*\otimes\ldots\otimes
V(\Lambda_{n-1})^*)
\end{equation}
This follows from ~\ref{scomod} and ~\ref{form} (b).

\section{Symmetrization}
\label{symmetr}

\subsection{}
\label{sym} Let us fix a finite set $J$ and a map $\pi:J\lra I$.
We set $\nu_{\pi}:=\sum_iN_ii\in\BZ[I]$ where
$N_i:=\card(\pi^{-1}(i))$.
The map $\pi$ induces a map $\BZ[J]\lra \BZ[I]$ also
to be denoted by $\pi$.
We will use the notation $\chi_K:=\sum_{j\in K}j\in\BN[J]$ for $K\subset J$.
Thus, $\pi(\chi_J)=\nu_{\pi}$.

We will denote also by $\mu,\mu'\mapsto
\mu\cdot\mu':=\pi(\mu)\cdot\pi(\mu')$ the bilinear form on $\BZ[J]$ induced
by the form on $\BZ[I]$.

We will denote by $\Sigma_{\pi}$ the group of all
bijections $\sigma:J\lra J$ preserving fibers of $\pi$.

Let $^{\pi}\fF$ be a free associative $B$-algebra with $1$ with generators
$\ttheta_j,\ j\in J$. It is evidently $\BN[J]$-graded.
For $\nu\in\BN[J]$ the corresponding homogeneous component will be denoted
$^{\pi}\fF_{\nu}$. The degree of a
homogeneous element $x\in\ ^{\pi}\fF$ will be denoted by $|x|\in\BN[J]$.
The group $\Sigma_{\pi}$ acts on algebras $^{\pi}\fF,\ ^{\pi}\fF^*$ by
permutation of generators.

\subsection{}
\label{average} In the sequel, if $G$ is a group and $M$ is a $G$-module,
$M^G$ will denote the subset of $G$-invariants in $M$.

Let us define a $B$-linear "averaging" mapping
\begin{equation}
\label{avermap}
^{\pi}a:\fF_{\nu_{\pi}}\lra (^{\pi}\fF_{\chi_J})^{\Sigma_{\pi}}
\end{equation}
by the rule
\begin{equation}
\label{avermapform}
^{\pi}a(\theta_{i_1}\cdot\ldots\cdot\theta_{i_N})=
\sum \ttheta_{j_1}\cdot\ldots\cdot \ttheta_{j_N},
\end{equation}
the sum being taken over the set of all sequences $(j_1,\ldots,j_N)$ such
that $\pi(j_p)=i_p$ for any $p$. Note that this set is naturally a
$\Sigma_{\pi}$-torsor.
Alternatively, $^{\pi}a$ may be defined as follows. Pick some sequence
$(j_1,\ldots,j_N)$ as above, and consider an element
$$
\sum_{\sigma\in\Sigma_{\pi}}\sigma(\ttheta_{j_1}\cdot\ldots\cdot
\ttheta_{j_N});
$$
this element obviously lies in $(^{\pi}\fF_{\chi_J})^{\Sigma_{\pi}}$ and
is equal to $^{\pi}a(\theta_{i_1}\cdot\ldots\cdot\theta_{i_N})$.

The map $\pi$ induces the map between homogeneous components
\begin{equation}
\label{ainv}
\pi:\ ^{\pi}\fF_{\chi_J}\lra \fF_{\nu_{\pi}}.
\end{equation}
It is clear that the composition $\pi\circ\ ^{\pi}a$ is equal to the
multiplication by $\card (\Sigma_{\pi})$, and $^{\pi}a\circ\pi$ ---
to the action of operator $\sum_{\sigma\in\Sigma_{\pi}}\sigma$.
As a consequence, we get

\subsubsection{}
{\bf Lemma,} ~\cite{sv1}, 5.11. {\em The map $^{\pi}a$ is an isomorphism.}
$\Box$

\subsection{}
\label{averdual} Let us consider the dual to the map ~(\ref{ainv}):
$\fF_{\nu_{\pi}}^*\lra\ ^{\pi}\fF_{\chi_J}^*$; it is obvious that it
lands in the subspace of $\Sigma_{\pi}$-invariant functionals. Let us
consider the induced map
\begin{equation}
\label{avermap*}
^{\pi}a^*:\fF_{\nu_{\pi}}^*\iso (^{\pi}\fF_{\chi_J}^*)^{\Sigma_{\pi}}
\end{equation}
It follows from the above discussion that $^{\pi}a^*$ is an isomorphism.

\subsection{}
\label{averhoch} Given a weight $\Lambda\in X=\Hom(\BZ[I],\BZ)$, we will denote
by $^{\pi}\Lambda$ the composition $\BZ[J]\overset{\pi}{\lra}\BZ[I]
\overset{\Lambda}{\lra}\BZ$, and by $V(^{\pi}\Lambda)$ the corresponding
Verma module over $^{\pi}\fF$.

Suppose we are given $n$ weights $\Lambda_0,\ldots,\Lambda_{n-1}$.
Let us consider the Hochschild complex
{}~{$C^{\bullet}_{^{\pi}\fF}(V(^{\pi}\Lambda_0)\otimes\ldots\otimes
V(^{\pi}\Lambda_{n-1}))$}. By definition, its $(-r)$-th term coincides with
the tensor power $^{\pi}\fF^{\otimes n+r}$. Therefore we can identify
the homogeneous component
$_{\chi_J}C^{-r}_{^{\pi}\fF}(V(^{\pi}\Lambda_0)\otimes\ldots\otimes
V(^{\pi}\Lambda_{n-1}))$ with $(^{\pi}\fF^{\otimes n+r})_{\chi_J}$ which in
turn is isomorphic to $^{\pi}\fF_{\chi_J}$, by means of the multiplication
map $^{\pi}\fF^{\otimes n+r}\lra\ ^{\pi}\fF$. This defines a map
\begin{equation}
\label{hochcompj}
_{\chi_J}C^{-r}_{^{\pi}\fF}(V(^{\pi}\Lambda_0)\otimes\ldots\otimes
V(^{\pi}\Lambda_{n-1}))\lra\ ^{\pi}\fF_{\chi_J}
\end{equation}
which is an embedding when restricted to polygraded components.
The $\Sigma_{\pi}$-action on $\fF$ induces the $\Sigma_{\pi}$-action on
$_{\chi_J}C^{-r}_{^{\pi}\fF}(V(^{\pi}\Lambda_0)\otimes\ldots\otimes
V(^{\pi}\Lambda_{n-1}))$.

In the same manner we define a map
\begin{equation}
\label{hochcomp}
_{\nu_{\pi}}C^{-r}_{\fF}(V(\Lambda_0)\otimes\ldots\otimes
V(\Lambda_{n-1}))\lra\fF_{\nu_{\pi}}
\end{equation}

Let us define an averaging map
\begin{equation}
\label{avermaphoch}
^{\pi}a:\ _{\nu_{\pi}}C^{-r}_{\fF}(V(\Lambda_0)\otimes\ldots\otimes
V(\Lambda_{n-1}))\lra\
_{\chi_J}C^{-r}_{^{\pi}\fF}(V(^{\pi}\Lambda_0)\otimes\ldots\otimes
V(^{\pi}\Lambda_{n-1}))^{\Sigma_{\pi}}
\end{equation}
as the map induced by ~(\ref{avermap}). It follows at once that this map is
is an isomorphism.

These maps for different $r$ are by definition compatible with
differentials in Hochschild complexes. Therefore we get

\subsubsection{}
\label{avehoch} {\bf Lemma.} {\em The maps ~(\ref{avermaphoch})
induce isomorphism of complexes
\begin{equation}
\label{averhochiso}
^{\pi}a:\ _{\nu_{\pi}}C^{\bullet}_{\fF}(V(\Lambda_0)\otimes\ldots\otimes
V(\Lambda_{n-1}))\iso\
_{\chi_J}C^{\bullet}_{^{\pi}\fF}(V(^{\pi}\Lambda_0)\otimes\ldots\otimes
V(^{\pi}\Lambda_{n-1}))^{\Sigma_{\pi}}.\ \Box
\end{equation}}

\subsection{}
\label{avercoact} {\bf Lemma.} {\em The averaging is compatible with coaction.
In other words, for any $\Lambda\in X$ the square
$$\begin{array}{ccccc}
\;&V(\Lambda)&\overset{\Delta_{\Lambda}}{\lra}&\fF\otimes V(\Lambda)&\;\\
\;&^{\pi}a\downarrow&\;&\downarrow\ ^{\pi}a&\;\\
\;&V(^{\pi}\Lambda)&\overset{\Delta_{^{\pi}\Lambda}}{\lra}&^{\pi}\fF
                                              \otimes V(^{\pi}\Lambda)&\;
\end{array}$$
commutes.}

{\bf Proof} follows at once by inspection of the definition ~(\ref{coactform}).
$\Box$

\subsection{}
\label{averhoch*} Consider the dual Hochschild complexes.
We have an obvious isomorphism
$$
_{\chi_J}C^{-r}_{^{\pi}\fF^*}(V(^{\pi}\Lambda_0)^*\otimes\ldots\otimes
V(^{\pi}\Lambda_{n-1})^*)\cong\
_{\chi_J}C^{-r}_{^{\pi}\fF}(V(^{\pi}\Lambda_0)\otimes\ldots\otimes
V(^{\pi}\Lambda_{n-1}))^*;
$$
using it, we define the isomorphism
$$
_{\chi_J}C^{-r}_{^{\pi}\fF^*}(V(^{\pi}\Lambda_0)^*\otimes\ldots\otimes
V(^{\pi}\Lambda_{n-1})^*)\iso\ ^{\pi}\fF^*_{\chi_J}
$$
as the dual to ~(\ref{hochcompj}). The $\Sigma_{\pi}$-action on the target
induces the action on
$_{\chi_J}C^{-r}_{^{\pi}\fF^*}(V(^{\pi}\Lambda_0)^*\otimes\ldots\otimes
V(^{\pi}\Lambda_{n-1})^*)$.
Similarly, the isomorphism
$$
_{\nu_{\pi}}C^{-r}_{\fF^*}(V(\Lambda_0)^*\otimes\ldots\otimes
V(\Lambda_{n-1})^*)\iso \fF^*_{\nu_{\pi}}
$$
is defined. We define the averaging map
\begin{equation}
\label{avermaphoch*}
^{\pi}a^*: _{\nu_{\pi}}C^{-r}_{\fF^*}(V(\Lambda_0)^*\otimes\ldots\otimes
V(\Lambda_{n-1})^*)\lra
_{\chi_J}C^{-r}_{^{\pi}\fF^*}(V(^{\pi}\Lambda_0)^*\otimes\ldots\otimes
V(^{\pi}\Lambda_{n-1})^*)^{\Sigma_{\pi}}
\end{equation}
as the map which coincides with ~(\ref{avermap*}) modulo the above
identifications. Again, this map is an isomorphism.

Due to Lemma ~\ref{avercoact} these maps for different $r$ are compatible
with the differentials in Hochschild complexes. Therefore we get

\subsubsection{}
\label{avehoch*} {\bf Lemma.} {\em The maps ~(\ref{avermaphoch*})
induce isomorphism of complexes
\begin{equation}
\label{averhochiso*}
^{\pi}a^*: _{\nu_{\pi}}C^{\bullet}_{\fF^*}(V(\Lambda_0)^*\otimes\ldots\otimes
V(\Lambda_{n-1})^*)\iso\
_{\chi_J}C^{\bullet}_{^{\pi}\fF^*}(V(^{\pi}\Lambda_0)^*\otimes\ldots\otimes
V(^{\pi}\Lambda_{n-1})^*)^{\Sigma_{\pi}}.\ \Box
\end{equation}}

\vspace{1cm}
{\em BILINEAR FORMS}
\vspace{1cm}

\subsection{} Using the bilinear form on $\BZ[J]$ introduced above,
we define the symmetric bilinear form
$S(\ ,\ )$ on $^{\pi}\fF$ exactly in the same way as the form $S$ on
$\fF$. Similarly, given $\Lambda\in X$, we define the bilinear
form $S_{^{\pi}\Lambda}$ on $V(^{\pi}\Lambda)$ as in ~\ref{slambda},
with $I$ replaced by $J$.

\subsubsection{}
\label{S-sym}
{\bf Lemma.} (i) {\em The square
$$\begin{array}{ccccc}
\;&\fF_{\nu_{\pi}}&\overset{S}{\lra}&\fF_{\nu_{\pi}}^*&\;\\
\;&^{\pi}a\downarrow&\;&\downarrow\ ^{\pi}a^*&\;\\
\;&^{\pi}\fF_{\chi_J}&\overset{S}{\lra}&^{\pi}\fF_{\chi_J}^*&\;
\end{array}$$
commutes.}

(ii) {\em For any $\Lambda\in X$ the square
$$\begin{array}{ccccc}
\;&V(\Lambda)_{\nu_{\pi}}&\overset{S_{\Lambda}}{\lra}&V(\Lambda)_{\nu_{\pi}}^*
&\;\\
\;&^{\pi}a\downarrow&\;&\downarrow\ ^{\pi}a^*&\;\\
\;&V(^{\pi}\Lambda)_{\chi_J}&\overset{S_{^{\pi}\Lambda}}{\lra}&
V(^{\pi}\Lambda)_{\chi_J}^*&\;
\end{array}$$
commutes.}

{\bf Proof.} (i) Let us consider an element
$\theta_{\vec{I}}=\theta_{i_1}\cdot\ldots\cdot\theta_{i_N}\in\fF_{\nu_{\pi}}$
(we assume that $N=\card (J)$).
The functional $^{\pi}a^*\circ S(\theta_{\vec{I}})$ carries a monomial
$\ttheta_{j_1}\cdot\ldots\cdot\ttheta_{j_N}$ to
$$
S(\theta_{i_1}\cdot\ldots\cdot\theta_{i_N},
\theta_{\pi(j_1)}\cdot\ldots\cdot\theta_{\pi(j_N)}).
$$
On the other hand,
$$
S\circ\ ^{\pi}a(\theta_{\vec{I}})(\ttheta_{j_1}\cdot\ldots\cdot\ttheta_{j_N})=
\sum S(\ttheta_{k_1}\cdot\ldots\cdot\ttheta_{k_N},
\ttheta_{j_1}\cdot\ldots\cdot\ttheta_{j_N}),
$$
the summation ranging over all sequences $\vec{K}=(k_1,\ldots,k_N)$ such that
$\pi(\vec{K})=\vec{I}$. It follows from Lemma ~\ref{formulas} that
both expressions are equal.

(ii) The same argument as in (i), using Lemma ~\ref{formulasv} instead
of ~\ref{formulas}. $\Box$

More generally, we have

\subsection{Lemma}
\label{averterms} {\em For every $m\geq 0$ and weights $\Lambda_0,\ldots,
\Lambda_{n-1}\in X$ the square
$$\begin{array}{ccccc}
\;&(\fF^{\otimes m}\otimes V(\Lambda_0)\otimes\ldots\otimes
V(\Lambda_{n-1}))_{\nu_{\pi}}&
\overset{S_{m;\Lambda_0,\ldots,\Lambda_{n-1}}}{\lra}&
(\fF^{*\otimes m}\otimes V(\Lambda_0)^*\otimes\ldots\otimes
V(\Lambda_{n-1})^*)_{\nu_{\pi}}
&\;\\
\;&^{\pi}a\downarrow&\;&\downarrow\ ^{\pi}a^*&\;\\
\;&(^{\pi}\fF^{\otimes m}\otimes V(^{\pi}\Lambda_0)\otimes\ldots\otimes
V(^{\pi}\Lambda_{n-1}))_{\chi_J}
&\overset{S_{m;^{\pi}\Lambda_0,\ldots,^{\pi}\Lambda_{n-1}}}{\lra}&
(^{\pi}\fF^{*\otimes m}\otimes V(^{\pi}\Lambda_0)^*\otimes\ldots\otimes
V(^{\pi}\Lambda_{n-1})^*)_{\chi_J}&\;
\end{array}$$
commutes.}

{\bf Proof} is quite similar to the proof of the previous lemma.
We leave it to the reader. $\Box$

\section{Quotient algebras}

\subsection{} Let us consider the map ~(\ref{formap})
$S:\fF\lra\fF^*$. Let us consider its kernel $\Ker (S)$.
It follows at once from ~(\ref{contrav}) that $\Ker (S)$ is a left
ideal in $\fF$. In the same manner, it is easy to see that
it is also a right ideal, cf. ~\cite{l}, 1.2.4.

We will denote by $\ff$ the quotient algebra $\fF/\Ker(S)$.
It inherits the $\BN[I]$-grading and the coalgebra structure from
$\fF$, cf. {\em loc.cit.} 1.2.5, 1.2.6.

\subsection{}
\label{llambda} In the same manner, given a weight $\Lambda$, consider
the kernel
of $S_{\Lambda}:V(\Lambda)\lra V(\Lambda)^*$. Let us denote
by $L(\Lambda)$ the quotient space $V(\Lambda)/\Ker(S_{\Lambda})$.
It inherits $\BN[I]$- and $X$-gradings from $V(\Lambda)$.
Due to Theorem ~\ref{coactshap} the structure of $\fF$-module
on $V(\Lambda)$ induces the structure of $\ff$-module
on $L(\Lambda)$.

More generally, due to the structure of a coalgebra on $\ff$, all
tensor products $L(\Lambda_0)\otimes\ldots\otimes L(\Lambda_{n-1})$
become $\ff$-modules (one should take into account the "sign rule"
{}~(\ref{multiter})).

\subsection{} We can consider Hochschild complexes
$C^{\bullet}_{\ff}(L(\Lambda_0)\otimes\ldots\otimes L(\Lambda_{n-1}))$.

\subsubsection{} {\bf Lemma.} {\em The map
$$
S:C^{\bullet}_{\fF}(V(\Lambda_0)\otimes\ldots\otimes V(\Lambda_{n-1}))\lra
C^{\bullet}_{\fF^*}(V(\Lambda_0)^*\otimes\ldots\otimes V(\Lambda_{n-1})^*)
$$
factors through the isomorphism
\begin{equation}
\label{ims}
\Ima(S)\iso
C^{\bullet}_{\ff}(L(\Lambda_0)\otimes\ldots\otimes L(\Lambda_{n-1}))
\end{equation}}

{\bf Proof.} This follows at once from the definitions. $\Box$

\newpage
\begin{center}
{\bf CHAPTER 2. Geometric discussion.}
\end{center}
\vspace{1cm}

\section{Diagonal stratification and related algebras}

\subsection{}
\label{diagsetup} Let us adopt notations of ~\ref{sym}. We set
$N:=\card(J)$.
Let $^{\pi}\BA_{\BR}$
denote a real affine space with coordinates $t_j,\ j\in J$, and
$^{\pi}\BA$ its complexification. Let us consider an arrangement $\CH_{\emp}$
consisting of all diagonals $\Delta_{ij},\ i,j\in J$. Let us denote
by $\CS_{\emp}$ the corresponding stratification;
$\CS_{\emp,\BR}$ will denote the corresponding real stratification of
$\BA_{\BR}$.

The stratification $\CS_{\emp}$ has a unique minimal stratum
\begin{equation}
\label{delta}
\Delta=\bigcap\ \Delta_{ij}
\end{equation}
--- main diagonal; it is one-dimensional.
We will denote by $^{\pi}\BAO_{\emp}$ (resp., $^{\pi}\BAO_{\emp,\BR}$)
the open stratum
of $\CS_{\emp}$ (resp., of $\CS_{\emp,\BR}$).

\subsection{}
\label{diagchamb} Let us describe the chambers of $\CS_{\emp,\BR}$.
If $C$ is a chamber and
$\bx=(x_j)\in C$, i.e. the embedding $J\hra\BR,\ j\mapsto x_j$,
it induces an obvious total order on $J$, i.e. a bijection
\begin{equation}
\label{tauc}
\tau_C:J\iso [N]
\end{equation}
Namely, $\tau_C$ is determined uniquely by the requirement
$\tau_C(i)<\tau_C(j)$ iff $x_i<x_j$; it does not depend on the choice
of $\bx$. This way we get  a one-to-one correspondence between
the set of chambers of $\CS_{\emp}$ and the set of all bijections
{}~(\ref{tauc}). We will denote by $C_{\tau}$ the chamber corresponding to
$\tau$.

Given $C$ and $\bx$ as above, suppose that we have $i,j\in J$ such that
$x_i<x_j$ and there is no $k\in J$ such that $x_i<x_k<x_j$. We will say that
$i,j$ are {\em neighbours} in $C$, more precisely that {\em $i$ is a left
neighbour of $j$}.

Let $\bx'=(x'_j)$ be a point with $x'_p=x_p$ for all $p\neq j$, and $x'_j$
equal to some number smaller than $x_i$ but greater than any $x_k$ such that
$x_k<x_i$.
Let $^{ji}C$ denote the chamber containing $\bx'$. Let us introduce a  homotopy
class of paths $^C\gamma_{ij}$ connecting $\bx$ and $\bx'$ as shown on
Fig. 1 below.

\begin{picture}(20,8)(-10,-4)

\put(0,0){\circle*{0.2}}
\put(0,-0.5){$i$}

\put(-4,0){\line(1,0){8}}

\put(1.5,0){\circle*{0.2}}
\put(1.5,-0.5){$j$}

\put(0,0){\oval(3,3)[t]}
\put(0,1.5){\vector(-1,0){0.5}}
\put(-0.5,2){$^C\gamma_{ij}$}

\put(-1.5,0){\circle*{0.2}}
\put(-1.5,-0.5){$j'$}

\put(3,0){\circle*{0.2}}
\put(-3,0){\circle*{0.2}}

\put(-0.5,-4){Fig. 1.}

\end{picture}

We can apply the discussion I.4.1
and consider the groupoid $\pi_1(^{\pi}\BAO_{\emp},^{\pi}\BAO_{\emp,\BR})$.
It has as the set of objects the set
of all chambers. The set of morphisms is generated by all morphisms
$^C\gamma_{ij}$
subject to certain evident braiding relations. We will need only the following
particular case.

To define a {\em one-dimensional} local system $\CL$ over $^{\pi}\BAO_{\emp}$
is
the same as to give a set of one-dimensional vector spaces
$\CL_C,\ C\in\pi_0(^{\pi}\BAO_{\emp,\BR})$,
together with arbitrary invertible linear operators
\begin{equation}
\label{halfmonodr}
^CT_{ij}:\CL_C\lra\CL_{^{ji}C}
\end{equation}
("half-monodromies")  defined for chambers having $i$ as a left neighbour of
$j$.

\subsection{}
\label{pici} We define a one-dimensional local system $^{\pi}\CI$ over
$^{\pi}\BAO_{\emp}$ as follows. Its fibers $^{\pi}\CI_{C}$
are one-dimensional
linear spaces with fixed basis vectors; they will be identified with $B$.

Half-monodromies are defined as
$$
^CT_{ij}=\zeta^{i\cdot j},\ i,j\in J
$$

\subsection{} Let $j:\ ^{\pi}\BAO_{\emp}\lra\ ^{\pi}\BA$ denote an open
embedding. We will study the
following objects of $\CM(^{\pi}\BA;\CS_{\emp})$:
$$
^{\pi}\CI_{?}=j_?^{\pi}\CI[N],
$$
where $?=!,*$. We have a canonical map
\begin{equation}
\label{map!*}
m:\ ^{\pi}\CI_!\lra\ ^{\pi}\CI_*
\end{equation}
and by definition $^{\pi}\CI_{!*}$ is its image, cf. I.4.5.

\subsection{}
\label{marking} For an integer $r$ let us denote by $\CP_r(J)$ the set of all
surjective mappings $J\lra [r]$. It is evident that $\CP_r(J)\neq\emp$ if and
only if $1\leq r\leq N$. To each $\rho\in\CP_r(J)$ let us assign a point
$w_{\rho}=(\rho(j))\in\ ^{\pi}\BA_{\BR}$. Let $F_{\rho}$ denote
the facet containing $w_{\rho}$. This way we get a bijection between
$\CP_r(J)$ and the set of $r$-dimensional facets. For $r=N$ we get the
bijection from ~\ref{diagchamb}.

At the same time we have defined a marking of $\CH_{\emp}$:
by definition, $^{F_{\rho}}w=w_{\rho}$. This defines cells $D_F,\ S_F$.

\subsection{}
The main diagonal $\Delta$ is a unique $1$-facet; it corresponds to the
unique element $\rho_0\in\CP_1(J)$.

We will denote by $\Ch$ the set of all chambers; it is the same as
$\Ch(\Delta)$ in notations of Part I.
Let $C_{\tau}$ be a chamber. The order $\tau$ identifies $C$ with
an open cone in the standard coordinate space $\BR^N$; we provide $C$
with the orientation induced from $\BR^N$.

\subsection{Basis in $\Phi_{\Delta}(^{\pi}\CI_*)$}
The construction
I.4.7 gives us the basis $\{c_{\Delta<C}\}$ in
$\Phi_{\Delta}(^{\pi}\CI_!)^*$ indexed
by $C\in\Ch$. We will use notation $c_{\tau,!}:=c_{\Delta<C_{\tau}}$.

A chain $c_{\tau,!}$ looks as follows.

\begin{picture}(20,6)(-10,-3)

\put(-4,0){\line(1,0){6}}

\put(-1,0){\vector(1,0){0.8}}
\put(-0.2,0){\circle*{0.15}}
\put(-0.3,-0.4){$j_1$}
\put(-1,0){\vector(1,0){1.6}}
\put(0.6,0){\circle*{0.15}}
\put(0.5,-0.4){$j_2$}

\put(2.2,0){$\ldots$}
\put(3,0){\line(1,0){2}}
\put(3,0){\vector(1,0){0.5}}
\put(3.5,0){\circle*{0.15}}
\put(3.5,-0.4){$j_N$}

\put(-1,-3){Fig. 2. A chain $c_{\tau,!}$.}

\end{picture}

Here $\tau(j_i)=i$. We will denote by $\{b_{\tau,!}\}$ the dual basis
in $\Phi_{\Delta}(\CI_!)$.

\subsection{Basis in $\Phi_{\Delta}(\CI_*)$}
Similarly, the definition I.4.9 gives us the basis $\{c_{\Delta<C}\},\
C\in\Ch$ in $\Phi_{\Delta}(\CI_*)^*$. We will use the notations
$c_{\tau,*}:=c_{\Delta<C_{\tau}}$.

If we specify
the definition I.4.9 and
its explanation I.4.12 to our arrangement, we get the following
picture for a dual chain $c_{\tau,*}$.

\begin{picture}(20,6)(-10,-3)

\put(-5,0){\line(1,0){9.5}}

\put(2,0){\circle*{0.15}}
\put(1.9,0.4){$j_1$}

\put(2,0){\oval(1.6,1.6)[b]}
\put(2,-0.8){\vector(1,0){0.3}}
\put(2.8,0){\circle*{0.15}}
\put(2.8,0.4){$j_2$}

\put(2,0){\oval(3.2,2.4)[b]}
\put(2,-1.2){\vector(1,0){0.8}}
\put(3.6,0){\circle*{0.15}}
\put(3.6,0.4){$j_3$}

\put(4.5,0){$\ldots$}
\put(5.5,0){\line(1,0){2}}

\put(2,0){\oval(10,4)[b]}
\put(2,-2){\vector(1,0){1.8}}
\put(7,0){\circle*{0.15}}
\put(7,0.4){$j_N$}

\put(-1,-3){Fig. 3. A chain $c_{\tau,*}$ .}

\end{picture}

This chain is represented by the section of a local system $\CI^{-1}$
over the cell in $^{\pi}\BAO_{\emp}$ shown above, which takes value $1$
at the point corresponding to the end of the travel (direction of travel
is shown by arrows).

To understand what is going on, it is instructive to treat the
case $N=2$ first, which essentially coincides with the Example I.4.10.

We will denote by $\{b_{\tau,*}\}$ the dual basis in $\Phi_{\Delta}(\CI_*)$.

\subsection{}
\label{signs} Obviously, all maps $\tau:J\lra [N]$ from $\CP_N(J)$ are
bijections. Given two such maps $\tau_1,\tau_2$, define the sign
$\sgn(\tau_1,\tau_2)=\pm 1$ as the sign of the permutation
$\tau_1 \tau_2^{-1}\in\Sigma_N$.

For any $\tau\in\CP_N(J)$ let us denote by $\vec{J}_{\tau}$ the sequence
$(\tau^{-1}(N),\tau^{-1}(N-1),\ldots,\tau^{-1}(1))$.

\subsection{} Let us pick $\eta\in\CP_N(J)$. Let us define
the following maps:
\begin{equation}
\label{phid!j}
^{\pi}\phi_{\Delta,!}^{(\eta)}:\Phi_{\Delta}(^{\pi}\CI_!)\lra\
^{\pi}\fF_{\chi_J}
\end{equation}
which carries $b_{\tau,!}$ to $\sgn(\tau,\eta)\cdot\theta_{\vec{J}_{\tau}}$,
and
\begin{equation}
\label{phid*j}
^{\pi}\phi_{\Delta,*}^{(\eta)}:\Phi_{\Delta}(^{\pi}\CI_*)\lra\
^{\pi}\fF^*_{\chi_J}
\end{equation}
which carries $b_{\tau,*}$ to
$\sgn(\tau,\eta)\cdot\theta^*_{\vec{J}_{\tau}}$.

\subsection{Theorem}
\label{diagiso} {\em (i) The maps $^{\pi}\phi_{\Delta,!}^{(\eta)}$ and
$^{\pi}\phi_{\Delta,*}^{(\eta)}$ are isomorphisms.
The square
$$\begin{array}{ccc}
\Phi_{\Delta}(^{\pi}\CI_!)&\overset{^{\pi}\phi_{\Delta,!}^{(\eta)}}{\iso}&
\ ^{\pi}\fF_{\chi_J}\\
m\downarrow&\;&\downarrow S\\
\Phi_{\Delta}(^{\pi}\CI_*)&\overset{^{\pi}\phi_{\Delta,*}^{(\eta)}}{\iso}&
\ ^{\pi}\fF^*_{\chi_J}
\end{array}$$
commutes.

(ii) The map $^{\pi}\phi_{\Delta,!}^{(\eta)}$ induces an isomorphism
\begin{equation}
\label{phi!*iso}
^{\pi}\phi_{\Delta,!*}^{(\eta)}:\Phi_{\Delta}(^{\pi}\CI_{!*})\iso\
^{\pi}\ff_{\chi_J}
\ \Box
\end{equation}}

{\bf Proof.} This theorem is particular case of I.14.16, I.4.17.
The claim about isomorphisms in (i) is clear. To prove the commutativity
of the square, we have to compute the action of the canonical map $m$
on our standard chains. The claim follows at once from their geometric
description given above. Note that here the sign in the definition
of morphisms $\phi$ is essential, due to orientations of our chains.
(ii) is a direct corollary of (i) $\Box$

\vspace{1cm}
{\em SYMMETRIZED CONFIGURATIONAL SPACES}
\vspace{1cm}

\subsection{Colored configuration spaces}
\label{color} Let us fix $\nu=\sum\nu_ii\in\BN[I]$, $\sum_i\nu_i=N$.
There exists a finite set $J$ and a morphism $\pi:J\lra I$ such that
$\card(\pi^{-1}(i))=\nu_i$ for all $i\in I$. Let us call such $\pi$
{\em an unfolding of $\nu$}. It is unique up to a non-unique isomorphism;
the automorphism group of $\pi$ is precisely $\Sigma_{\pi}$, and
$\nu=\nu_{\pi}$ in our previous notations.

Let us pick an unfolding $\pi$. As in the above discussion, we define
$^{\pi}\BA$ as a complex affine space with coordinates $t_j,\ j\in J$.
Thus, $\dim\ ^{\pi}\BA=N$.
The group
$\Sigma_{\pi}$ acts on the space $^{\pi}\BA$ by permutations
of coordinates.

Let us denote by $\CA_{\nu}$ the quotient manifold $^{\pi}\BA/\Sigma_{\pi}$.
As an
algebraic manifold, $\CA_{\nu}$ is also a complex $N$-dimensional affine
space. We have a canonical projection
\begin{equation}
\label{mappi}
\pi:\ ^{\pi}\BA\lra\CA_{\nu}
\end{equation}
The space $\CA_{\nu}$ does not depend on the choice of an unfolding $\pi$.
It will be called {\em the configuration space of $\nu$-colored points
on the affine line $\BA^1$}.

We will consider the stratification on $\CA_{\nu}$ whose strata are
$\pi(S),\ S\in\CS_{\emp}$; we will denote this stratification also by
$\CS_{\emp}$; this definition does not depend on the choice of $\pi$.
We will study the category $\CM(\CA_{\nu};\CS_{\emp})$.

We will denote by $\CAO_{\nu,\emp}$ the open stratum. It is clear that
$\pi^{-1}(\CAO_{\nu,\emp})=\ ^{\pi}\BAO_{\emp}$. The morphism $\pi$
is unramified over $\CAO_{\nu,\emp}$.

The action of $\Sigma_{\pi}$ on
$^{\pi}\BA$ may be extended in the evident
way to the local system $^{\pi}\CI$,
hence all our spaces of geometric origin ---
like $\Phi_{\Delta}(^{\pi}\CI_!)$, etc. ---
get an action of $\Sigma_{\pi}$.

\subsubsection{} If $M$ is an object with a $\Sigma_{\pi}$-action
(for example a vector
space or a sheaf), we will denote by $M^{\Sigma_{\pi},-}$ the subobject
$\{x\in M|\mbox{ for every }\sigma\in\Sigma_{\pi}\ \sigma x=\sgn(\sigma)x\}$
where $\sgn(\sigma)=\pm 1$ is the sign of a permutation.

A morphism $f:M\lra N$ between two objects with $\Sigma_{\pi}$-action will
be called {\em skew ($\Sigma_{\pi}$)-equivariant} if for any
$x\in M,\ \sigma\in\Sigma_{\pi},\ f(\sigma x)=\sgn(\sigma)\sigma f(x)$.

Let us define a local system over $\CA_{\nu}$
\begin{equation}
\label{cinu}
\CI_{\nu}=(\pi_*\ ^{\pi}\CI)^{\Sigma_{\pi},-}
\end{equation}

\subsection{}
\label{jcolor} Let $j:\ ^{\pi}\BAO_{\emp}\hra\ ^{\pi}\BA,\
j_{\CA}:\CAO_{\nu,\emp}\hra\CA_{\nu}$ be the open embeddings. Let us define
the following objects of $\CM(\CA_{\nu};\CS_{\emp})$:
\begin{equation}
\label{standca}
\CI_{\nu?}:=j_{\CA?}\CI_{\nu}[N]
\end{equation}
where $?=!,\ *$ or $!*$.
We have by definition
\begin{equation}
\label{inu!}
\CI_{\nu!}=(\pi_*\ ^{\pi}\CI_!)^{\Sigma_{\pi},-}
\end{equation}
The morphism $\pi$ is finite; consequently $\pi_*$ is $t$-exact
(see ~\cite{bbd}, 4.1.3) and commutes with the Verdier duality. Therefore,
\begin{equation}
\label{inu*}
\CI_{\nu*}=(\pi_*\ ^{\pi}\CI_*)^{\Sigma_{\pi},-};\
\CI_{\nu!*}=(\pi_*\ ^{\pi}\CI_{!*})^{\Sigma_{\pi},-}
\end{equation}

\subsection{} Let us define vector spaces
\begin{equation}
\label{phidelt}
\Phi_{\Delta}(\CI_{\nu?}):=(\Phi_{\Delta}(^{\pi}\CI_{?}))^{\Sigma_{\pi},-}
\end{equation}
where $?=!,*$ or $!*$.

Let us pick a $\Sigma_{\pi}$-equivariant marking of $\CH_{\emp}$, for
example the one from ~\ref{marking}; consider the corresponding
cells $D_{\Delta},\ S_{\Delta}$. It follows from ~(\ref{inu!}) and
{}~(\ref{inu*}) that
\begin{equation}
\label{relat}
\Phi_{\Delta}(\CI_{\nu?})=R\Gamma(\pi(D_{\Delta}),\pi(S_{\Delta});\
\CI_{\nu?})[-1]
\end{equation}
where $?=!,*$ or $!*$, cf. I.3.3.

\subsection{} The group $\Sigma_{\pi}$ is acting on
on $^{\pi}\fF$. Let us pick $\eta\in\CP_N(J)$.
It follows from the definitions that
the isomorphisms $^{\pi}\phi_{\Delta,!}^{(\eta)},\
^{\pi}\phi_{\Delta,*}^{(\eta)}$ are skew
$\Sigma_{\pi}$-equivariant. Therefore,
passing to invariants in Theorem ~\ref{diagiso} we get

\subsection{Theorem}
\label{symdiag} {\em The maps $^{\pi}\phi_{\Delta,!}^{(\eta)},
^{\pi}\phi_{\Delta,*}^{(\eta)}$ induce isomorphisms included into a commutative
square
$$\begin{array}{ccc}
\Phi_{\Delta}(\CI_{\nu!})&\overset{\phi_{\nu,!}^{(\eta)}}{\iso}&\fF_{\nu}\\
m\downarrow&\;&\downarrow S\\
\Phi_{\Delta}(\CI_{\nu*})&\overset{\phi_{\nu,*}^{(\eta)}}{\iso}&\fF^*_{\nu}
\end{array}$$
and
\begin{equation}
\label{phinu}
\phi_{\nu,!*}^{(\eta)}:\Phi_{\Delta}(\CI_{\nu!*})\iso\ff_{\nu}\ \ \ \Box
\end{equation}}

\section{Principal stratification}

The contents of this section is parallel to I, Section 3. However, we present
here certain modification of general constructions from {\em loc. cit.}

\subsection{}
\label{j} Let us fix a finite set $J$ of cardinality $N$. In this section
we will denote by $\BA_{\BR}$ a real affine space with fixed coordinates
$t_j:\BA_{\BR}\lra\BR,\ j\in J$, and by $\BA$ its complexification.
For $z\in \BC,\ i,j\in J$
denote by $H_j(z)\subset \BA$ a hyperplane $t_j=z$, and
by $\Delta_{ij}$ a hyperplane $t_i=t_j$.

Let us consider an arrangement
$\CH$ in $\Bbb A$ consisting of hyperplanes $H_i(0)$ and $\Delta_{ij}$,
$i,j\in J,\ i\neq j$. It is a complexification of an evident real arrangement
$\CH_{\BR}$ in $\BA_{\BR}$. As usual, the subscript $_{\BR}$ will denote
real points.

Denote by $\CS$ the corresponding stratification of $\BA$. To distinguish this
stratification from the diagonal stratification of the previous
section, we will call it {\em the principal stratification}.
To shorten the
notation, we will denote in this paper by $\CD(\CA,\CS)$ a category which would
be denoted $\CD^b(\BA;\CS)$ in I. In this section
we will study the category $\CM(\BA;\CS)$.

\subsubsection{}
\label{posfacets} Let us consider a positive cone
$$
\BA^+_{\BR}=\{ (t_j)|\mbox{ all }t_j\geq 0\}\subset\BA_{\BR}
$$
A facet will be called {\em positive} if it lies inside
$\BA^+_{\BR}$.

A flag $\bF$ is called {\em positive} if all its facets
are positive.

\subsection{}
\label{poscells} Let us fix a marking $\bw=\{\ ^Fw\}$ of $\CH_{\BR}$
(cf. I.3.2).
For a positive facet $F$ define
$$
D^+_F=D_F\cap\BA_{\BR}^+;\ S^+_F=S_F\cap\BA_{\BR}^+;\
\DO^+_F=D_F^+-S_F^+.
$$
Note that $D^+_F$ coincides with the union of $^{\bF}\Delta$ over all
positive flags beginning at $F$, and $S^+_F$ coincides with the union
of $^{\bF}\Delta$ as above with $\dim\ ^{\bF}\Delta<\codim\ F$. It follows
that only marking points $^Fw$ for positive facets $F$ take part in the
definition of cells $D^+_F,\ S^+_F$.

\subsection{} Let $\CK$ be an object of $\CD(\BA;\CS)$, $F$ a
positive facet of dimension $p$.
Let us introduce a notation
$$
\Phi_F^+(\CK)=\Gamma(D_F^+,S_F^+;\CK)[-p].
$$
This way we get a functor
\begin{equation}
\label{plus}
\Phi_F^+:\CD(\BA;\ \CS)\lra\CD^b(pt)
\end{equation}

\subsection{Theorem.}
\label{dual+} {\em Functors $\Phi_F^+$ commute
with Verdier duality. More precisely, we have canonical natural
isomorphisms
\begin{equation}
D\Phi_F^+(\CK)\iso\Phi_F^+(D\CK).
\end{equation}}

{\bf Proof} goes along the same lines as the proof of Theorem
I.3.5.

\subsection{}
\label{onedim} First let us consider the case $N=1$, cf. I.3.6. We will adopt
notations
from there and from I, Fig. 1. Our arrangement has one positive
$1$-dimensional facet $E=\BR_{>0}$, let $w\in E$ be a marking.

\begin{picture}(20,8)(-10,-4)

\put(0,0){\circle{0.2}}
\put(-0.5,-0.4){$F$}
\put(0,0){\oval(6,6)}
\put(-1.5,3.3){$S_{r''}$}

\put(0.1,0){\line(1,0){4}}
\put(0.5,0.3){$E$}

\put(0,3){\line(0,-1){2}}
\put(-0.4,1.8){$Y$}

\put(0,1){\circle*{0.2}}
\put(-0.5,0.5){$i\epsilon w$}

\put(2,0){\circle*{0.2}}
\put(2,-0.5){$w$}

\put(0,0){\oval(3,3)}
\put(-1.5,1.6){$S_{r'}$}

\put(-0.5,-4){Fig. 4}

\end{picture}

We have an isomorphism
\begin{equation}
\label{iso}
\Phi_{F}^+(\CK)\iso R\Gamma(\BA,\{ w\};\CK)
\end{equation}
Denote $j:=j_{\BA-\{ w\}}$. We have by Poincar\'{e} duality
\begin{equation}
D\Phi_F^+(\CK)\iso R\Gamma_c(\BA;j_*j^*D\CK)\iso
R\Gamma(\BA,\BA_{\geq r''};j_*j^*D\CK)
\end{equation}
Let us denote $D^{+opp}_F:=\BR_{\geq w}$, and $Y=\epsilon i\cdot D_F^{+opp}$.
We have
\begin{equation}
\label{homot}
R\Gamma(\BA,\BA_{\geq r''};j_*j^*D\CK)\iso
R\Gamma(\BA,Y\cup \BA_{\geq r''};j_*j^*D\CK)
\end{equation}
by homotopy. Consider the restriction map
\begin{equation}
\label{res}
res: R\Gamma(\BA,Y\cup \BA_{\geq r''};j_*j^*D\CK)\lra
R\Gamma(\BA_{\leq r'},Y\cap \BA_{\leq r'};D\CK)
\end{equation}
\subsubsection{}
\label{qiso} {\bf Claim.} {\em $res$ is an isomorphism.}

In fact, $\Cone (res)$ is isomorphic to
\begin{eqnarray}
R\Gamma(\BA,\BA_{\leq r'}\cup\BA_{\geq r''}\cup Y;j_*j^*D\CK)=
R\Gamma_c(\BA_{< r''},\BA_{\leq r'}\cup Y;j_*j^*D\CK)\cong\nonumber\\
\cong DR\Gamma(\BA_{< r''}-(\BA_{\leq r'}\cup Y);j_!j^*\CK)\nonumber
\end{eqnarray}
But
\begin{equation}
R\Gamma(\BA_{< r''}-(\BA_{\leq r'}\cup Y);j_!j^*\CK)=
R\Gamma(\BA_{< r''}-(\BA_{\leq r'}\cup Y),\{ w\};\CK)=0
\end{equation}
evidently. This proves the claim. $\Box$

A clockwise rotation by $\pi/2$ induces an isomorphism
$$
R\Gamma(\BA_{\leq r'},Y\cap \BA_{\leq r'};D\CK)\cong
R\Gamma(\BA_{\leq r'},\epsilon\cdot D_F^{opp};D\CK),
$$
and the last complex is isomorphic to $\Phi_F^+(D\CK)$. This proves the theorem
for $N=1$.

\subsection{}
\label{reduct} Now let us return to an arbitrary $J$. Let us prove
the theorem for $F$ equal to the unique
$0$-dimensional facet.

Let us introduce the following subspaces of
$\BA_{\BR}$ (as usually, a circle on the top will denote the interior).

$D_F^{opp}:=\BA_{\BR}-\DO_F$;
$D_F^{+opp}:=\BA^+_{\BR}-\DO_F^+$;
for each cell $D_{F<C},\ C\in \Ch(F)$ define $D_{F<C}^{opp}
:=C-\DO_{F<C}$.

It is easy to
see that the restriction induces isomorphism
$$
\Phi_F^+(\CK)\iso R\Gamma(\BA,D^{+opp}_F;\CK).
$$

We use the notations of I.3.8. Let us choose positive numbers
$r'<r'',\ \epsilon,$ such that
$$
\epsilon D_F\subset \BA_{<r'}\subset \DO_F\subset D_F\subset
\BA_{< r''}
$$
Define the subspace
$$
Y^+=\epsilon i\cdot D_F^{+opp};
$$
denote $j:=j_{\BA-S_F^+}$.
We have isomorphisms
\begin{eqnarray}
D\Phi_{F}^+(\CK)\cong DR\Gamma(\BA,S_{F}^+;\CK)
\cong R\Gamma_c(\BA; j_{*}j^*D\CK)\cong\\ \nonumber
\cong R\Gamma(\BA,\BA_{\geq r''};j_{*}j^*D\CK)\cong
R\Gamma(\BA,Y^+\cup \BA_{\geq r''};j_{*}j^*D\CK) \nonumber
\end{eqnarray}
Consider the restriction map
\begin{equation}
res: R\Gamma(\BA,Y^+\cup \BA_{\geq r''};j_*j^*D\CK)\lra
R\Gamma(\BA_{\leq r'},Y^+\cap \BA_{\leq r'};D\CK)
\end{equation}
$\Cone (res)$ is isomorphic to
\begin{eqnarray}
R\Gamma(\BA,\BA_{\leq r'}\cup\BA_{\geq r''}\cup Y^+;j_*j^*D\CK)=
R\Gamma_c(\BA_{< r''},\BA_{\leq r'}\cup Y^+;j_*j^*D\CK)\cong\nonumber\\
\cong DR\Gamma(\BA_{< r''}-(\BA_{\leq r'}\cup Y^+);j_!j^*\CK)=
R\Gamma(\BA_{(r',r'')}-Y^+,S_F^+;\CK)\nonumber
\end{eqnarray}

\subsubsection{}
\label{acycl+} {\bf Lemma.} (Cf. I.3.8.1.) {\em
$R\Gamma(\BA_{(r',r'')}-Y^+,S_F^+;\CK)=0$.}

{\bf Proof.} Let us define the following subspaces of $\BA$.

$A:=\{(t_j)|\mbox{ for all }j\ |t_j|<1;\ \mbox{ there exists }j: t_j\neq 0\}$;
$A_{\BR}^+:=A\cap\BA^+_{\BR}$. Note that $A_{\BR}^+\cap i\cdot A_{\BR}^+=
\emp$. Due to monodromicity, it is easy to see that
$$
R\Gamma(\BA_{(r',r'')}-Y^+,S_F^+;\CK)\cong
R\Gamma(A-i\cdot A_{\BR}^+,A_{\BR}^+;\CK).
$$
Therefore, it is enough to prove the following

\subsubsection{} {\bf Claim.} {\em The restriction map
\begin{equation}
\label{resmain}
R\Gamma(A-i\cdot A_{\BR}^+;\CK)\lra R\Gamma(A_{\BR}^+;\CK)
\end{equation}
is an isomorphism.}

{\bf Proof of the Claim.}
Let us introduce for each $k\in J$ open subspaces
$$
A_k=\{(t_j)\in A|t_k\not\in i\cdot\BR_{\geq 0}\}\subset A-i\cdot A^+_{\BR}
$$
and
$$
A'_k=\{(t_j)\in A_{\BR}^+|t_k>0\}\subset A_{\BR}^+
$$
Obviously $A'_k=A_k\cap A_{\BR}^+$.
For each subset $M\subset J$ set $A_M:=\bigcap_{k\in M}A_k;\
A'_M:=\bigcap_{k\in M}A'_k$.

For each non-empty $M$ define the spaces
$B_M:=\{(t_j)_{j\in M}|\mbox{ for all }j\ |t_j|<1,
t_j\not\in i\cdot\BR_{\geq 0}\}$ and
$B'_M:=\{(t_j)_{j\in M}|\mbox{ for all }j\ t_j\in\BR,\ 0<t_j<1\}$.
We have obvious projections
$f_M:A_M\lra B_M,\ f'_M:A'_M\lra B'_M$.

Let us look at fibers of $f_M$ and $f'_M$.
Given $b=(t_j)_{j\in M}\in B_M$, the fiber
$f^{-1}_M(b)$ is by definition $\{(t_k)_{k\in J-M}||t_k|<1\}$,
the possible singularities
of our sheaf $\CK$ are at the hyperplanes $t_k=t_j$ and $t_k=0$.
It is easy to see that $f_M$ is "lisse" with respect to $\CK$ which
means in particular that a stalk $(f_{M*}\CK)_b$ is equal to
$R\Gamma(f^{-1}_M(b);\CK)$. The same considerations apply to $f'_M$.
Moreover, it follows from I.2.12 that the restriction maps
$$
R\Gamma(f^{-1}_M(b);\CK)\lra R\Gamma((f'_M)^{-1}(b);\CK)
$$
are isomorphisms for every $b\in B'_M$. This implies that
$f'_{M*}\CK$ is equal to the restriction of $f_{M*}\CK$ to $B'_M$.

The sheaf $f_{M*}\CK$ is smooth along the diagonal stratification.
For a small $\delta>0$ let $U_{\delta}\subset B_M$ denote
an open subset $\{(t_j)\in B_M|\ |\arg(t_j)|<\delta\mbox{ for all }j\}$.
The restriction maps $R\Gamma(B_M;f_{M*}\CK)\lra R\Gamma(U_{\delta};
f_{M*}\CK)$
are isomorphisms. We have $B'_M=\bigcap_{\delta}U_{\delta}$, therefore
by I.2.12 the restriction
$R\Gamma(B_M;f_{M*}\CK)\lra R\Gamma(B_M';f_{M*}\CK)$ is an isomorphism.
This implies,
by Leray, that restriction maps
$$
R\Gamma(A_M;\CK)\lra R\Gamma(A'_M;\CK)
$$
are isomorphisms for every non-empty $M$.

Obviously $A-i\cdot A^+_{\BR}=\bigcup_{k\in J}A_k$ and
$A^+_{\BR}=\bigcup_{k\in J}A'_k$. Therefore,
by Mayer-Vietoris the map ~(\ref{resmain}) is an isomorphism.
This proves the claim, together with the lemma.
$\Box$

A clockwise rotation by $\pi/2$ induces an isomorphism
$$
R\Gamma(\BA_{\leq r'},Y^+\cap \BA_{\leq r'};D\CK)\cong
R\Gamma(\BA_{\leq r'},\epsilon\cdot D_F^{+opp};D\CK)\cong
R\Gamma(\BA_{\leq r'},\epsilon\cdot S^+_F;D\CK)
$$
and the last complex is isomorphic to $\Phi_F^+(D\CK)$ in view of
{}~\ref{reduct}. This proves the theorem for the case of the $0$-facet $F$.

\subsection{} Suppose that $F$ is an arbitrary positive facet.
{}From the description of positive facets (see {\em infra},
{}~\ref{listpos}) one sees that the cell $D^+_F$ is homeomorphic to a
cartesian product of the form
$$
D^+_{F_0}\times D^+_{F_1}\times\ldots\times D^+_{F_a}
$$
where $F_0$ is a $0$-facet of the principal
stratification in some affine space of smaller dimension, and
$D^+_{F_i}$ are the cells of {\em the diagonal} stratification
discussed in the previous section.

Using this remark, we apply a combination of the arguments of the previous
subsection (to the first factor) and of I.3.8 (to the remaining factors)
to get the required isomorphism. We leave details to the reader.

The theorem is proved. $\Box$

\subsection{} {\bf Theorem.} {\em All functors $\Phi^+_F$ are $t$-exact.
In other words,
for all positive facets $F$,
$$
\Phi^+_F(\CM(\BA;\CS))\subset\Vect\subset\CD^b(pt).
$$}

{\bf Proof.} The same as that of I.3.9. $\Box$.

\subsection{} Thus we get exact functors
\begin{equation}
\label{phi+}
\Phi^+:\CM(\CA;\CS)\lra\Vect
\end{equation}
commuting with Verdier duality.

We will also denote vector spaces $\Phi^+_F(\CM)$ by $\CM^+_F$.

\subsection{Canonical and variation maps} Suppose we have a positive facet $E$.
Let us denote by $^+\Fac^1(E)$ the set of all positive facets $F$ such that
$E<F$, $\dim\ F=\dim\ E+1$. We have
\begin{equation}
\label{unican}
S^+_E=\bigcup_{F\in\ ^+\Fac^1(E)}D^+_F
\end{equation}

Suppose we have $\CK\in\CD(\BA,\CS)$.

\subsubsection{}
\label{addcan} {\bf Lemma.} {\em We have a natural isomorphism
$$
R\Gamma(S^+_E,\bigcup_{F\in\ ^+\Fac^1(E)}S^+_F;\CK)\cong
\oplus_{F\in\ ^+\Fac^1(E)}\ R\Gamma(D^+_F,S^+_F;\CK)
$$}

{\bf Proof.} Note that $S^+_E-\bigcup_{F\in\ ^+\Fac^1(E)}S^+_F=
\bigcup_{F\in\ ^+\Fac^1(E)}\ \DO^+_F$ (disjoint union). The claim follows
now from the
Poincar\'{e} duality. $\Box$

Therefore, for any $F\in\ ^+\Fac^1(E)$ we get a natural inclusion map
\begin{equation}
\label{projcan}
i^F_E: R\Gamma(D^+_F,S^+_F;\CK)\hra
R\Gamma(S^+_E,\bigcup_{F'\in\Fac^1(E)}S^+_{F'};\CK)
\end{equation}

Let us define a map
$$
u^F_E(\CK):\Phi^+_F(\CK)\lra\Phi^+_E(\CK)
$$
as a composition
\begin{eqnarray}
R\Gamma(D^+_F,S^+_F;\CK)[-p]\overset{i^F_E}{\lra}
R\Gamma(S^+_E,\bigcup_{F'\in\ ^+\Fac^1(E)}S^+_{F'};\CK)[-p]\lra\\ \nonumber
R\Gamma(S^+_E;\CK)[-p]\lra
R\Gamma(D^+_E,S^+_E)[-p+1]\nonumber
\end{eqnarray}
where the last arrow is the coboundary map for the couple $(S^+_E,D^+_E)$,
and the second one is evident.

This way we get a natural transormation
\begin{equation}
\label{can}
^+u^F_E:\Phi^+_F\lra \Phi^+_E
\end{equation}
which will be called a {\em canonical map}.

We define a {\em variation morphism}
\begin{equation}
\label{var}
^+v^E_F:\Phi^+_E\lra\Phi^+_F
\end{equation}
as follows. By definition, $^+v^E_F(\CK)$ is the map dual to
the composition
$$
D\Phi^+_F(\CK)\iso\Phi^+_F(D\CK)\overset{^+u^F_E(D\CK)}{\lra}\Phi^+_E(D\CK)
\iso D\Phi^+_E(\CK).
$$

\subsection{Cochain complexes}
For each $r\in [0,N]$ and $\CM\in\CM(\BA,\CS)$
introduce vector spaces
\begin{equation}
^+C^{-r}(\BA;\CM)=\oplus_{F: F \mbox{ positive, }\dim\ F=r}\ \CM^+_F
\end{equation}
For $r>0$ or $r<-N$ set $^+C^r(\BA;\CM)=0$.

Define operators
$$
d:\ ^+C^{-r}(\BA;\CM)\lra\ ^+C^{-r+1}(\BA;\CM)
$$
having components $^+u^F_E$.

\subsubsection{}
\label{nilp+} {\bf Lemma.} {\em $d^2=0$}.

{\bf Proof.} The same as that of I.3.13.1. $\Box$

This way we get a complex $^+C^{\bullet}(\BA;\CM)$ lying in
degrees from $-N$ to $0$. It will be called the {\em
complex of positive cochains} of the sheaf $\CM$.

\subsection{Theorem}
\label{rgamma+} {\em (i) A functor
$$
\CM\mapsto \ ^+C^{\bullet}(\BA;\CM)
$$
is an exact functor from $\CM(\BA;\CS)$ to the category of complexes
of vector spaces.

(ii) We have a canonical natural isomorphism in
$\CD(\{ pt\})$
$$
^+C^{\bullet}(\BA;\CM)\iso R\Gamma(\BA;\CM)
$$}

{\bf Proof.} One sees easily that restriction maps
$$
R\Gamma(\BA;\CM)\lra R\Gamma(\BA_{\BR}^+;\CM)
$$
are isomorphisms. The rest of the argument is the same as in I.3.14. $\Box$

\subsection{} Let us consider a function $\sum_jt_j:\BA\lra\BA^1$,
and the corresponding vanishing cycles functor
\begin{equation}
\label{vanish}
\Phi_{\Sigma\ t_j}:\CD^b(\BA)\lra\CD^b(\BA_{(0)})
\end{equation}
where $\BA_{(0)}=\{(t_j)|\sum_j t_j=0\}$, cf. ~\cite{ks}, 8.6.2.

If $\CK\in\CD(\BA;\CS)$, it is easy to see that the complex
$\Phi_{\Sigma\ t_j}(\CK)$ has the support at the origin. Let us denote
by the same letter its stalk at the origin --- it is a complex
of vector spaces.

\subsubsection{} {\bf Lemma.} {\em We have a natural
isomorphism
\begin{equation}
\label{vanishin}
\Phi_{\Sigma\ t_j}(\CK)\iso\Phi^+_{\{0\}}(\CK)
\end{equation}}

{\bf Proof} is left to the reader. $\Box$

\subsection{}
\label{mainphi} Let us consider the setup of ~\ref{color}. Let us denote
by the same letter $\CS$ the stratification of $\CA_{\nu}$ whose strata
are subspaces $\pi(S)$, $S$ being a stratum of the stratification $\CS$
on $^{\pi}\BA$. This stratification will be called
{\em the principal stratification of $\CA_{\nu}$}.

The function $\Sigma\ t_j$ is obviously $\Sigma_{\pi}$-equivariant,
therefore it induces the function for which we will use the same notation,
\begin{equation}
\label{sum}
\Sigma\ t_j:\CA_{\nu}\lra\BA^1
\end{equation}
Again, it is easy to see that for $\CK\in\CD^b(\CA_{\nu};\CS)$ the complex
$\Phi_{\Sigma t_j}(\CK)$ has the support at the origin.
Let us denote by $\Phi_{\nu}(\CK)$ its stalk at the origin.

It is known that the functor of vanishing cycles is $t$-exact
with respect to the middle perversity; whence we get an exact functor
\begin{equation}
\label{phi0}
\Phi_{\nu}:\CM(\CA_{\nu};\CS)\lra\Vect
\end{equation}

\subsubsection{} {\bf Lemma.} {\em Suppose that $\CN$ is a $\Sigma_{\pi}$-
equivariant complex of sheaves over $^{\pi}\BA$ which belongs to
to $\CD^b(^{\pi}\BA;\CS)$ (after forgetting $\Sigma_{\pi}$-action).
We have a
natural isomorphism
\begin{equation}
\label{vanishinv}
\Phi_{\nu}((\pi_*\CN)^{\Sigma_{\pi},-})\iso(\Phi_{\Sigma\ t_j}(\CN))
^{\Sigma_{\pi},-}
\end{equation}}

{\bf Proof} follows from the proper base change for vanishing cycles
(see ~\cite{d}, 2.1.7.1) and
the exactness of the functor $(\cdot)^{\Sigma_{\pi},-}$.
We leave details to the reader. $\Box$

\subsubsection{} {\bf Corollary.} {\em For a $\Sigma_{\pi}$-equivariant
sheaf $\CN\in\CM(^{\pi}\BA_{\nu};\CS)$ we have
a natural isomorphism of vector spaces
\begin{equation}
\label{vanishface}
\Phi_{\nu}(\CM)\iso(\Phi^+_{\{0\}}(\CN))^{\Sigma_{\pi},-}\
\end{equation}
where $\CM=(\pi_*\CN)^{\Sigma_{\pi},-}$. $\Box$}

\section{Standard sheaves}
\label{standsheaves}

Let us keep assumptions and notations of ~\ref{sym}.

\subsection{}
\label{group}
Let us denote by $^{\pi}\BA$ the complex affine space with coordinates $t_j,\
j\in J$. We will consider its principal stratification as in ~\ref{j}.

Suppose we are given a positive chamber $C$ and a point
$\bx=(x_j)_{j\in J}\in C$.
There extists a unique bijection
\begin{equation}
\label{bij}
\sigma_C:J\iso [N]
\end{equation}
such that for any $i,j\in J$, $\sigma_C(i)<\sigma_C(j)$ iff
$x_i<x_j$. This bijection does not depend upon the choice
of $\bx$.
This way we can identify the set of all positive chambers
with the set of all isomorphisms ~(\ref{bij}),
or, in other words, with the set of all total orderings of $J$.

Given $C$ and $\bx$ as above, suppose that we have $i,j\in J$ such that
$x_i<x_j$ and there is no $k\in J$ such that $x_i<x_k<x_j$. We will say that
$i,j$ are {\em neighbours} in $C$, more precisely that $i$ is a left neighbour
of $j$.

Let $\bx'=(x'_j)$ be a point with $x'_p=x_p$ for all $p\neq j$, and $x'_j$
equal to some number smaller than $x_i$ but greater than any $x_k$ such that
$x_k<x_i$.
Let $^{ji}C$ denote the chamber containing $\bx'$. Let us introduce a  homotopy
class of paths $^C\gamma_{ij}$ connecting $\bx$ and $\bx'$ as shown on
Fig. 5(a) below.

On the other hand, suppose that $i$ and $0$ are neighbours in $C$, there
is no $x_j$ between $0$ and $x_i$. Then we introduce the homotopy class
of paths from $\bx$ to itself as shown on Fig. 5 (b).

\begin{picture}(20,8)(-10,-4)

\put(-5,0){\circle*{0.2}}
\put(-5,-0.5){$i$}

\put(-9,0){\line(1,0){8}}

\put(-3.5,0){\circle*{0.2}}
\put(-3.5,-0.5){$j$}

\put(-5,0){\oval(3,3)[t]}
\put(-5,1.5){\vector(-1,0){0.5}}
\put(-5.5,2){$^C\gamma_{ij}$}

\put(-6.5,0){\circle*{0.2}}
\put(-6.5,-0.5){$j'$}

\put(-2,0){\circle*{0.2}}
\put(-8,0){\circle{0.2}}
\put(-8,-0.5){$0$}

\put(-5.5,-3){$(a)$}


\put(5,0){\circle*{0.2}}
\put(5.3,-0.5){$i$}

\put(1,0){\line(1,0){8}}

\put(3.5,0){\circle{0.2}}
\put(3.5,-0.5){$0$}
\put(3.5,0){\oval(3,3)}
\put(3.5,1.5){\vector(-1,0){0.5}}
\put(3.4,2){$^C\gamma_{i0}$}

\put(6.5,0){\circle*{0.2}}

\put(8,0){\circle*{0.2}}

\put(4.5,-3){$(b)$}


\put(-0.5,-4){Fig. 5}

\end{picture}

All chambers are contractible. Let us denote by $^{\pi}\BAO^+_{\BR}$ the union
of all positive chambers.
We can apply the discussion I.4.1
and consider the groupoid $\pi_1(^{\pi}\BAO,\ ^{\pi}\BAO^+_{\BR})$.
It has as the set of objects the set
of all positive chambers. The set of morphisms is generated by all morphisms
$^C\gamma_{ij}$ and $^C\gamma_{i0}$
subject to certain evident braiding relations. We will need only the following
particular case.

To define a {\em one-dimensional} local system $\CL$ over $^{\pi}\BAO$ is
the same as to give a set of one-dimensional vector spaces
$\CL_C,\ C\in\pi_0(^{\pi}\BAO^+_{\BR})$,
together with arbitrary invertible linear operators
\begin{equation}
\label{halfmon}
^CT_{ij}:\CL_C\lra\CL_{^{ji}C}
\end{equation}
("half-monodromies")  defined for chambers where $i$ is a left neighbour of
$j$ and
\begin{equation}
\label{mon}
^CT_{i0}:\CL_C\lra\CL_C
\end{equation}
defined for chambers with neighbouring $i$ and $0$.

\subsection{} Let us fix a weight $\Lambda\in X$.
We define a one-dimensional local system $\CI(^{\pi}\Lambda)$ over
$^{\pi}\BAO$ as follows. Its fibers $\CI(^{\pi}\Lambda)_{C}$ are
one-dimensional
linear spaces with fixed basis vectors; they will be identified with $B$.

Monodromies are defined as
$$
^CT_{ij}=\zeta^{i\cdot j},\
^CT_{i0}=\zeta^{-2\langle \Lambda,\pi(i)\rangle},
$$
for $i,j\in J$.

\subsection{} Let $j:\ ^{\pi}\BAO\lra\ ^{\pi}\BA$ denote an open embedding.
We will study the
following objects of $\CM(^{\pi}\BA;\CS)$:
$$
\CI(^{\pi}\Lambda)_{?}=j_?\CI(^{\pi}\Lambda)[N],
$$
where $?=!,*$. We have a canonical map
\begin{equation}
\label{!*}
m: \CI(^{\pi}\Lambda)_!\lra\CI(^{\pi}\Lambda)_*
\end{equation}
and by definition $\CI(^{\pi}\Lambda)_{!*}$ is its image, cf. I.4.5.

\vspace{1cm}
{\em COMPUTATIONS FOR $\CI(^{\pi}\Lambda)_!$}
\vspace{0.8cm}

\subsection{}
\label{listpos} We will use the notations ~\ref{gh} with $I=J$ and $n=1$.
For each $r\in [0,N]$, let us assign to a map $\varrho\in\CP_r(J;1)$ a point
$w_{\varrho}=(\varrho(j))_{j\in J}\in\ ^{\pi}\BA_{\BR}$. It is easy to see that
there exists a unique positive facet $F_{\varrho}$ containing $w_{\varrho}$,
and the rule
\begin{equation}
\label{rhofac}
\varrho\mapsto\ F_{\varrho}
\end{equation}
establishes an isomorphism between $\CP_r(J;1)$ and the set of all
positive facets of $\CS_{\BR}$. Note that $\CP_0(J;1)$ consists
of one element --- the unique map $J\lra [0]$; our stratification
has one zero-dimensional facet.

At the same time we have picked a point $^{F_{\varrho}}w:=w_{\varrho}$ on
each positive facet $F_{\varrho}$; this defines cells $D^+_F,\ S^+_F$
(cf. the last remark in ~\ref{poscells}).

\subsection{}
\label{listcouples} Given $\varrho\in\CP_r(J;1)$ and $\tau\in\CP_N(J;1)$,
it is easy to see that the chamber $C=F_{\tau}$ belongs to
$\Ch(F_{\varrho})$ iff $\tau$ is a refinement of $\varrho$ in the sense of
{}~\ref{refine}. This defines a bijection between the set of all refinements
of $\varrho$ and the set of all positive chambers containing $F_{\varrho}$.

We will denote the last set by $\Ch^+(F_{\varrho})$.

\subsubsection{Orientations} Let $F=F_{\varrho}$ be a positive facet and
$C=F_{\tau}\in\Ch^+(F)$. The map $\tau$ defines an isomorphism denoted
by the same letter
\begin{equation}
\label{tau}
\tau:J\iso [N]
\end{equation}
Using $\tau$, the natural order on $[N]$ induces a total order on $J$.
For each $a\in [r]$, let $m_a$ denote the minimal element
of $\varrho^{-1}(a)$, and set
$$
J'=J_{\varrho\leq \tau}:=J-\{ m_1,\ldots, m_r\}\subset J
$$
Let us consider the map
$$
(x_j)\in D_{F<C}\mapsto \{ x_j-m_{\varrho(j)}|j\in J'\}\in\BR^{J'}.
$$
It is easy to see that this mapping establishes an isomorphism
of the germ of the cell $D_{F<C}$ near the point $^Fw$ onto a germ
of the cone
$$
\{ 0\leq u_1\leq\ldots\leq u_{N-r}\}
$$
in $\BR^{J'}$ where we have denoted for a moment by $(u_i)$ the coordinates
in $\BR^{J'}$ ordered by the order induced from $J$.

This isomorphism together with the above order defines an orientation
on $D_{F<C}$.

\subsection{Basis in $\Phi_F^+(\CI(^{\pi}\Lambda)_!)^*$}
\label{basis!} We follow the pattern of
I.4.7. Let $F$ be a positive facet of dimension $r$. By definition,
$$
\Phi^+_F(\CI(^{\pi}\Lambda)_!)=H^{-r}(D^+_F,S^+_F;\CI(^{\pi}\Lambda)_!)=
H^{N-r}(D^+_F,S^+_F;j_!\CI(^{\pi}\Lambda))\cong
H^{N-r}(D^+_F,S^+_F\cup (_{\CH}H_{\BR}\cap D^+_F);j_!\CI(^{\pi}\Lambda)).
$$
Note that we have
$$
D^+_F-(S^+_F\cup\ _{\CH}H_{\BR})=\bigcup_{C\in\Ch^+(F)}\DO_{F<C}
$$
(disjoint union), therefore by additivity
$$
H^{N-r}(D^+_F,S^+_F\cup (_{\CH}H_{\BR}\cap D^+_F);j_!\CI(^{\pi}\Lambda))\cong
\oplus_{C\in\Ch^+(F)}
H^{N-r}(D_{F<C},S_{F<C};j_!\CI(^{\pi}\Lambda)).
$$
By Poincar\'{e} duality,
$$
H^{N-r}(D_{F<C},S_{F<C};j_!\CI(^{\pi}\Lambda))^*\cong
H^0(\DO_{F<C};\CI(^{\pi}\Lambda)^{-1})
$$
--- here we have used the orientations of cells $D_{F<C}$ introduced
above. By definition of the local system $\CI$, the last space is
canonically identified with $B$.

If $F=F_{\varrho},\ C=F_{\tau}$, we will denote by $c_{\varrho\leq \tau}
\in\Phi^+_F(\CI(^{\pi}\Lambda)_!)^*$
the image of $1\in H^0(\DO_{F<C};\CI(^{\pi}\Lambda)^{-1})$. Thus the
chains $c_{\varrho\leq\tau},\ \tau\in\Ord(\varrho)$, form
a basis of the space $\Phi^+_F(\CI(^{\pi}\Lambda)_!)^*$.

\subsection{Diagrams} It is convenient to use the following diagram
notations for chains $c_{\varrho\leq\tau}$.

Let us denote elements of $J$ by letters $a,b,c,\ldots$. An $r$-dimensional
chain $c_{\varrho\leq\tau}$ where $\varrho:J\lra [0,r]$, is represented
by a picture:

\begin{picture}(20,6)(-10,-3)

\put(-8,0){\line(1,0){10}}
\put(3,0){$\ldots$}
\put(4,0){\line(1,0){4}}

\put(-6,0){\circle{0.2}}
\put(-6.1,0.4){$0$}
\put(-6,0){\vector(1,0){0.8}}
\put(-5.2,0){\circle*{0.1}}
\put(-5.3,-0.4){$a$}

\put(-3,0){\circle{0.2}}
\put(-3.1,0.4){$1$}
\put(-3,0){\vector(1,0){0.8}}
\put(-2.2,0){\circle*{0.1}}
\put(-2.3,-0.4){$b$}
\put(-3,0){\vector(1,0){1.6}}
\put(-1.4,0){\circle*{0.1}}
\put(-1.5,-0.4){$c$}

\put(5,0){\circle{0.2}}
\put(4.9,0.4){$r$}
\put(5,0){\vector(1,0){0.8}}
\put(5.8,0){\circle*{0.1}}
\put(5.7,-0.4){$d$}
\put(5,0){\vector(1,0){1.6}}
\put(6.6,0){\circle*{0.1}}
\put(6.5,-0.4){$e$}

\put(-1,-3){Fig. 6. Chain $c_{\varrho\leq\tau}$.}

\end{picture}

A picture consists of $r+1$ fragments:

\begin{picture}(20,6)(-10,-3)

\put(-4,0){\line(1,0){8}}

\put(-1,0){\circle{0.2}}
\put(-1.1,0.4){$i$}
\put(-1,0){\vector(1,0){0.8}}
\put(-0.2,0){\circle*{0.1}}
\put(-0.3,-0.4){$a$}
\put(-1,0){\vector(1,0){1.6}}
\put(0.6,0){\circle*{0.1}}
\put(0.5,-0.4){$b$}
\put(-1,0){\vector(1,0){2.4}}
\put(1.4,0){\circle*{0.1}}
\put(1.3,-0.4){$c$}

\put(-1,-3){Fig. 7. Set $\varrho^{-1}(i)$.}

\end{picture}

where $i=0,\ldots, r$, the $i$-th fragment being a blank circle with
a number of small vectors going from it. These vectors are
in one-to-one correspondence with the set $\varrho^{-1}(i)$; their ends
are labeled by elements of this set. Their order
(from left to right) is determined by the order on $\varrho^{-1}(i)$
induced by $\tau$. The point $0$ may have no vectors (when $\varrho^{-1}(0)=
\emp$); all other points have at least one vector.

\subsection{} Suppose we have $\varrho\in\CP_r(J;1),\
\varrho'\in\CP_{r+1}(J;1)$. It is easy to see that $F_{\varrho}<
F_{\varrho'}$ if and only if there exists $i\in [0,r]$ such that
$\varrho=\delta_i\circ\varrho'$ where
$\delta_i:[0,r+1]\lra [0,r]$ carries $a$ to $a$ if $a\leq i$ and to $a-1$
if $a\geq i+1$. We will write in this case that $\varrho<\varrho'$.

Let us compute the dual to the canonical map
\begin{equation}
\label{can!}
^+u^*: \Phi^+_{F_{\varrho}}(\CI(^{\pi}\Lambda)_!)^*\lra
\Phi^+_{F_{\varrho'}}(\CI(^{\pi}\Lambda)_!)^*.
\end{equation}
Suppose we have $\tau\in\Ord(\varrho')$; set $C=F_{\tau}$, thus
$F_{\varrho}<F_{\varrho'}<C$. Let us define the sign
\begin{equation}
\label{signum}
\sgn (\varrho<\varrho'\leq\tau)=
(-1)^{\sum_{j=i+1}^{r+1}\card((\varrho')^{-1}(j)-1)}
\end{equation}
This sign has the following geometrical meaning. The cell
$D_{F_{\varrho'<C}}$ lies in the boundary of $D_{F_{\varrho<C}}$.
We have oriented these cells above. Let us define the compatibility
of these orientations as follows.
Complete an orienting basis of the smaller cell
by a vector directed outside the larger cell --- if we get the orientation
of the larger cell, we say that the orientations are compatible, cf. I.4.6.1.

It is easy to see from the definitions that the sign ~(\ref{signum}) is
equal $+1$ iff the orientations of $D_{F_{\varrho'<C}}$ and
$D_{F_{\varrho<C}}$ are compatible. As a consequence, we get

\subsubsection{} {\bf Lemma.}
\label{matrixu} {\em The map ~(\ref{can!}) has the following
matrix:
$$
^+u^*(c_{\varrho\leq\tau})=\sum\sgn(\varrho<\varrho'\leq\tau)c_{\varrho'\leq
\tau},
$$
the summation over all $\varrho'$ such that $F_{\varrho}<F_{\varrho'}<
F_{\tau}$ and $\dim F_{\varrho'}=\dim F_{\varrho}+1$.} $\Box$

\subsection{Isomorphisms $^{\pi}\phi_{!}$} We will use notations
of ~\ref{refine}
with $I$ replaced by $J$, $\fF$ by $^{\pi}\fF$, with $n=1$ and
$\Lambda_0=\ ^{\pi}\Lambda$.

Thus, for any $r\in [0,N]$ the set
$\{ b_{\varrho\leq\tau}|\varrho\in\CP_r(J;1),\
\tau\in\Ord(\varrho)\}$ is a basis of
$^+C^{-r}(^{\pi}\BA;\CI(^{\pi}\Lambda)_!)$.

Let us pick $\eta\in\CP_N(J)$. Any $\tau\in\CP_N(J;1)$ induces the bijection
$\tau':J\iso [N]$. We will denote by $\sgn(\tau,\eta)=\pm 1$ the
sign of the permutation $\tau'\eta^{-1}$.

Let us denote by
$\{ b_{\varrho\leq\tau}|\tau\in\Ord(\varrho)\}$
the basis in
$\Phi_{F_{\varrho}}^+(\CI(^{\pi}\Lambda)_!)$ dual to
$\{ c_{\varrho\leq\tau}\}$.

Let us define isomorphisms
\begin{equation}
\label{phirho}
^{\pi}\phi_{\varrho,!}^{(\eta)}:\Phi^+_{F_{\varrho}}(\CI(^{\pi}\Lambda)_!)\iso\
_{\varrho}C^{-r}_{^{\pi}\fF}(V(^{\pi}\Lambda))
\end{equation}
by the formula
\begin{equation}
\label{phirhoform}
^{\pi}\phi_{\varrho,!}^{(\eta)}(b_{\varrho\leq\tau})=\sgn(\tau,\eta)
\sgn(\varrho)\theta_{\varrho\leq\tau}
\end{equation}
(see ~(\ref{monom})) where
\begin{equation}
\label{signrho}
\sgn(\varrho)=(-1)^{\sum_{i=1}^r(r-i+1)\cdot(\card(\varrho^{-1}(i))-1)}
\end{equation}
for $\varrho\in\CP_r(J;1)$. Taking the direct sum of $^{\pi}\phi_{\varrho,!}
^{(\eta)},\
\varrho\in\CP_r(J;1)$, we get isomorphisms
\begin{equation}
\label{phir}
^{\pi}\phi_{r,!}^{(\eta)}:\ ^+C^{-r}(^{\pi}\BA;\CI(^{\pi}\Lambda)_!)\iso\
_{\chi_J}C^{-r}_{^{\pi}\fF}(V(^{\pi}\Lambda))
\end{equation}
A direct computation using ~\ref{matrixu} shows that the maps
$^{\pi}\phi_{r,!}^{(\eta)}$ are
compatible with differentials. Therefore, we arrive at

\subsection{}
\label{phi} {\bf Theorem.} {\em The maps $^{\pi}\phi_{r,!}^{(\eta)}$ induce
an isomorphism
of complexes
\begin{equation}
\label{phieq}
^{\pi}\phi_{!}^{(\eta)}:\ ^+C^{\bullet}(^{\pi}\BA;\CI(^{\pi}\Lambda)_!)\iso\
_{\chi_J}C^{\bullet}_{^{\pi}\fF}(V(^{\pi}\Lambda))\ \Box
\end{equation}}

\vspace{1cm}
{\em COMPUTATIONS FOR $\CI(^{\pi}\Lambda)_*$}
\vspace{0.8cm}

\subsection{Bases} The Verdier dual to $\CI(^{\pi}\Lambda)_*$ is canonically
isomorphic to $\CI(^{\pi}\Lambda)^{-1}_!$. Therefore, by Theorem
{}~\ref{dual+} for
each positive facet $F$ we
have natural isomorphisms
\begin{equation}
\label{duality}
\Phi^+_F(\CI(^{\pi}\Lambda)_*)^*\cong \Phi^+_F(\CI(^{\pi}\Lambda)^{-1}_!)
\end{equation}
Let $\{ \tc_{\varrho\leq\tau}|\tau\in\Ord(\varrho)\}$ be
the basis in $\Phi^+_F(\CI(^{\pi}\Lambda)^{-1}_!)^*$ defined in ~\ref{basis!},
with $\CI(^{\pi}\Lambda)$ replaced by $\CI(^{\pi}\Lambda)^{-1}$.
We will denote by
$\{ c_{\varrho\leq\tau,*}|\tau\in\Ord(\varrho)\}$ the dual basis
in $\Phi^+_F(\CI(^{\pi}\Lambda)_*)^*$. Finally, we will denote by
$\{ b_{\varrho\leq\tau,*}|\tau\in\Ord(\varrho)\}$ the basis in
$\Phi^+_F(\CI(^{\pi}\Lambda)_*)$ dual to the previous one.

Our aim in the next subsections will be the description
of canonical morphisms $m:\Phi^+_F(\CI(^{\pi}\Lambda)_!)
\lra\Phi^+_F(\CI(^{\pi}\Lambda)_*)$ and
of the cochain complex $^+C^{\bullet}(^{\pi}\BA;\CI(^{\pi}\Lambda)_*)$ in
terms of our bases.

\subsection{Example.} Let us pick an element, $i\in I$ and let
$\pi:J:=\{ i\}\hra I$. Then we are in a one-dimensional situation,
cf. ~\ref{onedim}. The space $^{\pi}\BA$ has one coordinate $t_i$.
By definition, the local system $\CI(^{\pi}\Lambda)$ over
$^{\pi}\BAO=\ ^{\pi}\BA-\{ 0\}$ with the base point
$w: t_i=1$ has the fiber $\CI_w=B$ and monodromy equal to $\zeta^{-2\langle
\Lambda,i\rangle}$ along a counterclockwise loop.

The stratification $\CS_{\BR}$ has a unique $0$-dimensional
facet $F=F_{\varrho}$ corresponding to the unique element $\varrho
\in\CP_0(J;1)$,
and a unique $1$-dimensional positive facet $E=F_{\tau}$ corresponding
to the unique element $\tau\in\CP_1(J;1)$.

Let us construct the dual chain $c^*_1:=c_{\varrho\leq\tau,*}$. We adopt the
notations of ~\ref{onedim}, in particular of Fig.1.
The chain
$\tc_1:=\tc_{\varrho\leq\tau}\in
H^1(^{\pi}\BA,\{ 0,w\};j_!\CI(^{\pi}\Lambda)^{-1})^*$
is shown on Fig. 8(a) below. As a first step, we define the dual chain
$c^0_1\in H^1(^{\pi}\BA-\{ 0,w\},^{\pi}\BA_{\geq r''};\CI(^{\pi}\Lambda))^*$
--- it is also
shown on Fig. 8(a).

\begin{picture}(20,9)(-10,-4.5)

\put(-5,0){\circle{0.2}}
\put(-5.5,-0.4){$F$}
\put(-5,0){\oval(6,6)}
\put(-5.5,3.3){$S_{r''}$}

\put(-4.9,0){\line(1,0){4}}
\put(-1.5,0.3){$E$}

\put(-4.9,0){\vector(1,0){1}}
\put(-4,0){\circle*{0.1}}
\put(-4.7,-0.4){$\tc_1$}

\put(-4,-3){\line(0,1){6}}
\put(-4,0){\vector(0,1){1}}
\put(-3.8,2){$c^0_1$}

\put(-3,0){\circle*{0.2}}
\put(-3.1,-0.5){$w$}

\put(-5,0){\oval(3,3)}
\put(-6.5,1.6){$S_{r'}$}

\put(-5,-4){(a)}


\put(5,0){\circle{0.2}}
\put(5,0){\oval(6,6)}
\put(4.5,3.3){$S_{r'}$}

\put(5.1,0){\line(1,0){4}}

\put(5.1,0){\vector(1,0){1}}
\put(5.3,-0.4){$\epsilon\tc_1$}

\put(5,0){\oval(3,3)}
\put(5.5,1.5){\vector(-1,0){0.5}}
\put(5.5,1.7){$c^*_1$}

\put(6.5,0){\circle*{0.2}}
\put(6.6,-0.5){$\epsilon w$}

\put(5,-4){(b)}


\put(-0.5,-4.5){Fig. 8.}
\end{picture}

Next, we make a clockwise rotation of $c_1^0$ on $\pi/2$, and make a homotopy
inside the disk $^{\pi}\BA_{\leq r'}$ to a chain $c^*_1\in
H^1(^{\pi}\BA_{\leq r'}-\{ 0\},\epsilon\cdot w;\CI(^{\pi}\Lambda))^*$ beginning
and
ending at $\epsilon\cdot w$. This chain is shown on Fig. 8(b) (in a bigger
scale). Modulo "homothety" identification
$$
H^1(^{\pi}\BA_{\leq r'}-\{ 0\},\epsilon\cdot w;\CI(^{\pi}\Lambda))^*\cong
H^1(^{\pi}\BA-\{ 0\},w;\CI(^{\pi}\Lambda))^*=\Phi_F(\CI(^{\pi}\Lambda)_*)^*
$$
this chain coincides with $c_{\varrho\leq\tau,*}$.

The canonical map
$$
m^*: \Phi_F^+(\CI(^{\pi}\Lambda)_*)^*\lra\Phi_F^+(\CI(^{\pi}\Lambda)_!)^*
$$
carries $c_{\varrho\leq\tau,*}$ to
$[\langle\Lambda,i\rangle]_{\zeta}\cdot c_{\varrho\leq\tau}$. The
boundary map
$$
d^*: \Phi^+_F(\CI(^{\pi}\Lambda)_*)^*\lra\Phi^+_E(\CI(^{\pi}\Lambda)_*)^*
$$
carries $c_{\varrho\leq\tau,*}$ to
$[\langle\Lambda,i\rangle]_{\zeta}\cdot c_{\tau\leq\tau,*}$.

\subsection{Vanishing cycles at the origin}
\label{orig} Let us return
to the general situation. First let us consider an important case
of the unique $0$-dimensional facet --- the origin. Let $0$ denote the
unique element of $\CP_0(J;1)$. To shorten the notation,
let us denote $\Phi^+_{F_0}$ by $\Phi^+_0$.

The bases in $\Phi^+_0(\CI(^{\pi}\Lambda)_!)$,
etc. are
numbered by all bijections $\tau:J\iso [N]$.
Let us pick such a bijection. Let $\tau^{-1}(i)=\{j_i\},\ i=1,\ldots, N$.
The chain $c_{0\leq\tau}\in\Phi_{0}^+(\CI(^{\pi}\Lambda)_!)^*$
is depicted as follows:

\begin{picture}(20,6)(-10,-3)

\put(-4,0){\line(1,0){4}}

\put(-3,0){\circle{0.2}}
\put(-3.1,0.4){$0$}
\put(-3,0){\vector(1,0){0.8}}
\put(-2.2,0){\circle*{0.1}}
\put(-2.3,-0.4){$j_1$}
\put(-3,0){\vector(1,0){1.6}}
\put(-1.4,0){\circle*{0.1}}
\put(-1.5,-0.4){$j_2$}
\put(-3,0){\vector(1,0){2.4}}
\put(-0.6,0){\circle*{0.1}}
\put(-0.7,-0.4){$j_3$}

\put(0.5,0){$\ldots$}
\put(1.5,0){\line(1,0){2}}
\put(2,0){\vector(1,0){0.5}}
\put(2.5,0){\circle*{0.1}}
\put(2.4,-0.4){$j_N$}

\put(-1,-3){Fig. 9. $c_{0\leq\tau}$}

\end{picture}

Using considerations completely analogous to the one-dimensional case above,
we see that the dual chain
$c_{0\leq\tau,*}\in\Phi_{0}^+(\CI(^{\pi}\Lambda)_*)^*$ is portayed as
follows:

\begin{picture}(20,6)(-10,-3)

\put(-10,0){\line(1,0){10}}

\put(-3,0){\circle{0.2}}
\put(-3.1,0.2){$0$}

\put(-3,0){\oval(1.6,1.6)}
\put(-3,0.8){\vector(-1,0){0.3}}
\put(-2.2,0){\circle{0.2}}
\put(-2.1,-0.4){$j_1$}

\put(-3,0){\oval(3.2,2.4)}
\put(-3,1.2){\vector(-1,0){0.8}}
\put(-1.4,0){\circle{0.2}}
\put(-1.3,-0.4){$j_2$}

\put(-3,0){\oval(4.8,3.2)}
\put(-3,1.6){\vector(-1,0){1.3}}
\put(-0.6,0){\circle{0.2}}
\put(-0.5,-0.4){$j_3$}

\put(0.5,0){$\ldots$}
\put(1.5,0){\line(1,0){2}}

\put(-3,0){\oval(10,4)}
\put(-3,2){\vector(-1,0){1.8}}
\put(2,0){\circle{0.2}}
\put(2.1,-0.4){$j_N$}

\put(-1,-3){Fig. 10. $c_{0\leq\tau,*}$}

\end{picture}

The points $t_{j_i}$ are travelling independently along the corresponding
loops, in the indicated directions. The section of $\CI(^{\pi}\Lambda)^{-1}$
over
this cell is determined by the requierement to be equal to $1$ when
all points are equal to marking points, at the end of their the travel
(coming from below).

\subsection{Isomorphisms $\phi_{0,*}^{(\eta)}$} We will use notations of
{}~\ref{refine} and
{}~\ref{dualbases} with $I$ replaced by $J$, $\fF$ --- by $^{\pi}\fF$,
with $n=1$ and
$\Lambda_0=\ ^{\pi}\Lambda$.
By definition, $C^0_{^{\pi}\fF}(V(^{\pi}\Lambda))=V(^{\pi}\Lambda)$.
The space
$V(^{\pi}\Lambda)_{\chi_J}$
admits as a basis (of cardinality $N!$) the set consisting of all monomials
$$
\theta_{0\leq\tau}=\theta_{\tau^{-1}(N)}\theta_{\tau^{-1}(N-1)}\cdot
\ldots\cdot\theta_{\tau^{-1}(1)}v_{^{\pi}\Lambda},
$$
where $\tau$ ranges through the set of all bijections $J\iso [N]$.
By definition, $\{\theta^*_{0\leq\tau}\}$ is the dual basis of
$V(^{\pi}\Lambda)^*_{\chi_J}$.

Let us pick $\eta\in\CP_N(J)$.
Let us define an isomorphism
\begin{equation}
\label{piphio*}
^{\pi}\phi_{0,*}^{(\eta)}: \Phi^+_0(\CI(^{\pi}\Lambda)_*)\iso\
V(^{\pi}\Lambda)^*_{\chi_J}
\end{equation}
by the formula
\begin{equation}
\label{piphioform*}
^{\pi}\phi_{0,*}^{(\eta)}(b_{0\leq\tau})=\sgn(\tau,\eta)\theta^*_{0\leq\tau}
\end{equation}

\subsection{}
\label{shapov} {\bf Theorem.} {\em The square
$$\begin{array}{rcccl}
\;&\Phi^+_0(\CI(^{\pi}\Lambda)_!)&\overset{^{\pi}\phi_{0,!}^{(\eta)}}{\iso}&
         V(^{\pi}\Lambda)_{\chi_J}&\;\\
\;&m\downarrow&\;&\downarrow S_{^{\pi}\Lambda}&\;\\
\;&\Phi^+_0(\CI(^{\pi}\Lambda)_*)&\overset{^{\pi}\phi_{0,*}^{(\eta)}}{\iso}&
      V(^{\pi}\Lambda)^*_{\chi_J}&\;
\end{array}$$
commutes.}

{\bf Proof.} This follows directly from the discussion of ~\ref{orig}
and the definition of the form $S_{^{\pi}\Lambda}$. $\Box$

\subsection{}
\label{letus} Let us pass to the setup of ~\ref{color} and ~\ref{mainphi}.
Let $j_{\nu}:\CAO_{\nu}\hra\CA_{\nu}$ denote the embedding of the open
stratum of the principal stratification.

Set
\begin{equation}
\label{inulambda}
\CI_{\nu}(\Lambda)=(\pi_*\CI(^{\pi}\Lambda))^{\Sigma_{\pi},-}
\end{equation}
It is a local system over $\CAO_{\nu}$.

Let us define objects
\begin{equation}
\label{inula?}
\CI_{\nu}(\Lambda)_?:=j_{\nu?}\CI_{\nu}(\Lambda)[N]\in\CM(\CA_{\nu};\CS)
\end{equation}
where $?=!,*$ or $!*$. These objects will be called
{\em standard sheaves over $\CA_{\nu}$}.

The same reasoning as in ~\ref{jcolor} proves

\subsubsection{} {\bf Lemma.} {\em We have natural isomorphisms
\begin{equation}
\label{jstand}
\CI_{\nu}(\Lambda)_?\cong(\pi_*\CI(^{\pi}\Lambda)_?)^{\Sigma_{\pi},-}
\end{equation}
for $?=!,*$ or $!*$. $\Box$}

\subsection{} For a given $\eta\in\CP_N(J)$ the isomorphisms
$^{\pi}\phi_{0,*}^{(\eta)}$ and
$^{\pi}\phi_{0,!}^{(\eta)}$ are skew $\Sigma_{\pi}$-
equivariant. Therefore, after passing to invariants in Theorem ~\ref{shapov}
we get

\subsection{}
\label{shaposym} {\bf Theorem.} {\em The maps $^{\pi}\phi_{0,*}^{(\eta)}$  and
$^{\pi}\phi_{0,!}^{(\eta)}$ induce isomorphisms included into a commutative
square
$$\begin{array}{rcccl}
\;&\Phi_{\nu}(\CI_{\nu}(\Lambda)_!)&
\overset{\phi_{\nu,\Lambda,!}^{(\eta)}}{\iso}&
         V(\Lambda)_{\nu}&\;\\
\;&m\downarrow&\;&\downarrow S_{\Lambda}&\;\\
\;&\Phi_{\nu}(\CI_{\nu}(\Lambda)_*)&
\overset{\phi_{\nu,\Lambda,*}^{(\eta)}}{\iso}&
         V(\Lambda)^*_{\nu}&\;
\end{array}$$
and
\begin{equation}
\label{lnu}
\phi_{\nu,\Lambda,!*}^{(\eta)}:\Phi_{\nu}(\CI_{\nu}(\Lambda)_{!*})\iso
L(\Lambda)_{\nu}
\end{equation}}

{\bf Proof} follows from the previous theorem and Lemma ~\ref{S-sym}(ii).
$\Box$

\subsection{} Now suppose we are given an arbitrary $r$, $\varrho\in
\CP_r(J;1)$ and $\tau\in\Ord(\varrho)$. The picture of the dual chain
$c_{\varrho\leq\tau,*}$ is a combination of Figures 10 and 3. For example,
the chain dual to the one on Fig. 6 is portayed as follows:

\begin{picture}(20,6)(-10,-3)

\put(-8,0){\line(1,0){8}}
\put(0.5,0){$\ldots$}
\put(1.5,0){\line(1,0){6.5}}

\put(-6,0){\circle{0.2}}
\put(-6.1,0.3){$0$}
\put(-6,0){\oval(1.6,1.6)}
\put(-6,0.8){\vector(-1,0){0.5}}
\put(-5.2,0){\circle*{0.1}}
\put(-5.1,-0.4){$a$}

\put(-3,0){\circle{0.2}}
\put(-3.1,0.4){$1$}
\put(-3,0){\oval(1.6,1.6)[b]}
\put(-3,-0.8){\vector(1,0){0.5}}
\put(-2.2,0){\circle*{0.1}}
\put(-2.1,-0.4){$b$}

\put(-3,0){\oval(3.2,3.2)[b]}
\put(-3,-1.6){\vector(1,0){0.5}}
\put(-1.4,0){\circle*{0.1}}
\put(-1.2,-0.4){$c$}

\put(5,0){\circle{0.2}}
\put(4.9,0.4){$r$}
\put(5,0){\oval(1.6,1.6)[b]}
\put(5,-0.8){\vector(1,0){0.5}}
\put(5.8,0){\circle*{0.1}}
\put(5.9,-0.4){$d$}

\put(5,0){\oval(3.2,3.2)[b]}
\put(5,-1.6){\vector(1,0){0.5}}
\put(6.6,0){\circle*{0.1}}
\put(6.7,-0.4){$e$}

\put(-1,-3){Fig. 11. Chain $c_{\varrho\leq\tau,*}$.}

\end{picture}

\subsection{Isomorphisms $^{\pi}\Phi_*$} Let us pick $\eta\in\CP_N(J)$.
Let us define isomorphisms
\begin{equation}
\label{phirho*}
^{\pi}\phi_{\varrho,*}^{(\eta)}:\Phi^+_{F_{\varrho}}(\CI(^{\pi}\Lambda)_*)\iso\
_{\varrho}C^{-r}_{^{\pi}\fF^*}(V(^{\pi}\Lambda)^*)
\end{equation}
by the formula
\begin{equation}
\label{phirhoform*}
^{\pi}\phi_{\varrho,*}^{(\eta)}(b_{\varrho\leq\tau,*})=\sgn(\tau,\eta)
\sgn(\varrho)\theta_{\varrho\leq\tau}^*
\end{equation}
where $\sgn(\varrho)$ is defined in ~(\ref{signrho}).

Taking the direct sum of $^{\pi}\phi_{\varrho,*}^{(\eta)},\
\varrho\in\CP_r(J;1)$, we get isomorphisms
\begin{equation}
\label{phir*}
^{\pi}\phi_{r,*}^{(\eta)}:\ ^+C^{-r}(^{\pi}\BA;\CI(^{\pi}\Lambda)_*)\iso\
_{\chi_J}C^{-r}_{^{\pi}\fF^*}(V(^{\pi}\Lambda)^*)
\end{equation}

\subsection{Theorem}
\label{phi*} {\em The maps $^{\pi}\phi_{r,*}^{(\eta)}$ induce
an isomorphism
of complexes
\begin{equation}
\label{phieq*}
^{\pi}\phi_{*}^{(\eta)}:\ ^+C^{\bullet}(^{\pi}\BA;\CI(^{\pi}\Lambda)_*)\iso\
_{\chi_J}C^{\bullet}_{^{\pi}\fF^*}(V(^{\pi}\Lambda)^*)
\end{equation}
which makes the square
$$\begin{array}{ccc}
^+C^{\bullet}(^{\pi}\BA;\CI(^{\pi}\Lambda)_!)&\overset{^{\pi}\phi_{!}^{(\eta)}}
\iso&\
_{\chi_J}C^{\bullet}_{^{\pi}\fF}(V(^{\pi}\Lambda))\\
m\downarrow&\;&\downarrow S\\
^+C^{\bullet}(^{\pi}\BA;\CI(^{\pi}\Lambda)_*)&\overset{^{\pi}\phi_{*}^{(\eta)}}
\iso&\
_{\chi_J}C^{\bullet}_{^{\pi}\fF^*}(V(^{\pi}\Lambda)^*)
\end{array}$$
commutative.}

{\bf Proof.} Compatibility with differentials is verified directly and
commutativity of the square are checked directly from the
geometric description of our chains (actually, it is sufficient to
check one of these claims --- the other one follows formally).

Note the geometric meaning of operators $t_i$ from ~(\ref{ti}) ---
they correspond to the deletion of the $i$-th loop on Fig. 10. $\Box$

\subsection{} Now let us pass to the situation ~\ref{letus}.
It follows from Theorem ~\ref{rgamma+} (after passing to skew $\Sigma_{\pi}$-
invariants) that the complexes $^+C^{\bullet}(^{\pi}\BA;\CI_{\nu?})
^{\Sigma_{\pi},-}$ compute the stalk of $\CI_{\nu?}$ at the origin.
Let us denote
this stalk by $\CI_{\nu?,0}$.

Therefore, passing to $\Sigma_{\pi}$-invariants in the previous theorem,
we get

\subsection{Theorem}
\label{stalks} {\em The isomorphisms $^{\pi}\phi^{(\eta)}_?$ where
$?=!,*$ or $!*$, induce isomorphisms in $\CD^b(pt)$ included into a
commutative square
$$\begin{array}{ccc}
\CI_{\nu}(\Lambda)_{!,0}&\overset{_{\Lambda}\phi^{(\eta)}_{\nu,!,0}}
{\iso}&_{\nu}C^{\bullet}_{\fF}(V(\Lambda))\\
m\downarrow&\;&\downarrow S\\
\CI_{\nu}(\Lambda)_{*,0}&\overset{_{\Lambda}\phi^{(\eta)}_{\nu,*,0}}
{\iso}&_{\nu}C^{\bullet}_{\fF^*}(V(\Lambda)^*)
\end{array}$$
and
\begin{equation}
\label{stalkirr}
_{\Lambda}\phi^{(\eta)}_{\nu,!*,0}:
\CI_{\nu}(\Lambda)_{!*,0}\iso\ _{\nu}C^{\bullet}_{\ff}(L(\Lambda))\ \ \Box
\end{equation}}

\newpage
\begin{center}
{\bf CHAPTER 3. Fusion.}
\end{center}
\vspace{1cm}

\section{Additivity theorem}

\subsection{} Let us start with the setup of ~\ref{j}. For a non-negative
integer $n$ let us denote by $(n)$ the set $[-n,0]$.
Let us introduce the following spaces.
$\BA^{(n)}$ - a complex affine space with a fixed system
of coordinates $(t_i),\ i\in (n)$.
Let $^nJ$ denote the disjoint union $(n)\cup J$;
$^n\BA$ --- a complex affine space with coordinates $t_j,\ j\in\ ^nJ$.

In general, for an affine space with a distinguished coordinate system of,
we will denote by $\CS_{\Delta}$ its diagonal stratification
as in ~\ref{diagsetup}.

Let $^n\BAO\subset\ ^n\BA$,
$\BAO\ ^{(n)}\subset\BA^{(n)}$ be the open strata of $\CS_{\Delta}$.

Let $^np:\ ^n\BA\lra\BA^{(n)}$ be the evident projection;
$^n\BB=\ ^np^{-1}(\BAO\ ^{(n)})$.
Given a point $\bz=(z_i)\in\BA^{(n)}$, let us denote by $^{\bz}\BA$ the
fiber $^np^{-1}(\bz)$ and by $^{\bz}\CS$ the stratification induced
by $\CS_{\Delta}$. We will consider $t_j,\ j\in J,$ as coordinates on
$^{\bz}\BA$.

The subscript $(.)_{\BR}$ will mean as usually "real points".

\subsection{}
Let us fix a point $\bz=(z_{-1},z_0)\in\BA^{(2)}_{\BR}$ such that
$z_{-1}<z_{0}$.
Let us concentrate on the fiber $^{\bz}\BA$. As an abstract variety
it is canonically isomorphic to $\BA$ --- a complex affine space
with coordinates $t_j,\ j\in J$; so we will suppress $\bz$ from its notation,
keeping it in the notation for the stratification $^{\bz}\CS$ where
the dependence on $\bz$ really takes place.

Let us fix a real $w>z_0$.
Let us pick two open non-intersecting disks $D_i\subset\BC$ with centra at
$z_i$ and not containing $w$. Let us pick two real numbers $w_{i}>z_i$
such that $w_i\in D_i$,
and paths $P_i$ connecting $w$ with $w_i$ as shown
on Fig. 12 below.

\begin{picture}(20,6)(-10,-3)

\put(-6,0){\circle{0.2}}
\put(-6.9,0){$z_{-1}$}


\put(-6,1){\oval(3,6)[b]}
\put(-6,0){\oval(2,2)}
\put(-6.4,-.6){$D_{-1}$}

\put(-.5,1){\oval(8,2)[t]}
\put(-.5,1){\oval(14,4)[t]}

\put(5,0){\circle*{0.2}}
\put(5.2,0){$w$}
\put(5,1){\oval(3,6)[b]}
\put(-6.3,1.5){$U_{-1}$}


\put(0.1,1){\oval(9.8,2.8)[t]}
\put(-5.5,1){\oval(1.4,2)[br]}
\put(5,0){\line(0,1){1}}
\put(-5.5,0){\circle*{0.2}}
\put(-5.9,.2){$w_{-1}$}
\put(.1,2.5){$P_{-1}$}


\put(0,0){\circle{0.2}}
\put(-.5,0.2){$z_{0}$}
\put(0,0){\oval(2,2)}
\put(-.2,-.6){$D_0$}


\put(.5,0){\circle*{.2}}
\put(0.2,.2){$w_0$}
\put(0.5,0){\line(1,0){4.5}}
\put(1.75,0.3){$P_0$}


\put(2.25,0){\oval(7.5,3)}
\put(1.75,-1){$U_0$}

\put(-1,-3){Fig. 12.}
\end{picture}

Let us denote by $\CQ_2(J)$ the set of all
maps $\rho:J\lra [-1,0]$. Given such a map, let us denote by $\BA_{\rho}
\subset\ \BA$ an open subvariety consisting of points
$(t_j)_{j\in J}$ such that $t_j\in D_{\rho(j)}$ for all $j$.
We will denote by the same letter $^{\bz}\CS$ the stratification of this
space induced by $^{\bz}\CS$.

Set $H_w:=\cup_j\ H_j(w)$; $P=P_{-1}\cup P_0$; $\tP=\{(t_j)\in\ \BA|
\mbox{ there exists $j$ such that }t_j\in P\}$.

Given $\CK\in\CD(\BA;\ ^{\bz}\CS)$, the restriction map
$$
R\Gamma(\BA,H_w;\CK)\lra R\Gamma(\BA,\tP;\CK)
$$
is an isomorphism by homotopy. On the other hand, we have restriction maps
$$
R\Gamma(\BA,\tP;\CK)\lra R\Gamma(\BA_{\rho},
\tP_{\rho};\CK)
$$
where $\tP_{\rho}:=\tP\cap\BA_{\rho}$.
Therefore we have canonical maps
\begin{equation}
\label{resrho}
r_{\rho}:R\Gamma(\BA,H_w;\CK)\lra
R\Gamma(\BA_{\rho},\tP_{\rho};\CK)
\end{equation}

\subsection{Theorem}
\label{addtheor} {\em For every $\CK\in\CD(\BA;\ ^{\bz}\CS)$ the
canonical map
\begin{equation}
\label{rest}
r=\sum r_{\rho}: R\Gamma(\BA,H_w;\CK)\lra\oplus_{\rho\in\CQ_2(J)}
R\Gamma(\BA_{\rho},\tP_{\rho};\CK)
\end{equation}
is an isomorphism.}

\subsection{Proof.} Let us pick two
open subsets $U_{-1},U_0\subset\BC$ as shown on Fig. 12.
Set $U=U_{-1}\cup U_0$, $\BA_U=\{(t_j)\in\BA|t_j\in U\mbox{ for all }
j\}$. It is clear that
the restriction morphism
$$
R\Gamma(\BA,\tP;\CK)\lra R\Gamma(\BA_U,\tP_U;\CK)
$$
where $\tP_U:=\tP\cap \BA_U$, is an isomorphism.

For each $\rho\in\CQ_2(J)$ set
\begin{equation}
\label{aurho}
\BA_{U,\rho}:=\{(t_j)\in \BA_U|t_j\in U_{\rho(j)}\mbox{ for all }j\};\
\tP_{U,\rho}:=\tP\cap\BA_{U,\rho}
\end{equation}
We have $\BA_U=\bigcup_{\rho}\BA_{U,\rho}$.

\subsubsection{} {\bf Lemma.} {\em For every $\CK\in\CD(\BA;^{\bz}\CS)$
the sum of restriction maps
\begin{equation}
\label{mapq}
q:R\Gamma(\BA_U,\tP_U;\CK)\lra\sum_{\rho\in\CQ_2(J)}
R\Gamma(\BA_{U,\rho},\tP_{U,\rho};\CK)
\end{equation}
is an isomorphism.}

{\bf Proof.}
Suppose we have distinct $\rho_1,\ldots,\rho_m$ such that
$\BA_{U;\rho_1,\ldots,\rho_m}:=\BA_{U,\rho_1}\cap\ldots\cap \BA_{U,\rho_m}
\neq\emp$; set
$\tP_{U;\rho_1,\ldots,\rho_m}:=\tP\cap\BA_{U;\rho_1,\ldots,\rho_m}$.

Our lemma follows at once by Mayer-Vietoris argument from
the following

\subsubsection{} {\bf Claim.} {\em For every $m\geq 2$
$$
R\Gamma(\BA_{U;\rho_1,\ldots,\rho_m},
\tP_{U;\rho_1,\ldots,\rho_m};\CK)=0.
$$}

{\bf Proof of the claim.} It is convenient to use the following notations.
If $J=A\cup B$ is a disjoint union, we will denote by
$\rho_{A,B}$ the map $J\lra [-1,0]$ such that
$\rho^{-1}(-1)=A,\ \rho^{-1}(0)=B$, and
by $U_{A;B}$ the subspace $\BA_{U,\rho_{A,B}}$.

Let us prove the claim for the case $N=2$.
Let $J=\{i,j\}$.
In this case it is easy to see that the only non-trivial intersections are
$U^{(1)}=U_{ij;\emp}\cap U_{j;i}$ and $U^{(2)}=U_{ij;\emp}\cap U_{\emp;ij}$.

To prove our claim for $U^{(1)}$ we will use
{\em a shrinking neighbourhood argument}
based on Lemma I.2.12. Let $\CK'$ denote the sheaf on $U^{(1)}$ obtained
by extension by zero of $\CK|_{U^{(1)}-\tP}$. For each $\epsilon>0$
let us denote by $U_{-1,\epsilon}\subset\BC$ an open domain consisting
of points having distance $<\epsilon$ from $D_{-1}\cup P$. Set
$$
U^{(1)}_{\epsilon}=\{(t_i,t_j)|\ t_i\in U_{-1,\epsilon},\ t_j\in U_0\cap
U_{-1,\epsilon}\}.
$$
It is clear that restriction maps
$R\Gamma(U^{(1)};\CK')\lra R\Gamma(U^{(1)}_{\epsilon};\CK')$ are isomorphisms.
On the other hand, it follows from I.2.12 that
$$
\dirlim_{\epsilon}\ R\Gamma(U^{(1)}_{\epsilon};\CK')\cong
R\Gamma(\bigcap_{\epsilon}\ U^{(1)}_{\epsilon};\CK'),
$$
and the last complex is zero by the definition of $\CK'$
(the point $t_j$ is confined to $P$ in
$\bigcap_{\epsilon}\ U^{(1)}_{\epsilon}$).

The subspace $U^{(2)}$ consists of $(t_i,t_j)$ such that both $t_i$ and $t_j$
lie in $U_{-1}\cap U_{0}$. This case is even simpler. The picture is
homeomorphic to an affine plane with a sheaf smooth along
the diagonal stratification; and we are interested in its
cohomology modulo the coordinate cross. This is clearly equal to zero, i.e.
$R\Gamma(U^{(2)},\tP\cap U^{(2)};\CK)=0$.

This proves the claim
for $N=2$.
The case of an arbitrary $N$ is treated in a similar manner, and we leave
it to the reader. This completes the proof of the claim and of the lemma.
$\Box$

\subsubsection{} {\bf Lemma.} {\em For every $\rho\in\CQ_2(J)$ the restriction
map
$$
R\Gamma(\BA_{U,\rho},\tP_{U,\rho};\CK)\lra
R\Gamma(\BA_{\rho},\tP_{\rho};\CK)
$$
is an isomorphism.}

{\bf Proof.} Again let us consider the case $J=\{i,j\}$.
If $\rho=\rho_{ij;\emp}$ or $\rho_{\emp;ij}$ the statement is obvious.
Suppose $\rho=\rho_{i;j}$. Let us denote by $\BA'_{U,\rho}\subset\BA_{U,\rho}$
the subspace
$\{(t_i,t_j)|t_i\in D_{-1},\ t_j\in U_0\}$. It is clear that
the restriction map
$$
R\Gamma(\BA'_{U,\rho},\tP'_{U,\rho};\CK)\lra
R\Gamma(\BA_{\rho},\tP_{\rho};\CK)
$$
where $\tP'_{U,\rho}:=\tP\cap\BA'_{U,\rho}$,
is an isomorphism. Let us consider the restriction
\begin{equation}
\label{resaa}
R\Gamma(\BA_{U,\rho},\tP_{U,\rho};\CK)\lra
R\Gamma(\BA'_{U,\rho},\tP'_{U,\rho};\CK)
\end{equation}
The cone of this map is isomorphic to
$R\Gamma(\BA_{U,\rho},\tP_{U,\rho};\CM)$ where the sheaf $\CM$ has the same
singularities as $\CK$ and in addition is $0$ over the closure
$\bar{D}_{-1}$. Now, consider a system of shrinking neighbourhoods
of $P_{-1}\cup\bar{D}_{-1}$ as in the proof af the claim above,
we see that $R\Gamma(\BA_{U,\rho},\tP_{U,\rho};\CM)=0$, i.e.
{}~(\ref{resaa})
is an isomorphism. This implies the lemma for this case.

The case of arbitrary $J$ is treated exactly in the same way. $\Box$

Our theorem is an obvious consequence of two previous lemmas. $\Box$

\section{Fusion and tensor products}
\label{fus}

\subsection{Fusion functors} The constructions below were inspired by
{}~\cite{dr}.

For each integer $n\geq 1$ and $i\in [n]$, let
us define functors
\begin{equation}
\label{fusfunct}
^n\psi_i:\CD(^n\BA;\CS_{\Delta})\lra\ \CD(^{n-1}\BA;\CS_{\Delta})
\end{equation}
as follows. We have the $t$-exact nearby cycles functors
(see ~\cite{d} or ~\cite{ks}, 8.6, but note the shift by [-1]!)
$$
\Psi_{t_{-i}-t_{-i+1}}[-1]:\CD(^n\BA;\CS_{\Delta})\lra
\CD(\BA';\CS_{\Delta})
$$
where $\BA'$ denotes (for a moment) an affine space with coordinates
$t_j,\ j\in ((n)\cup J)-\{-i\}$. We can identify the last space with
$^{n-1}\BA$ simply by renaming coordinates $t_j$ to $t_{j+1}$ for
$-n\leq j\leq -i-1$. By definition, $^n\psi_i$ is equal to
$\Psi_{t_{-i}-t_{-i+1}}[-1]$ followed by this identification.

\subsection{Lemma}
\label{assoc} {\em (i) For each $n\geq 2$, $i\in [n]$ have canonical
isomorphisms
\begin{equation}
\label{alphai}
^n\alpha_i:\ ^{n-1}\psi_i\circ\ ^n\psi_i\iso\ ^{n-1}\psi_i\circ\ ^n\psi_{i+1}
\end{equation}
and equalities
$$
^{n-1}\psi_{j}\circ\ ^n\psi_i=\ ^{n-1}\psi_i\circ\ ^n\psi_{j+1}
$$
for $j>i$,
such that

(ii) ({\em "Stasheff pentagon" identity})
the diagram below commutes:

\vspace{.8cm}

\begin{picture}(20,12)(-10,-6)

\put(0,6){$^{n-2}\psi_i\circ\ ^{n-1}\psi_i\circ\ ^n\psi_i$}

\put(5,2){$^{n-2}\psi_i\circ\ ^{n-1}\psi_i\circ\ ^n\psi_{i+1}$}
\put(5,-2){$^{n-2}\psi_i\circ\ ^{n-1}\psi_{i+1}\circ\ ^n\psi_{i+1}$}
\put(0,-6){$^{n-2}\psi_i\circ\ ^{n-1}\psi_{i+1}\circ\ ^n\psi_{i+2}$}

\put(-5,2){$^{n-2}\psi_i\circ\ ^{n-1}\psi_{i+1}\circ\ ^n\psi_i$}
\put(-5,-2){$^{n-2}\psi_i\circ\ ^{n-1}\psi_i\circ\ ^n\psi_{i+2}$}

\put(3.6,5.6){\vector(1,-1){2.6}}
\put(5,4.2){$^{n-2}\psi_i\circ\ ^n\alpha_i$}

\put(7,1){\vector(0,-1){2}}
\put(7.4,0){$^{n-1}\alpha_i\circ\ ^n\psi_{i+1}$}

\put(7,-2.3){\vector(-1,-1){3}}
\put(6,-4){$^{n-2}\psi_i\circ\ ^n\alpha_{i+1}$}

\put(0.4,5.6){\vector(-1,-1){2.4}}
\put(-3.8,4.6){$^{n-1}\alpha_i\circ\ ^n\psi_i$}

\put(-3.1,1){\line(0,-1){2}}
\put(-2.9,1){\line(0,-1){2}}

\put(-3,-2.3){\vector(1,-1){3}}
\put(-4.6,-4){$^{n-1}\alpha_i\circ\
^n\psi_{i+2}$}

\end{picture}
$\Box$}

\subsection{} Let us define a $t$-exact functor
\begin{equation}
\label{psi}
\psi:\CD(^n\BA;\CS_{\Delta})\lra\CD(\BA;\CS)
\end{equation}
as a composition $i_0^*[-1]\circ\ ^n\psi_1\circ\ ^{n-1}\psi_1\circ\ldots
\circ\ ^1\psi_1$, where
$$
i^*_0:\CD(^0\BA;\CS_{\Delta})\lra\CD(\BA;\CS)
$$
denotes the restriction to the subspace $t_0=0$. Note that $i^*_0[-1]$ is
a $t$-exact equivalence. (Recall that $\BA$ and $\CS$ denote the same as in
{}~\ref{j}).

\vspace{.5cm}
{\em STANDARD SHEAVES}
\vspace{.5cm}

The constructions and computations below generalize Section
{}~\ref{standsheaves}.

\subsection{} Let us make the following assumptions. Let us denote by $A$
the $I\times I$-matrix $(i\cdot j)$. Let us suppose that $\det\ A\neq 0$.
There exists a unique $\BZ[\frac{1}{\det\ A}]$-valued symmetric
bilinear form on $X$ (to be denoted by $\lambda,\mu\mapsto\lambda\cdot\mu$)
such that the map $\BZ[I]\lra X,\ \nu\mapsto\lambda_{\nu}$ respects
scalar products.

Let us suppose that our field $B$ contains an element $\zeta'$ such that
$(\zeta')^{\det\ A}=\zeta$, and fix such $\zeta'$. For $a=\frac{c}{\det\ A},\
c\in\BZ$, we set by definition $\zeta^a:=(\zeta')^c$.

\subsection{} Let us fix $\nu=\sum\nu_ii\in \BN[I]$ and its unfolding
$\pi:J\lra I$ as in ~\ref{color}, and an integer $n\geq 1$.
We will use the preceding notations
with this $J$.

Let us fix $n+1$ weights $\Lambda_0,
\Lambda_{-1},\ldots,\Lambda_{-n}\in X$.
We define a one-dimensional local system
$\CI(\Lambda_0,\ldots,\Lambda_{-n};\nu)$ over $^n\BAO$ exactly in the
same manner as in ~\ref{pici}, with half-monodromies equal to
$\zeta^{\pi(i)\cdot\pi(j)}$ if $i,j\in J$,
to $\zeta^{-\langle\Lambda_i,\pi(j)\rangle}$ if $i\in (n),j\in J$ and
to $\zeta^{\Lambda_i\cdot\Lambda_j}$ if $i,j\in (n)$.

Let
$j:\ ^n\BAO\lra\ ^n\BA$ be the embedding. Let us introduce the sheaves
\begin{equation}
\label{ilambdas}
\CI(\Lambda_0,\ldots,\Lambda_{-n};\nu)_?:=
j_?\CI(\Lambda_0,\ldots,\Lambda_{-n};\nu)[-n-N-1]\in\CM(^n\BA;\CS_{\Delta})
\end{equation}
where $?=!,*$ or $!*$.
Applying the functor $\psi$, we get the sheaves
$\psi\CI(\Lambda_0,\ldots,\Lambda_{-n};\nu)_?\in\CM(\BA;\CS)$.

All these objects are naturally $\Sigma_{\pi}$-equivariant. We define
the following sheaves on $\CA_{\nu}$:
\begin{equation}
\label{tensanu}
^{\psi}\CI_{\nu}(\Lambda_0,\ldots,\Lambda_{-n})_?:=
(\pi_*\psi\CI(\Lambda_0,\ldots,\Lambda_{-n};\nu)_?)^{\Sigma_{\pi},-}
\end{equation}
(cf. ~\ref{letus}).

The following theorem generalizes Theorem ~\ref{shaposym}.

\subsection{Theorem}
\label{shsymtens} {\em Given a bijection $\eta:J\iso [N]$, we have
natural isomorphisms included into a commutative square
$$\begin{array}{ccc}
\Phi_{\nu}(^{\psi}\CI_{\nu}(\Lambda_0,\ldots,\Lambda_{-n})_!)&
\overset{\phi_{!}^{(\eta)}}{\iso}&
         (V(\Lambda_0)\otimes\ldots\otimes V(\Lambda_{-n}))_{\nu}\\
m\downarrow&\;&\downarrow S_{\Lambda}\\
\Phi_{\nu}(^{\psi}\CI_{\nu}(\Lambda_0,\ldots,\Lambda_{-n})_*)&
\overset{\phi_{*}^{(\eta)}}{\iso}&
         (V(\Lambda_0)^*\otimes\ldots\otimes V(\Lambda_{-n})^*)_{\nu}
\end{array}$$
and
\begin{equation}
\label{lnucohtens}
\phi_{!*}^{(\eta)}:
\Phi_{\nu}(^{\psi}\CI_{\nu}(\Lambda_0,\ldots,\Lambda_{-n})_{!*})\iso
(L(\Lambda_0)\otimes\ldots\otimes L(\Lambda_{-n}))_{\nu}
\end{equation}
A change of $\eta$ multiplies these isomorphisms by the sign
of the corresponding permutation of $[N]$.}

{\bf Proof.} We may suppose that $\pi$ is injective, i.e. all
$\nu_i=0$ or $1$. The general case immediately follows from this one
after passing to $\Sigma_{\pi}$-skew invariants.

Let us consider the case $n=1$. In this case one sees easily
from the definitions that we have a canonical isomorphism
$$
\Phi_{\nu}(^{\psi}\CI_{\nu}(\Lambda_0,\Lambda_{-1})_?)\cong
R\Gamma(\BA,H_w;\CI(\Lambda_0,\Lambda_{-1};\nu)_?)
$$
in the notations of Additivity Theorem ~\ref{addtheor}. On the other hand,
the set $\CQ_2(J)$ is in one-to-one correspondence with the set of
all decompositions $\nu=\nu_0+\nu_{-1},\ \nu_i\in\BN[I]$, and if $\rho$
corresponds to such
a decomposition, we have a natural isomorphism
\begin{equation}
\label{factorstand}
R\Gamma(\BA_{\rho},\tP\cap\BA_{\rho};\CI_{\nu}(\Lambda_0,\Lambda_{-1})_?)\cong
\Phi_{\nu_0}(^{\psi}\CI(\Lambda_0;\nu_0)_?)\otimes
\Phi_{\nu_{-1}}(^{\psi}\CI(\Lambda_{-1};\nu_{-1})_?)
\end{equation}
by the K\"{u}nneth formula. Therefore,  Additivity Theorem implies
isomorphisms
$$
\Phi_{\nu}(^{\psi}\CI_{\nu}(\Lambda_0,\Lambda_{-1})_?)\cong
\oplus_{\nu_0+\nu_{-1}=\nu}
\Phi_{\nu_0}(^{\psi}\CI_{\nu_0}(\Lambda_0)_?)\otimes
\Phi_{\nu_{-1}}(^{\psi}\CI_{\nu_{-1}}(\Lambda_{-1})_?)
$$
which are the claim of our theorem.

The case $n>2$ is obtained similarly by the iterated use of the Additivity
Theorem. $\Box$

\subsection{} Next we will consider the stalks $(?)_0$ of our sheaves at $0$,
or what is the same (since they are $\BR^{+*}$-homogeneous), the complexes
$R\Gamma(\CA_{\nu};?)$.

The next theorem generalizes Theorem ~\ref{stalks}.

\subsection{Theorem}
\label{stalktens} {\em Given a bijection $\eta:J\iso [N]$, we have
natural isomorphisms included into a commutative square
$$\begin{array}{ccc}
^{\psi}\CI_{\nu}(\Lambda_0,\ldots,\Lambda_{-n})_{!0}&
\overset{\phi_{!,0}^{(\eta)}}{\iso}&\
_{\nu}C^{\bullet}_{\fF}(V(\Lambda_0)\otimes\ldots\otimes V(\Lambda_{-n}))\\
m\downarrow&\;&\downarrow S_{\Lambda}\\
^{\psi}\CI_{\nu}(\Lambda_0,\ldots,\Lambda_{-n})_{*0}&
\overset{\phi_{*,0}^{(\eta)}}{\iso}&\
_{\nu}C^{\bullet}_{\fF^*}(V(\Lambda_0)^*\otimes\ldots\otimes V(\Lambda_{-n})^*)
\end{array}$$
and
\begin{equation}
\label{lnustens}
\phi_{!*,0}^{(\eta)}:\
^{\psi}\CI_{\nu}(\Lambda_0,\ldots,\Lambda_{-n})_{!*0}\iso\
_{\nu}C^{\bullet}_{\ff}(L(\Lambda_0)\otimes\ldots\otimes L(\Lambda_{-n}))
\end{equation}
A change of $\eta$ multiplies these isomorphisms by the sign
of the corresponding permutation of $[N]$.}

{\bf Proof.} It is not hard to deduce from the previous theorem that
we have natural isomorphisms of complexes included into a
commutative square
$$\begin{array}{ccc}
^+C(\BA;\psi\CI(\Lambda_0,\ldots,\Lambda_{-n};\nu)_{!})&\iso&\
_{\chi_J}C^{\bullet}_{^{\pi}\fF}(V(^{\pi}\Lambda_0)\otimes\ldots
\otimes V(^{\pi}\Lambda_{-n}))\\
m\downarrow&\;&\downarrow S_{\Lambda}\\
^+C(\BA;\psi\CI(\Lambda_0,\ldots,\Lambda_{-n};\nu)_{*})&\iso&\
_{\chi_J}C^{\bullet}_{^{\pi}\fF^*}(V(^{\pi}\Lambda_0)^*\otimes\ldots
\otimes V(^{\pi}\Lambda_{-n})^*)
\end{array}$$
This implies our claim after passing to $\Sigma_{\pi}$-(skew) invariants.
$\Box$

\newpage
\begin{center}
{\bf CHAPTER 4. Category $\CC$.}
\end{center}

\section{Simply laced case}
\label{compar}

\subsection{}
\label{zeta} From now on untill the end of the paper we will assume,
in addition to the assumptions
of ~\ref{notations}, that $\zeta$ is
a primitive $l$-th root of unity, where $l$ is a fixed integer $l>3$
prime to $2$ and $3$.

\subsection{}
\label{Lusztig}
We will use notations of ~\cite{l}, Chapters 1, 2, which we briefly recall.

\subsubsection{}
Let $(I,\cdot)$ be {\em a simply laced Cartan datum of finite type}
(cf. {\em loc.cit.}, 1.1.1, 2.1.3),
that is, a finite set $I$ and a nondegenerate symmetric
bilinear form $\alpha,\beta\mapsto\alpha\cdot\beta$ on the free abelian
group $\BZ[I]$. This form satisfies conditions

(a) $i\cdot i=2$ for any $i\in I$;

(b) $i\cdot j\in\{0,-1\}$ for any $i\not=j$ in $I$.

\subsubsection{}
We will consider {\em the simply connected root datum of type $(I,\cdot)$},
that is (see {\em loc. cit.}, 2.2.2),
two free abelian groups $Y=\BZ[I]$ and $X=\Hom_{\BZ}(Y,\BZ)$ together with

(a) the canonical bilinear pairing
$\langle,\rangle:\;Y\times X\lra {\Bbb Z};$

(b) an obvious embedding $I\hookrightarrow Y\;(i\mapsto i)$ and an embedding
$I\hookrightarrow X\;(i\mapsto i')$, such that
$\langle i,j'\rangle=i\cdot j$ for any $i,j\in I$.

We will call $X$ {\em the lattice of weights}, and $Y$
{\em the lattice of coroots}.
An element of $X$ will be typically denoted by $\lambda,\mu,\nu,\ldots$;
and an element of $Y$ will be typically denoted by
$\alpha,\beta,\gamma,\ldots$.

\subsection{}
\label{C}
We consider the finite dimensional algebra $U$ over the field $B$
defined as in the section 1.3 of ~\cite{ajs}. We also consider the
category $\CC$ of finite dimensional $X$-graded $U$-modules
defined as in the section 2.3 of ~\cite{ajs}.

\subsubsection{}
\label{u}
The algebra $U$ is given by generators
$E_i,F_i,K^{\pm1}_i,\ i\in I$, subject to relations

(z) $K_i\cdot K_i^{-1}=1;\ K_iK_j=K_jK_i$;

(a) $K_jE_i=\zeta^{j\cdot i}E_iK_j$;

(b) $K_jF_i=\zeta^{-j\cdot i}F_iK_j$;

(c) $E_iF_j-F_jE_i=\delta_{ij}\frac{K_i-K_i^{-1}}{\zeta-\zeta^{-1}}$;

(d) $E_i^l=F_j^l=0$;

(e) $E_iE_j-E_jE_i=0$ if $i\cdot j=0$;
    $E_i^2E_j-(\zeta+\zeta^{-1})E_iE_jE_i+E_iE_j^2=0$ if $i\cdot j=-1$;

(f) $F_iF_j-F_jF_i=0$ if $i\cdot j=0$;
    $F_i^2F_j-(\zeta+\zeta^{-1})F_iF_jF_i+F_iF_j^2=0$ if $i\cdot j=-1$.

The algebra $U$ has a unique $B$-algebra $X$-grading $U=\oplus U_{\mu}$ for
which $|E_i|=i',\ |F_i|=-i',\ |K_i|=0$.

We define a comultiplication
\begin{equation}
\label{coprodu}
\Delta:U\lra U\otimes U
\end{equation}
as a unique $B$-algebra mapping such that
$$
\Delta(K^{\pm 1})=K^{\pm 1}\otimes K^{\pm 1};\
$$
$$
\Delta(E_i)=E_i\otimes K_i+1\otimes E_i;\
$$
$$
\Delta(F_i)=F_i\otimes 1+K_i^{-1}\otimes F_i.
$$
This makes $U$ a Hopf algebra (with obvious unit and counit).

\subsubsection{} The category $\CC$ is by definition a category
of finite dimensional $X$-graded $B$-vector spaces
$V=\oplus_{\mu\in X}V_{\mu}$,
equipped with a left action of $U$ such that the
$U$-action is compatible with the $X$-grading,
and
$$
K_ix=\zeta^{\langle i,\mu\rangle}x
$$
for $x\in V_{\mu},\ i\in I$.

Since $U$ is a Hopf algebra, $\CC$ has a canonical structure of a {\em tensor
category}.

\subsection{}
\label{fu}
We define an algebra $\fu$ having generators $\theta_i,\epsilon_i,
{K_i}^{\pm1},\ i\in I,$ subject to relations

(z) $K_i\cdot K_i^{-1}=1;\ K_iK_j=K_jK_i$;

(a) ${K_j}\epsilon_i=
{\zeta}^{j\cdot i}\epsilon_i{K_j}$;

(b) ${K_j}\theta_i=
{\zeta}^{-j\cdot i}\theta_i{K_j}$;

(c) $\epsilon_i\theta_j-\zeta^{i\cdot j}\theta_j\epsilon_i=
\delta_{ij}(1-K^{-2}_i)$

(d) if $f\in\Ker (S)\subset\fF$ (see ~(\ref{formap})) then $f=0$;

(e) the same as (d) for the free algebra $\fE$ on the generators $\epsilon_i$.

\subsubsection{}
\label{comultu} Let us define the comultiplication
\begin{equation}
\label{comueq}
\Delta:\fu\lra\fu\otimes\fu
\end{equation}
by the formulas
$$
\Delta({K_i}^{\pm 1})={K_i}^{\pm 1}\otimes{K_i}^{\pm 1};\
$$
$$
\Delta(\theta_i)=\theta_i\otimes 1+K_i^{-1}\otimes\theta_i;\
$$
$$
\Delta(\epsilon_i)=\epsilon_i\otimes 1+K_i^{-1}\otimes\epsilon_i
$$
and the condition that $\Delta$ is a morphism of $B$-algebras.

This makes $\fu$ a Hopf algebra (with obvious unit and counit).

$\fu$ is an $X$-graded $B$-algebra, with an $X$-grading defined uniquely
by the conditions $|{K_i}^{\pm 1}|=0;\
|\theta_i|=-i';\ |\epsilon_i|=i'$.

\subsection{} We define $\tCC$ as a category of finite dimensional
$X$-graded vector spaces $V=\oplus V_{\lambda}$, equipped with a structure
of a left
$\fu$-module compatible with $X$-gradings and such that
$$
{K_i}x={\zeta}^{\langle i,\lambda\rangle}x
$$
for $x\in V_{\lambda},\ i\in I$.

Since $\fu$ is a Hopf algebra, $\tCC$ is a tensor category.

\subsection{}
Recall that for any $\Lambda\in X$ we have defined in ~\ref{llambda}
the $X$-graded $\ff$-module $L(\Lambda)$. It is a quotient-module of the
Verma module $V(\Lambda)$,
and it inherits its $X$-grading from the one of $V(\Lambda)$
(see ~\ref{verma}).
Thus $L(\Lambda)=\oplus L(\Lambda)_\lambda$, and we define the action of
generators ${K_i}$ on $L(\Lambda)_\lambda$ as multiplication by
${\zeta}^{\langle i,\lambda\rangle}$. Finally, we define the action of
generators $\epsilon_i$ on $V(\Lambda)$ as in ~\ref{slambda}. These operators
on $V(\Lambda)$ satisfy the relations (a) --- (c) above.

We
check immediately that this action descends to the quotient $L(\Lambda)$.
Moreover, it follows from Theorem ~\ref{coactshap}
that these operators acting on $L(\Lambda)$ satisfy the
relations (a) --- (e) above. So we have constructed the action of $\fu$
on $L(\Lambda)$, therefore we can regard it as an object of $\tCC$.

\subsubsection{} {\bf Lemma.} {\em $L(\Lambda)$ is an
irreducible object in $\tCC$.}

{\bf Proof.}
Let $I(\Lambda)$ be the maximal proper homogeneous (with respect to
$X$-grading)
submodule of $V(\Lambda)$ (the sum of all homogeneous submodules not containing
$v_\Lambda$). Then $V(\Lambda)/I(\Lambda)$ is irreducible, so it suffices to
prove that $I(\Lambda)=\ker(S_\Lambda)$. The inclusion $\ker(S_\Lambda)
\subset I(\Lambda)$ is obvious. Let us prove the opposite inclusion. Let
$y\in I(\Lambda)$. It is enough to check that $S_\Lambda(y,x)=0$ for any
$x\in V(\Lambda)$ of the form $\theta_{i_1}\ldots\theta_{i_n}v_\Lambda$.
By (22) we have $S_\Lambda(y,\theta_{i_1}\ldots\theta_{i_n}v_\Lambda)=
S_\Lambda(\epsilon_{i_n}\ldots\epsilon_{i_1}y,v_\Lambda)=0$ since
$\epsilon_{i_n}\ldots\epsilon_{i_1}y\in I(\Lambda)$. $\Box$

\subsection{}
Let us consider elements $E_i,F_i\in\fu$ given by the following formulas:

\begin{equation}
\label{formr}
E_i=\frac{\zeta^2}{\zeta-\zeta^{-1}}
\epsilon_i{K_i};\
F_i=\theta_i
\end{equation}
It is immediate to check that these elements satisfy the relations
{}~\ref{u} (a) --- (c).

Moreover, one checks without difficulty that
$$
\theta_i\theta_j-\theta_j\theta_i\in\Ker(S)\ \mbox{if}\ i\cdot j=0,
$$
and
$$
\theta_i^2\theta_j-(\zeta+\zeta^{-1})\theta_i\theta_j\theta_i+
\theta_j\theta_i^2\in\Ker(S)\ \mbox{if}\ i\cdot j=-1
$$
(cf. ~\cite{sv2}, 1.16).
Also, it is immediate that
$$
S(\theta_i^a,\theta_i^a)=\prod_{p=1}^a\frac{1-\zeta^{2p}}{1-\zeta^2}.
$$
It follows that $\theta_i^l\in\Ker(S)$ for all $i$.

It follows easily that the formulas ~(\ref{formr}) define a surjective
morphism of algebras
\begin{equation}
\label{mapr}
R:U\lra\fu
\end{equation}
Moreover, one checks at once that $R$ is a map of Hopf algebras.

Therefore, $R$ induces a tensor functor
\begin{equation}
\label{Q}
Q:\tCC\lra\CC
\end{equation}
which is an embedding of a full subcategory.

\subsection{Theorem}
\label{ceqcc}
{\em $Q$ is an equivalence.}

{\bf Proof.} It is enough to check that $\tilde{\CC}$ contains enough
projectives for $\CC$ (see e.g. Lemma A.15. of ~\cite{kl} IV).

First of all, we know from ~\cite{ajs}, section 4.1,
that the simple $\fu$-modules
$L(\Lambda),\Lambda\in X$ exhaust the list of simple objects of $\CC$.
Second, we know from ~\cite{apk}, Theorem 4.6 and Remark 4.7,
that the module $L(-\rho)$ is projective
where $-\rho\in X$ is characterized by the property $\langle i,-\rho\rangle=-1$
for any $i\in I$.
Finally we know, say, from ~\cite{apk}, Remark 4.7, and ~\cite{apw}, Lemma
9.11,
that the set of
modules $\{L(\Lambda)\otimes L(-\rho), \Lambda\in X\}$ is an ample system of
projectives for $\CC$.
$\Box$

\subsection{}
Denote by $\fu^0$ (resp., $U^0$) the subalgebra of $\fu$ (resp., of $U$)
generated by $K_i^{\pm 1},\ i\in I$. Obviously, both algebras
are isomorphic to the ring of Laurent polynomials in $K_i$, and the map
$R$ induces an identity isomorphism between them.

Denote by $\fu^{\leq 0}$ (resp., by $\fu^-$) the subalgebra of
$\fu$ generated by $\theta_i, {K_i}^{\pm 1}$
(resp., by $\theta_i$), $i\in I$. The last algebra may be
identified with $\ff$.
As a vector space, $\fu^{\leq 0}$
is isomorphic to $\ff\otimes\fu^0$.

Denote by $U^-\subset\ U$ the subalgebra generated by $F_i,i\in I$,
and by $U^{\leq 0}\subset\ U$ the subalgebra generated by $F_i,\
{K_i}^{\pm 1},i\in I$. As a vector space it is isomorphic to
$U^-\otimes U^0$.

\subsection{Theorem.} {\em (a) $R$ is an isomorphism;

(b) $R$ induces an isomorphism $U^-\iso\ff$.}

{\bf Proof.} Evidently it is enough to prove b). We know that $R$
is surjective, and
that $U^-$ is finite dimensional. So it suffices to prove that
$\dim\ U^-\leq\dim\ \ff$. We know by ~\cite{l} 36.1.5. that
$\dim\ U^-=\dim L(-\rho)$. On the other hand, the map $\ff\lra L(-\rho),
f\mapsto f(v_{-\rho})$ is surjective by construction. $\Box$

\section{Non-simply laced case}

In the non-simply laced case all main results of this paper hold true as well.
However, the definitions need some minor modifications.
In this section we will describe them.

\subsection{} We will use terminology and notations from ~\cite{l},
especially from Chapters  1-3.
Let us fix a Cartan datum $(I,\cdot)$ of finite type, not necessarily
simply laced, cf. {\em loc. cit.}, 2.1.3.
Let $(Y=\BZ[I],X=\Hom (Y,\BZ),\langle,\rangle,I\hra X,I\hra Y)$
be the simply connected root datum associated with $(I,\cdot)$,
{\em loc.cit.}, 2.2.2.

We set $d_i:=\frac{i\cdot i}{2},\ i\in I$; these numbers are positive
integers. We set $\zeta_i:=\zeta^{d_i}$.

The embedding $I\hra X$ sends $i\in I$ to $i'$ such that
$\langle d_jj,i'\rangle =j\cdot i$ for all $i,j\in I$.

\subsection{} The category $\CC$ is defined in the same way as in the simply
laced case, where the definition of the Hopf algebra $U$ should be modified as
follows (cf. ~\cite{ajs}, 1.3).

By definition, $U$ has generators
$E_i,F_i,K^{\pm1}_i,\ i\in I$, subject to relations

(z) $K_i\cdot K_i^{-1}=1;\ K_iK_j=K_jK_i$;

(a) $K_jE_i=\zeta^{\langle j, i'\rangle}E_iK_j$;

(b) $K_jF_i=\zeta^{-\langle j, i'\rangle}F_iK_j$;

(c) $E_iF_j-F_jE_i=\delta_{ij}\frac{\tK_i-\tK_i^{-1}}{\zeta_i-\zeta_i^{-1}}$;

(d) $E_i^l=F_j^l=0$;

(e) $\sum_{p=0}^{1-\langle i,j'\rangle}(-1)^pE_i^{(p)}E_j
E_i^{(1-\langle i,j'\rangle-p)}=0$ for $i\neq j$;

(f) $\sum_{p=0}^{1-\langle i,j'\rangle}(-1)^pF_i^{(p)}F_j
F_i^{(1-\langle i,j'\rangle-p)}=0$ for $i\neq j$,

where we have used the notations:
$\tK_i:=K_i^{d_i}$;
$G_i^{(p)}:=G_i^p/[p]_i^!,\ G=E$ or $F$,
$$
[p]_i^!:=\prod_{a=1}^p\frac{\zeta_i^a-\zeta_i^{-a}}{\zeta_i-\zeta^{-1}_i}.
$$
The $X$-grading on $U$ is defined in the same way as in the simply laced case.

The comultiplication is defined as
$$
\Delta(K_i)=K_i\otimes K_i;
$$
$$
\Delta(E_i)=E_i\otimes\tK_i+1\otimes E_i;
$$
$$
\Delta(F_i)=F_i\otimes 1+\tK_i^{-1}\otimes F_i.
$$

\subsubsection{Remark} This algebra is very close to
(and presumably coincides with)
the algebra
$\bU$ from ~\cite{l}, 3.1, specialized to $v=\zeta$. We use the opposite
comultiplication, though.

\subsection{} The definition of the Hopf algebra $\fu$ should be modified as
follows. It has generators
$\theta_i,\epsilon_i,
{K_i}^{\pm1},\ i\in I,$ subject to relations

(z) $K_i\cdot K_i^{-1}=1;\ K_iK_j=K_jK_i$;

(a) ${K_j}\epsilon_i=
{\zeta}^{\langle j, i'\rangle}\epsilon_i{K_j}$;

(b) ${K_j}\theta_i=
{\zeta}^{-\langle j, i'\rangle}\theta_i{K_j}$;

(c) $\epsilon_i\theta_j-\zeta^{i\cdot j}\theta_j\epsilon_i=
\delta_{ij}(1-\tK^{-2}_i)$

(d) if $f\in\Ker (S)\subset\fF$ (see ~(\ref{formap})) then $f=0$;

(e) the same as (d) for the free algebra $\fE$ on the generators $\epsilon_i$.

The comultiplication is defined as
$$
\Delta({K_i}^{\pm 1})={K_i}^{\pm 1}\otimes{K_i}^{\pm 1};\
$$
$$
\Delta(\theta_i)=\theta_i\otimes 1+\tK_i^{-1}\otimes\theta_i;\
$$
$$
\Delta(\epsilon_i)=\epsilon_i\otimes 1+\tK_i^{-1}\otimes\epsilon_i
$$

The category $\tCC$ is defined in the same way as in the simply laced
case.

\subsection{} We define an $X$-grading on the free $\fF$-module
$V(\Lambda)$ with generator $v_{\Lambda}$ by setting
$V(\Lambda)_{\Lambda}=B\cdot v_{\Lambda}$ and assuming that
operators $\theta_i$ decrease the grading by $i'$.

The definition of the form $S_{\Lambda}$ on $V(\Lambda)$ should be
modified as follows.
It is a unique bilinear from such that
$S_{\Lambda}(v_{\Lambda},v_{\Lambda})=1$ and
$S(\theta_ix,y)=S(x,\epsilon_iy)$ where the operators
$\epsilon_i:V(\Lambda)\lra V(\Lambda)$ are defined by the
requirements $\epsilon_i(v_{\Lambda})=0$,
$$
\epsilon_i(\theta_jx)=\zeta^{i\cdot j}\theta_j\epsilon_i(x)+
\delta_{ij}[\langle i,\lambda\rangle]_{\zeta_i}x
$$
for $x\in V(\Lambda)_{\lambda}$.

We define $L(\Lambda)$ as a quotient $V(\Lambda)/\Ker(S)$. As in the simply
laced case, $L(\Lambda)$ is naturally an object of $\CC$, and
the same argument proves that it is irreducible.

\subsection{} We define the morphism
\begin{equation}
\label{maprnsl}
R: U\lra\fu
\end{equation}
by the formulas
\begin{equation}
\label{formrnsl}
R(E_i)=\frac{\zeta_i^2}{\zeta_i-\zeta_i^{-1}}
\epsilon_i{\tK_i};\
R(F_i)=\theta_i;\ R(K_i)=K_i
\end{equation}
Using ~\cite{l}, 1.4.3, one sees immediately that it is correctly
defined morphism of algebras. It follows at once from the definitions
that $R$ is a morphism of Hopf algebras.

Hence, we get a tensor functor
\begin{equation}
\label{Qnsl}
Q:\tCC\lra\CC
\end{equation}
and the same proof as in ~\ref{ceqcc} shows that $Q$ is an {\em equivalence
of categories}.

It is a result of primary importance for us.
It implies in particular that all irreducibles in $\CC$
(as well as their tensor products), come from $\tCC$.

\subsection{} Suppose we are given $\Lambda_0,\ldots,\Lambda_{-n}\in X$
and $\nu\in \BN[I]$. Let $\pi:J\lra I$, $\pi:\ ^{\pi}\BA\lra\CA_{\nu}$
denote the same as in ~\ref{color}. We will use the notations
for spaces and functors from Section ~\ref{fus}.

The definition of the local system $\CI(\Lambda_0,\ldots,\Lambda_{-n};
\nu)$ from {\em loc. cit.} should be modified: it should
have half-monodromies
$\zeta_j^{-\langle \pi(j),\Lambda_i\rangle}$ for
$i\in (n),\ j\in J$, the other formulas stay without change.

After that, the standard sheaves are defined as in {\em loc. cit.}
Now we are arriving at the main results of this paper.
The proof is  the same as the proof of theorems ~\ref{shsymtens}
and ~\ref{stalktens}, taking into account the previous
algebraic remarks.

\subsection{Theorem}
\label{shsymnsl} {\em Let $L(\Lambda_0),\ldots,L(\Lambda_{-n})$ be
irreducibles of $\CC$, $\lambda\in X$,
$\lambda=\sum_{m=0}^n\Lambda_{-m}\ -\ \sum_i\nu_ii'$ for some
$\nu_i\in\BN$.  Set $\nu=\sum\nu_ii\in\BN[I]$.

Given a bijection $\eta:J\iso [N]$, we have
natural isomorphisms
\begin{equation}
\label{lnucohnsl}
\phi_{!*}^{(\eta)}:
\Phi_{\nu}(^{\psi}\CI_{\nu}(\Lambda_0,\ldots,\Lambda_{-n})_{!*})\iso
(L(\Lambda_0)\otimes\ldots\otimes L(\Lambda_{-n}))_{\lambda}
\end{equation}
A change of $\eta$ multiplies these isomorphisms by the sign
of the corresponding permutation of $[N]$. $\Box$}

\subsection{Theorem}
\label{stalknsl} {\em In the notations of the previous theorem we have
natural isomorphisms
\begin{equation}
\label{lnusnsl}
\phi_{!*,0}^{(\eta)}:\
^{\psi}\CI_{\nu}(\Lambda_0,\ldots,\Lambda_{-n})_{!*0}\iso\
C^{\bullet}_{\ff}(L(\Lambda_0)\otimes\ldots\otimes L(\Lambda_{-n}))_{\lambda}
\end{equation}
where we used the notation $C^{\bullet}(\ldots)_{\lambda}$ for
$_{\nu}C^{\bullet}(\ldots)$.
A change of $\eta$ multiplies these isomorphisms by the sign
of the corresponding permutation of $[N]$. $\Box$}




\begin{thebibliography}{MMMMM}


\bibitem[AJS]{ajs} H.Andersen, J.Jantzen, W.Soergel, Representations of
quantum groups at $p$-th root of unity and of semisimple groups in
characteristic $p$: independence of $p$, {\em Ast\'{e}risque} {\bf 220}(1994).

\bibitem[APK1]{apw} H.Andersen, P.Polo, W.Kexin, Representations of quantum
algebras, {\em Invent. math.} {\bf 104} (1991), 1-59.

\bibitem[APK2]{apk} H.Andersen, P.Polo, W.Kexin, Injective modules for
quantum algebras, {\em Amer. J. Math.} {\bf 114}(1992), 571-604.

\bibitem[BBD]{bbd} A.Beilinson, J.Bernstein, P.Deligne, Faisceaux Pervers,
{\em Ast\'{e}risque} {\bf 100}(1982).

\bibitem[D]{d} P.Deligne, Le formalisme des cycles \'{e}vanescents,
in: "Groupes de monodromie en G\'{e}om\'{e}trie Alg\'{e}brique (SGA 7)",
Lect. Notes Math. {\bf 340}, Springer-Verlag, Berlin et al., 1973,
82-115.

\bibitem[Dr]{dr} V.Drinfeld, Quasihopf algebras, {\em Algebra and Analysis},
{\bf 1}, no. 6(1989), 114-149(russian).

\bibitem[FS]{fs} M.Finkelberg, V.Schechtman, Localization of $\fu$-modules. I.
Intersection cohomology of real arrangements, Preprint hep-th/9411050
(1994), 1-23.

\bibitem[FW]{fw} G.Felder, C.Wieczerkowski, Topological representations
of the quantum group $U_q(sl_2)$, {\em Comm. Math. Phys.,}
{\bf 138}(1991), 583-605.

\bibitem[KL]{kl} D.Kazhdan, G.Lusztig, Tensor structures arising from
affine Lie algebras. I-IV, {\em Amer. J. Math.}, {\bf 6}(1993), 905-947;
{\bf 6}(1993), 949-1011; {\bf 7}(1994), 335-381; {\bf 7}(1994), 383-453.

\bibitem[KS]{ks} M.Kashiwara, P.Schapira, Sheaves on Manifolds,
{\em Grund. math. Wiss.} {\bf 292}, Springer-Verlag, Berlin et al., 2nd
Printing, 1994.

\bibitem[L]{l} G.Lusztig, Introduction to quantum groups, Birkh\"{a}user,
Boston et al., 1993.

\bibitem[M]{m} S.MacLane, Homology, Springer-Verlag, Berlin et al., 1963.

\bibitem[S]{s} V.Schechtman, Vanishing cycles and quantum groups II,
{\em Int. Math. Res. Notes}, {\bf 10}(1992), 207-215.

\bibitem[SV1]{sv1} V.Schechtman, A.Varchenko, Arrangements of hyperplanes
and Lie algebra homology, {\em Inv. Math.}, {\bf 106}(1991), 139-194.

\bibitem[SV2]{sv2} V.Schechtman, A.Varchenko, Quantum groups and homology
of local systems, in: Algebraic Geometry and Analytic Geometry,
A. Fujiki et al. (eds.), Springer-Verlag, Tokyo et al., 1991, 182-197.

\bibitem[V]{v} A.Varchenko, Multidimensional hypergeometric functions
and representation theory of  Lie algebras
and quantum groups,  Advanced Series in Mathematical Physics
- Vol. 21, World Scientific Publishers, to appear



\end{thebibliography}
\end{document}